\documentclass[11pt,a4paper]{article}
\usepackage{ifthen}
\usepackage{amsmath}
\usepackage{amssymb}
\usepackage{latexsym}
\usepackage{bbm}
\usepackage{graphicx}
\usepackage{pdfsync}
\usepackage{verbatim}
\usepackage{fancyhdr}
\usepackage{url}               
\usepackage[english]{babel}
\usepackage[applemac]{inputenc}

\usepackage{dsfont}

\usepackage[pdftex,%
a4paper=true,%
bookmarks=true,%
bookmarksnumbered=true,%
colorlinks=true,%
urlcolor=blue,%
pdfstartview=FitH%
]{hyperref}

\hypersetup{
	colorlinks=true,     
	citecolor=red,
	linkcolor=blue,     
	}

\usepackage{tikz}

\newtheorem{theorem}{Theorem}[section]
\newtheorem{lemma}[theorem]{Lemma}

\usepackage{enumitem}
\newlist{assumptions}{enumerate}{1}
\setlist[assumptions,1]{label=(A\arabic*), ref=(A\arabic*)}
\newlist{assb}{enumerate}{1}
\setlist[assb,1]{label=(B\arabic*), ref=(B\arabic*)}
\newtheorem{proposition}[theorem]{Proposition}
\newtheorem{corollary}[theorem]{Corollary}
\newtheorem{predefinition}[theorem]{Definition}
\newenvironment{definition}{\begin{predefinition}\upshape}{\end{predefinition}}
\newtheorem{preremark}[theorem]{Remark}
\newenvironment{remark}{\begin{preremark}\upshape}{\end{preremark}}

\numberwithin{equation}{section}
\numberwithin{figure}{section}

\newtheorem{Complementary lemma}[theorem]{Complementary lemma}

\newcommand{\bp}{{\it Proof. }}
\newcommand{\ep}{\hfill $\square$\\}

\newcommand{\bpl}{{\it Proof of Lemma }} 
\newcommand{\epl}{\hfill $\square$\\}

\newcommand{\bpp}{{\it Proof of Proposition }} 
\newcommand{\epp}{\hfill $\square$\\}

\newcommand{\bpt}{{\it Proof of Theorem }}
\newcommand{\ept}{\hfill $\square$\\}

\newcommand{\Imm}{\mbox{Im}}
\renewcommand{\Re}{\mbox{Re}}
\newcommand{\be}{\begin{equation}}
\newcommand{\ee}{\end{equation}}
\newcommand{\bea}{\begin{eqnarray}}
\newcommand{\eea}{\end{eqnarray}}
\newcommand{\bee}{\begin{eqnarray*}}
\newcommand{\eee}{\end{eqnarray*}}

\def\pa{\partial}

\def\na{\nabla}

\def\CC{\mathbb{C}}

\def\RR{\mathbb{R}}

\def\dps{\displaystyle}
\def\ni{\noindent}

\def\eps{\vare}

\def\eps{\varepsilon}

\catcode`@=11
\def\supess{\mathop{\operator@font Sup\,ess}}
\catcode`@=12
\def\CC{\mathbb{C}}

\def\RR{\mathbb{R}}
\def\CC{\mathbb{C}}

\def\ds{\displaystyle}

\def\ni{\noindent}

\def\bar#1{{\overline #1}}

\def\R2+{\RR ^2_+}

\def\pa{\partial}
\def\na{\nabla}

\def\lim{\mathop{\rm lim}}

\def\sup{\mathop{\rm sup}}

\def\log{{\rm log}}

\def\pa{\partial}

\def\pa{\partial}

\newcommand{\ud}{\mathrm{d}}
\newcommand{\ue}{\mathrm{e}}
\newcommand{\ui}{\mathrm{i}}

\begin{document}

\renewcommand{\refname}{References}
\bibliographystyle{abbrv}

\pagestyle{fancy}
\fancyhead[L]{ }
\fancyhead[R]{}
\fancyfoot[C]{}
\fancyfoot[L]{ }
\fancyfoot[R]{}
\renewcommand{\headrulewidth}{0pt} 
\renewcommand{\footrulewidth}{0pt}


\newcommand{\montitre}{On the spectral problem and fractional diffusion limit for Fokker-Planck with(-out) drift and for a general heavy tail equilibrium}

\newcommand{\auteur}{\textsc{Dahmane Dechicha}}
\newcommand{\affiliation}{Laboratoire de recherche AGM, CY Cergy Paris Universit\'e. UMR CNRS 8088 \\

2 Avenue Adolphe Chauvin 95302 Cergy-Pontoise Cedex, France \\
\url{dahmane.dechicha@cyu.fr} }

 \begin{center}
{\bf  {\LARGE \montitre}}\\ \bigskip \bigskip
 {\large\auteur}\\ \bigskip \smallskip
 \affiliation \\ \bigskip
\today
 \end{center}
 \begin{abstract} 
This paper is devoted to the study of a kinetic Fokker-Planck equation with general heavy-tailed equilibrium without an explicit formula, such as $C_\beta \langle v \rangle^{-\beta}$,  in particular non-symmetric and non-centred. This work extends the results obtained in \cite{DP} and \cite{DP-2}. We prove that if the equilibrium behaves like $\langle v \rangle^{-\beta}$ at infinity with $\beta > d$, along with an other assumption, there exists a unique eigenpair solution to the spectral problem associated with the Fokker-Planck operator, taking into account the advection term. As a direct consequence of this construction, and under the hypothesis of the convergence of the rescaled equilibrium, we obtain the fractional diffusion limit for the kinetic Fokker-Planck equation, with or without drift, depending on the decay of the equilibrium and whether or not the first moment is finite. This latter result generalizes all previous results on the fractional diffusion limit for the Fokker-Planck equation and rigorously justifies the remarks and results mentioned in \cite[Section 9]{BM}.
\end{abstract}

\tableofcontents

\newpage

 \pagestyle{fancy}
\fancyhead[R]{\thepage}
\fancyfoot[C]{}
\fancyfoot[L]{}
\fancyfoot[R]{}
\renewcommand{\headrulewidth}{0.2pt} 
\renewcommand{\footrulewidth}{0pt} 

\section{Introduction and main results}
\subsection{Setting of the problem}
In this paper we study the kinetic Fokker-Planck equation with a genral heavy tail equilibrium. This equation decomposes into two parts: a transport part, which represents the movement of particles in straight lines with a velocity $ v $, and a part given by a linear collisional operator, which acts as an external force. Thus, the equation is given by
\begin{equation}\label{FP} 
\left\{\begin{array}{l}
 \partial_t f + v\cdot \nabla_x f = \mathsf{Q} (f), \; \; \; t \geqslant 0 , \ x \in \mathbb{R}^d , \ v \in \mathbb{R}^d ,\\
\\
    f(0,x,v) = f_0(x,v), \hspace{1.8cm} x \in \mathbb{R}^d , \ v \in \mathbb{R}^d ,
\end{array}\right.
\end{equation}
where the collisional Fokker-Planck operator $ \mathsf{Q} $ is given by
\begin{equation}\label{defQ} 
\mathsf{Q}(f)=\nabla_v\cdot\bigg(F\nabla_v\bigg(\frac{f}{F}\bigg)\bigg) = \Delta_v f - \mathrm{div}_v \bigg(\frac{\na_v F}{F} f\bigg) ,
\end{equation}
and $ F $ is the equilibrium of $ \mathsf{Q}$, a fixed function that depends only on $ v $ and satisfies  
$$ \mathsf{Q}(F)=0 \qquad \text{and} \qquad \int_{\mathbb{R}^d} F(v)\ \mathrm{d}v = 1 . $$ 
For an initial condition $ f_0 \geqslant 0 $, the unknown $ f(t,x,v) \geqslant 0 $ in \eqref{FP} can be interpreted as the density of a collection of particles at time $ t \geqslant 0 $, at position $ x \in \mathbb{R}^d $ and with velocity $v \in \mathbb{R}^d$. The Fokker-Planck equation \eqref{FP}-\eqref{defQ} deterministically describes the Brownian motion of this collection of particles, hence the presence of the Laplacian.
The assumption $ F \in L^1 $ reflects the interest in finite mass distributions, i.e., non-negative functions in $L^1$. In this case, the kernel of the differential operator $\mathsf{Q}$ is one-dimensional and is generated by the equilibrium $ F $.

In stochastic language, one can study \eqref{FP}-\eqref{defQ} through the corresponding trajectories
$$ X(t) = x + \int_0^t V(s) \, \ud s \qquad \text{and} \qquad V(t) = v + 2 B_t + \int_0^t \frac{\nabla_v F (V(s))}{F(V(s))} \, \ud s , $$
with $(B_t)_{t \geqslant 0} $ being the $d$-dimensional Brownian motion. Then, for $(V(t),X(t))$ solving the previous system, the family of time-marginals $f(t) = \text{Law}(X(t),V(t))$ solves \eqref{FP}-\eqref{defQ} in the sense of distributions. We refer to \cite{FT} for the probabilistic approach to the Fokker-Planck equation. 

Among the objectives of this paper is the study of the diffusion limit for the kinetic equation \eqref{FP}-\eqref{defQ}. The interest in the diffusion approximation lies in its ability to simplify kinetic models with generally complex collision operators, as well as in the typically large size of the phase space. Theoretically, this limit falls within the scope of the hydrodynamic limit, as it involves transitioning from a kinetic equation, written at the mesoscopic scale, to the diffusion equation, written at the macroscopic scale. Since the time and space scales at the macroscopic level are much larger than those at the mesoscopic one,  the mean free path parameter $\eps \ll 1$ is introduced to transition between these regimes. This parameter represents the characteristic distance between two collisions. The time and space variables are rescaled by setting $t'=\theta(\eps) t$ and $x'= \eps x$, where $\theta(\eps)$ tends to $0$ as $\eps \to 0$, to be determined later. Thus, the solution of \eqref{FP} satisfies, after rescaling and omitting primes, the following equation:
$$ \theta(\varepsilon) \partial_t f^\varepsilon + \eps v \cdot \nabla_x f^\varepsilon = \mathsf{Q}(f^\varepsilon). $$ 
In this work, we focus on heavy-tailed equilibria, which do not decay faster than a polynomial, more generally satisfying
$$ \langle v \rangle^{-\beta} \lesssim F(v) \lesssim \langle v \rangle^{-\beta} ,$$
with $\beta > d$.  The symmetry of the equilibrium plays a crucial role in studying the limit of the solution $f^\eps$ as $\eps \to 0$. For a non-symmetric equilibrium, this limit may not exist. This leads us to consider instead the following rescaled equation with a drift term, given by the first moment of the equilibrium, instead of the previous equation:
\begin{equation}\label{fp-theta}
\theta(\varepsilon) \partial_t f^\varepsilon + \varepsilon (v - j^\eps)\cdot \nabla_x f^\varepsilon = \mathsf{Q}(f^\varepsilon) , 
\end{equation}
where
\begin{equation}\label{j_eps}
j^\eps := \left\{ \begin{array}{l}
0  \hspace{3.75cm} \text{ if } \quad \beta < d+1 ,   \\
\\
 \ds \int_{|v|\leqslant \eps^{-\frac{1}{3}}} v F(v) \ud v \hspace{1.12cm} \text{ if } \quad \beta = d+1 ,    \\
 \\
\ds  \int_{\RR^d} v F(v) \ud v  \hspace{1.9cm} \text{ if } \quad \beta > d+1 .
\end{array}\right.  
\end{equation}
Note that the initial condition written in non-rescaled variables is a well-prepared condition.

In the case of non-symmetrical equilibrium, the integral $\int v F \, \ud v$ is not always finite, as it depends on $\beta$. This is why $j^\eps$ is defined by several expressions. The symmetric case corresponds to $j^\eps = 0$ for all $\beta > d$. The value $\beta = d+1$ is the threshold where the integral $\int v F \, \ud v$ becomes infinite. The truncation at $\eps^{-\frac{1}{3}}$ corresponds to the velocities where all terms of the Fokker-Planck operator, including the advection term, are of the same order. In this case, $j^\eps$ is of the order of $|\ln(\eps)|$, more precisely, $j^\eps \sim c \ln\big(\eps^{-\frac{1}{3}}\big)$ for a certain $c \in \RR^d$ given in Lemma \ref{j_eps equiv ln eps}.

The study of the diffusion limit for kinetic models such as linear Boltzmann (LB), Fokker-Planck (FP), and L\'evy-Fokker-Planck (LFP) has been the subject of several papers since the end of the last century. This limit depends on the properties of the operator $\mathsf{Q}$, the nature of the equilibrium $F$, and the scaling considered. Broadly, two classes of equilibria are distinguished: Gaussian equilibria and heavy-tailed equilibria. Results for the first class of equilibrium, which exhibit rapid decay, show that a classical diffusion equation can be derived with the classical time scale $\theta(\eps) = \eps^2$. Its diffusion coefficient involves an integral of the velocity moment of the equilibrium, of order $4$ in the case of FP and an order depending on the decay of the cross section in the case of LB. See, for instance, \cite{BaSaSe,BeLi,DeGoPo,LaKe,De2} for the LB case, and \cite{DeMa-Ga} for FP. For the case of slow decay equilibrium, it depends on the decay rate since the integral giving the diffusion coefficient (in the previous case) can be infinite, necessitating a change in $\theta(\eps)$ as $\eps^2$ is no longer appropriate. However, as long as the coefficient remains well-defined, the classical diffusion is obtained until the critical case, where the coefficient is defined as a principal value with the anomalous time scale $\theta(\eps) = \eps^2|\ln(\eps)|$. Refer to \cite{MMM} for the LB equation and \cite{NaPu,CNP} for the FP equation.

When $F$ decays very slowly and does not have sufficiently bounded moments, the study is more complicated. This study has been the interest of many papers in the last two decades, with different methods and for different collision operators. A fractional diffusion limit has been obtained in the case of LB when the cross section is such that the operator has a spectral gap; see \cite{MMM} for the pioneering paper in the case of space-independent cross section,  and see \cite{M} for a weak convergence result obtained by the moment method, which also applies to cross sections that depend on the position variable with bounded collision frequency. The case of a more general cross-section, which depends on the spatial variable and degenerates for large velocities, is given in \cite[Chapter 1]{dechicha2023fractional}. See also \cite{BMP-FD} for a strong convergence obtained via a Hilbert expansion, and \cite{MKO} for a probabilistic approach.  Further results on fractional-diffusion-advection limits for kinetic models have been obtained, see \cite{aceves2016fractional,aceves2019fractional} and the references therein.  All the previous results were in the case $x \in \RR^d$. The study on bounded domains is significantly more difficult, as it depends on the definition of the operator on the boundary, which is related to the trajectory of the particles. This generates issues with trace problems and the uniqueness of solutions for boundary value problems. For an analysis methods, see for example \cite{cesbron2020fractional, cesbron2022fractional, cesbron2024fractional}, and for a probabilistic approach, see \cite{bethencourt2022fractional}.

In the case of FP, despite similarities with LB, even if the same results are obtained, the methods used are sometimes entirely different, and the difficulties are not the same. The diffusion limit for the FP equation seems more complicated than the LB one, and the main difficulty is due to the absence of a spectral gap and the fact that all terms contribute to the limit, including the advection term. This can be seen for velocities close to $\eps^{-\frac{1}{3}}$, where all terms are of order 2. \\
The first result in the phase $\beta < d + 4$, where classical diffusion no longer holds, is due to Lebeau and Puel \cite{LebPu} in the case $d = 1$, with a spectral method based on the reconnection of two branches on $\RR_-$ and $\RR_+$, and ODE techniques which provided many explicit formulas, including the diffusion coefficient. The proof based on ODE techniques and domain geometry posed a great difficulty for Puel and the author of this paper in generalizing this result to higher dimensions. This led us to seek a new, alternative method that could more easily adapt to higher dimensions. This was first done in dimension $1$ \cite{DP} using ODE methods, then in higher dimensions \cite{DP-2} using PDE and functional analysis tools, while maintaining the same strategy. This method seems robust for dealing with nonlinear problems; in fact, the method was inspired by work by Koch on the nonlinear KdV equation \cite{Koch}, and also seems interesting for potentials that depend on the spatial variable.

Between the work of \cite{LebPu} and \cite{DP,DP-2}, the fractional diffusion limit for the FP equation was obtained using a purely probabilistic approach by Fournier and Tardif, in dimension $1$ \cite{FT-d1} and then in higher dimensions \cite{FT}, in a unified presentation even for more general values of $\beta$. The same result was found using a quasi-spectral approach by Bouin and Mouhot \cite{BM}, unified for the three models LB, FP, and LFP. The approach in this latter paper, combined with the work of Ellis and Pinski \cite{EP}, was used very recently by Bouin et al.  \cite{BKM} to study hydrodynamic limits for some transport equations which admit enough conserved quantities, such as mass, momentum, and energy. Also, in the work of Gervais and Lods \cite{GL} for hydrodynamic limits for other kinetic models preserving the aforementioned three quantities, via a spectral and unified approach in the presence of a spectral gap. For the type of equations considered in this paper, only the mass is conserved. Furthermore, in the study of diffusion limits for the FP equation, in all existing results, whether in the classical or the fractional case and for all values of $\beta$, it was for a symmetric equilibrium, sometimes even given explicitly, as in \cite{NaPu,CNP,LebPu,BM,DP,DP-2}, except in the papers by probabilists \cite{FT-d1,FT} where it is under certain assumptions on $F$, including symmetry in dimension 1.  The symmetry significantly simplifies the calculations and plays a role in many aspects, particularly by giving $\int v F \ud v = 0$, even in the case where $v F \notin L^1$. In this work, we address the case of more general equilibrium, particularly non-symmetric and non-centered, using the spectral approach developed in \cite{DP-2}, while simultaneously constructing an eigenpair for the FP operator. Since the work \cite{LebPu}, it is known that if one solves the spectral problem associated with equation \eqref{FP}, then the diffusion limit becomes straightforward and is, in fact, a direct consequence. The eigenvalue provides the correct time scale as well as the diffusion coefficient, and the eigenfunction is used as a test function in the moment method. The non-symmetry of the equilibrium introduces a drift term in the kinetic equation. Such a study has been conducted in the context of the homogenization of an eigenvalue problem with drift, see for instense \cite{capdeboscq}. 

\subsection{Notations and assumptions}
Before presenting the hypotheses and stating the main results of this paper, we establish some notations that will be used throughout the paper. \\

\noindent \textbf{Notations.} 
We will keep the same notations introduced in \cite{LebPu} and used in \cite{DP-2} for the higher-dimensional case. Indeed, in order to simplify the calculations and work in the Hilbert space $L^2$, we proceed with the following change of unknown by writing
$$ f = F^\frac{1}{2}g =  M g  \qquad \mbox{ with } \qquad  M := F^\frac{1}{2} . $$ 
If $F$ is in $L^1(\RR^d)$,  the new equilibrium $M $ belongs to $ L^2(\RR^d)$.  Then, equation \eqref{fp-theta} becomes
$$ \theta(\eps) \partial_t g^\eps + \eps (v-j^\eps)\cdot \nabla_x g^\eps = -Q(g^\eps)  ,
$$
where $Q$ is a self-adjoint operator given by
$$ Q  := -\frac{1}{M}\nabla_v \cdot  \bigg(M^2\nabla_v\bigg(\frac{\cdot}{M}\bigg)\bigg) =  -\Delta_v + W(v) \qquad \text{ with } \qquad W(v) := \frac{\Delta_vM}{M}  \cdot $$
We see the equation as 
$$
\theta(\eps) \partial_t g^\eps -\eps j^\eps \cdot \na_x g^\eps = - \mathcal{L}_\eps g^\eps ,
$$
where
$$ \mathcal{L}_\eps :=  Q + \eps v \cdot \nabla_x =  -\Delta_v + W(v)  + \eps v \cdot \nabla_x \cdot $$
We operate a Fourier transform in $x$ and since the operator $Q$ has coefficient that do not depend on $x$, we get:
\begin{equation}\label{rescaled}
\theta(\eps) \partial_t \hat g^\eps - \ui \eta \ j^\eps_1 \ \hat g^\eps = - \mathcal{L}_\eta \hat g^\eps ,
\end{equation}
where
$$ \mathcal{L}_\eta := Q + \ui \eta v_1 = -\Delta_v + W(v) + \mathrm{i} \eta v_1 ,  $$
with 
\begin{equation}\label{def eta, j-1,v-1}
\eta := \eps |\xi|, \qquad j^\eps_1 := j^\eps \cdot \frac{\xi}{|\xi|}  \qquad \mbox{ and } \qquad v_1 := v \cdot \frac{\xi}{|\xi|} ,
\end{equation}
such that $\xi$ is the space Fourier variable.  \\ 

\ni \textbf{Assumptions}
\begin{assumptions}
\item \label{a1} There exists $\gamma := \frac{\beta}{2} \in (\frac{d}{2},\frac{d+4}{2})$ and there exist two positive constants $C_1 \leqslant C_2$ such that:   
$$
C_1 \langle v \rangle^{-\gamma} \leqslant M(v) \leqslant C_2 \langle v \rangle^{-\gamma} , \qquad \forall \ v \in \RR^d .
$$
\item \label{a2} The function $v \mapsto \langle v \rangle \left|\na_v\left(\frac{M}{\langle v \rangle^{2}}\right)\right|$ belongs to $L^2(\RR^d)$.
\item \label{a3} There exists a function $m$ such that, for all $s \in \RR^d\setminus \{0\}$,  one has
$$ \underset{\lambda \to 0}{\lim} \ m_\lambda (s) :=  \underset{\lambda \to 0}{\lim} \ \lambda^{-\gamma} M(\lambda^{-1} s) = m(s) . $$
\item \label{a5} There exists $\sigma \in [2,+\infty[$ and a positive constant $C$ such that, $ |W(v)| \leqslant C \langle v \rangle^{-\sigma}$. 
\end{assumptions}

\ni \textbf{Comments on the assumptions and consequences}  \medskip

\ni 1.  Assumption \ref{a1} allows us to recover the Hardy-Poincar\'e inequality which is crucial in this study.  \smallskip

\ni 2.  Assumption \ref{a1} implies that, after the change of variable $v = \eta^{-\frac{1}{3}}s$ we have
\begin{equation}\label{inegalite sur m_eta}
C_1 |s|_\eta^{-\gamma} \leqslant m_\eta(s) \leqslant C_2 |s|_\eta^{-\gamma} ,  \qquad \forall s \in \RR^d , \  \forall \eta >0 ,
\end{equation}
where $m_\eta(s) := \eta^{-\frac{\gamma}{3}}M(\eta^{-\frac{1}{3}}s)$ and $|s|_\eta^2 := (\eta^\frac{2}{3}+|s|^2)$.  Similarly, for the limit $m$, we have
\begin{equation}\label{inegalite sur m}
C_1 |s|^{-\gamma} \leqslant m(s) \leqslant C_2 |s|^{-\gamma} ,  \qquad \forall s \in \RR^d\setminus \{0\} .
\end{equation}
3.  Assumption  \ref{a2} is used to establish the compactness of an operator, which then allows the Fredholm alternative to be applied.  \smallskip

\ni 4.  \label{m_eta -->m L^2} Assumption \ref{a3} along with inequality \eqref{inegalite sur m_eta} imply that $m_\eta$ converges, when $\eta$ goes to $0$,  to $m$ in $L^2(\{|s|\geqslant r\})$ for all $r > 0$.  \smallskip

\ni 5.  Thanks to \ref{a3}, the limit of the rescaled equilibrium satisfies:   
\begin{equation}\label{m sur la sphere}
\lambda^\gamma m(\lambda s) = m(s) ,  \qquad \forall \lambda > 0 ,  \ \forall s \in \RR^d\setminus \{0\} .
\end{equation}
In particular, in the critical case $2\gamma = d+1$, the measure $s_1 m^2(s) \, \mathrm{d}s$ is invariant under dilation, that is, for any $\lambda > 0$ we have
\begin{equation} \label{mesure invariante} 
\lambda s_1 \ m^2(\lambda s) \ \ud (\lambda s) = s_1 \ m^2(s) \ \ud s . 
\end{equation}
6.  The scaling  depends on the potential $W$. In general, we distinguish three cases. For the two cases where the potential has the same scaling or is negligible compared to the Laplacian, we consider the same scaling $\eta^{-\frac{1}{3}}$. This corresponds to $\langle v \rangle^{\sigma} W \in L^\infty$ with $\sigma \in [2,+\infty[$. In the opposite case, i.e., for $\sigma < 2$, which means that the Laplacian is negligible compared to the potential, the correct scaling should be $\eta^{-\frac{1}{1+\sigma}}$. This case is discussed in Subsection \ref{subsection W}.

\subsection{Main results and idea of the proof}
The first result concerns the existence of a unique eigen-solution to the spectral problem associated with the FP operator, which bifurcates from the equilibrium $M$ with eigenvalue $0$.
\begin{theorem}[Eigen-solution for the Fokker-Planck operator with/without drift]\label{main}
\item Assume \ref{a1}-\ref{a2}.  Let  $\alpha := (\beta-d+2)/3$,  and let $\eta_0>0$ and  $\lambda_0>0$ small enough.  Then, for all $\eta\in [0,\eta_0]$,  there exists a unique eigen-couple $\big(\mu(\eta),M_\eta\big)$ in $\{\mu\in \mathbb C, |\mu|\leqslant \eta^{\frac{2}{3}}\lambda_0\}\times L^2(\RR^d,\CC)$, solution to the spectral problem 
\begin{equation}\label{M_mu,eta}
\mathcal{L}_\eta(M_{\mu,\eta})(v)=\big[-\Delta_v +W(v) + \mathrm{i} \eta v_1 \big] M_{\mu,\eta} (v)= \mu M_{\mu,\eta}(v) , \quad  v\in \RR^d .
\end{equation}
Moreover,  \smallskip

\ni 1.  the eigenfunction $M_\eta$ converges to the equilibrium $M$ in $L^2(\RR^d,\CC)$ when $\eta$ goes to $0$.  \smallskip

\ni 2.  If furthermore we assume \ref{a3} and \ref{a5} then, the  eigenvalue $\mu(\eta)$ satisfies
 \begin{equation}\label{mu(eta)}
\mu(\eta) - \ui \eta  j^\eps_1 = \kappa \eta^\alpha \big(1+O(\eta^{\alpha})\big) ,
\end{equation}
where $j^\eps_1$ is defined in \eqref{def eta, j-1,v-1}
and $\kappa$ is a positive constant given by 
\begin{equation}\label{def kappa}
\kappa := - \int_{\{|s_1|>0\}} s_1 m(s) \mathrm{Im} H_0(s)  \ud s ,
\end{equation}
and where $H_0$ is the unique solution to  
\begin{equation}\label{eq de H_0}
- \frac{1}{m} \na_s \cdot \bigg(m^2 \ \na_s\bigg(\frac{H_0}{m}\bigg)\bigg) + \mathrm{i} s_1 H_0 = 0 , \quad  \  s \in\RR^d\setminus\{0\}  ,
\end{equation}
satisfying \begin{equation}\label{condition H_0}
\int_{\{|s_1|\geqslant 1\}}|H_0(s)|^2\mathrm{d}s <+\infty \qquad  \mbox{ and } \qquad   H_0(s)\underset{0}{\sim} m(s).
\end{equation}
\end{theorem}

\begin{remark} 
 1.  The proof given for the existence of the eigen-solution $(\mu(\eta),M_\eta)$ does not depend on the potential $W$, and all we assume is Assumptions \ref{a1} and \ref{a2}. For the approximation of the eigenvalue, which gives the diffusion coefficient, we assume \ref{a3} and \ref{a5}.   \\
2.  For $M \in H^k(\RR^d,\RR)$, the eigenfunction $M_\eta$ converges to the equilibrium $M$ in $H^k(\RR^d,\CC)$. In other words, the eigenfunction $M_\eta$ has the same Sobolev regularity as its limit $M$. \\
\ni 3.  With Assumption \ref{a2},  $M$ is not necessarily in $H^1$ since
$ \na_v M = 2 \frac{v}{\langle v \rangle^2} M + \langle v \rangle^2 \na_v\left(\frac{M}{\langle v \rangle^2}\right) .$
\end{remark}

The second result concerns the diffusion limit for a kinetic Fokker-Planck equation with a general heavy-tailed equilibrium. In order to state the theorem,  we need to introduce some functional spaces. Let $\mathcal V$ be the space defined by
$$ \mathcal V := \left\{f: \mathbb{R}^d \longrightarrow \mathbb{R}; \int_{\mathbb{R}^{d}} \frac{|f|^2}{F} \ \mathrm{d}v< \infty \ \mbox{ and } \  \int_{\mathbb{R}^{d}}  \bigg|\nabla_v\bigg(\frac{f}{F}\bigg)\bigg|^2  F \mathrm{d}v < \infty \right\} ,$$ 
with $\mathcal V'$ being its dual, and let
$$ Y := \left\{ f \in L^2\big([0,T]\times\RR^d; \mathcal V\big); \ \theta(\eps) \pa_t f + \eps (v -j^\eps)\cdot \nabla_x f \in L^2\big([0,T]\times\RR^d; \mathcal V'\big) \right\} . $$

\begin{theorem}[Fractional diffusion for the Fokker-Planck equation with(-out) drift]\label{main2}
Assume \ref{a1}--\ref{a5} and assume that  $f_0$ is a non-negative function in  $ L^2_{F^{-1}}(\RR^{2d})\cap L^\infty_{F^{-1}}(\RR^{2d})$. Let $f^\varepsilon$ be the solution of \eqref{fp-theta} in $Y$ with initial data $f_0$ and $\theta(\varepsilon)=\varepsilon^{\alpha} $. Let $ \kappa$ be the constant given by \eqref{def kappa}. Then, $f^\varepsilon$ converges weakly star in $L^\infty\big([0,T]; L^2_{F^{-1}}(\RR^{2d})\big)$ towards $\rho(t,x) F(v)$, where $\rho$ is the unique solution to 
\begin{equation}
\partial_t\rho +\kappa (-\Delta_x)^{\frac{\alpha}{2}}\rho =0,\quad \rho(0,x)=\int_{\RR^d} f_0 \ \mathrm{d}v .
 \end{equation}
\end{theorem}

\subsubsection*{Ideas of the proof and outline of the paper.}
The strategy and main steps of the proof are the same as those in \cite{DP-2}. The differences lie in the proof based on the Lax-Milgram theorem, used to show the existence of solutions, and in the approximation of the eigenvalue.

The proof of Theorem \ref{main} is carried out in two main steps, both based on the implicit function theorem (IFT). First, we consider the \emph{penalized equation}:
\begin{equation}\label{eq penalisee1}
\left\{ \begin{array}{l}
\big[-\Delta_v + W(v) + \mathrm{i} \eta v_1 \big] M_{\mu,\eta} (v)= \mu M_{\mu,\eta}(v) - \langle M_{\mu,\eta}-M,\Phi\rangle \Phi(v), \quad  v\in \mathbb{R}^d,\\
\\
M_{\mu,\eta} \in L^2(\mathbb{R}^d,\mathbb{C}),
\end{array}\right.
\end{equation}
where $\Phi$ is given by
$$ \Phi  := \bigg(\int_{\mathbb{R}^d}\frac{M^2}{\langle v \rangle^2} \mathrm{d} v\bigg)^{-1} \frac{M}{\langle v \rangle^2} \cdot $$

The aim of the first step is to show the existence of a unique solution for equation \eqref{eq penalisee1} for fixed $\eta$ and $\mu$, which is the purpose of Section \ref{section penalized eq}. The idea is to transform equation \eqref{eq penalisee1} into a fixed point problem in the form of the identity plus a compact operator.  To achieve this, we need to look at the equation differently, adding a term to both sides to obtain an invertible operator with a continuous and compact inverse. This is where we use the Lax-Milgram theorem.  When applying this theorem, we do not use the same norm as in \cite{DP-2}. In the latter reference, the potential $W$ was given explicitly, and the definition of the norm as well as the invertibility of the operator were based on the decomposition of $W$.  However, in this paper, we work with a norm that does not depend on $W$, i.e., it does not depend on its expression or its decomposition.

In the second step, in order to ensure that the penalized term vanishes, we have to choose $\mu(\eta)$  via the IFT around the point $(\mu,\eta) = (0, 0)$. The study of this constraint is the subject of Subsection \ref{subsection L2}. Subsections \ref{subsection H_eta} and \ref{subsection vp} are devoted to the approximation of the eigenvalue and the computation of the diffusion coefficient. The fact that $M$ is not necessarily symmetric in this paper makes the calculations more complex compared to \cite{DP-2} and sometimes requires adding a term, corresponding to a drift in the kinetic equation, to close the estimates and obtain a finite diffusion coefficient.  Additionally, in this paper, we address the case $\beta = d + 1$ which has a different scaling.

Section \ref{section diff frac} is devoted to the proof of Theorem \ref{main2}. It consists of two subsections: a priori estimates and limiting process in the weak formulation of equation \eqref{rescaled}.  Finally,  Section \ref{section comments} contains some comments on the assumptions with examples of equilibria and potentials, as well as the link to the work of \cite{FT}, and a discussion of the case where the potential dominates the Laplacian.

\section{Existence of solutions for the penalized equation}\label{section penalized eq}
The objective of this section is to prove the existence of a unique solution for the penalized equation \eqref{eq penalisee1}, for fixed $\mu$ and $\eta$ that are sufficiently small. We start this section by some notations and definition of the considered operators. Let $\mu =\lambda\eta^\frac{2}{3}$ with $\lambda \in \mathbb{C}$ and let denote by $L_{\lambda,\eta}$ the operator
$$ L_{\lambda,\eta} := Q + V + \mathrm{i} \eta v_1 -  \lambda \eta^\frac{2}{3} ,$$
where $$  Q :=  -\frac{1}{M}\na_v\left(M^2\na_v\left(\frac{\cdot}{M}\right)\right) \quad \mbox{ and } \quad V(v) := \frac{1}{\langle v \rangle^{\sigma}} \quad \mbox{ with } \quad \sigma \in ]2,2\gamma-d+4[ .$$
We rewrite equation \eqref{eq penalisee1} as follows
\begin{equation}\label{eq penalisee2}
\left\{ \begin{array}{l}
L_{\lambda,\eta} (M_{\lambda,\eta})  =  V(v) M_{\lambda,\eta} - \langle M_{\lambda,\eta}-M,\Phi\rangle \Phi , \ v\in \RR^d ,\\
\\
M_{\lambda,\eta} \in L^2(\RR^d,\CC).
\end{array}\right.
\end{equation}
The two equations \eqref{eq penalisee1} and \eqref{eq penalisee2} are equivalent. 
\begin{remark}\label{symetrie de la sol}
\item 1.  Since $L_{\lambda,0}$ does not depend on $\lambda$, let's denote it by $L_0$: $\ L_0:=L_{\lambda,0}$.  \smallskip

\ni 2.  The term $V$ was introduced to make the operator $L_{0}$ invertible. See Remark \ref{tilde Q = Q + V}.
\end{remark}

\subsection{A Lax-Milgram theorem and invertibility of the operator}
The purpose of this subsection is to show that the operator $L_{\lambda,\eta}$ defined above is invertible in an appropriate Hilbert space. Before presenting the strategy, we first define the suitable Hilbert space and the sesquilinear form associated with $L_\eta$.
\begin{definition} 
\item 1.  We define the Hilbert space $\mathcal{H}_\eta$ as being the completion of the space $C_c^\infty(\RR^d,\CC)$ for the norm $\|\cdot\| _{\mathcal{H}_\eta}$ induced from the inner product $\langle \cdot,\cdot\rangle_{\mathcal{H}_\eta}$
$$
\mathcal{H}_\eta := \overline{\big\{ \psi \in C_c^\infty(\RR^d,\CC) ; \ \| \psi \|_{\mathcal{H}_\eta}^2:= \langle \psi,\psi\rangle_{\mathcal{H}_\eta} < +\infty \big\} } ,
$$
where
$$
\langle \psi,\phi\rangle_{\mathcal{H}_\eta}:=  \int_{\RR^d} \nabla_v \bigg(\frac{\psi}{M}\bigg) \cdot \nabla_v \bigg(\frac{\bar\phi}{M}\bigg) \ M^2 \mathrm{d}v + \int_{\RR^d}   \frac{\psi \bar \phi }{\langle v \rangle^2} \ \mathrm{d}v +  \eta \int_{\RR^d} |v_1| \psi \bar \phi \ \mathrm{d}v .
$$
We have the embeddings
$$ \mathcal{H}_\eta \subseteq \mathcal{H}_{\eta^*} \subseteq \mathcal{H}_{0}  ,   \quad \forall \ 0 \leqslant \eta^* \leqslant \eta $$
since $\ \| \cdot \|_{\mathcal{H}_0} \leqslant \| \cdot \|_{\mathcal{H}_{\eta^*}} \leqslant \| \cdot \|_{\mathcal{H}_\eta} \ $  for all $\ 0\leqslant \eta^* \leqslant \eta$.   \smallskip

\ni 2.  We define the sesquilinear form $a$ on $\mathcal{H}_\eta\times\mathcal{H}_\eta$ by
$$
a(\psi,\phi):= \int_{\RR^d} \nabla_v \bigg(\frac{\psi}{M}\bigg) \cdot \nabla_v \bigg(\frac{\bar\phi}{M}\bigg) \ M^2 \mathrm{d}v + \int_{\RR^d}  V \psi \bar \phi  \ \mathrm{d}v +  \mathrm{i} \eta \int_{\RR^d} v_1 \psi \bar \phi \ \mathrm{d}v  - \lambda \eta^{\frac{2}{3}} \int_{\RR^d} \psi \bar \phi \ \mathrm{d}v  .
$$
\end{definition}
\begin{remark}\label{tilde Q = Q + V}
\item 1.  Note that the sesquilinear form $a$ depends on $\lambda$ and $\eta$ and in order to simplify the notation, we omit the subscript when no confusion is possible. \smallskip

\ni 2.  The operator $\tilde Q := Q + V$ is dissipative and one has:
$$  \int_{\RR^d} \tilde Q(\psi)\psi \ \mathrm{d}v = \int_{\RR^d} Q(\psi)\psi \ \mathrm{d}v + \int_{\RR^d} V |\psi|^2 \ \mathrm{d}v =  \int_{\RR^d} \bigg| \nabla_v \bigg(\frac{\psi}{M}\bigg)\bigg|^2 M^2 + V |\psi|^2 \ \mathrm{d}v\geqslant 0 . $$
3.  For $\eta=0$ in $\|\cdot\|_{\mathcal{H}_\eta}$, the term arising from the operator $Q$ alone, i.e. $\int_{\RR^d} \big|\nabla_v\big(\frac{\varphi}{M}\big)\big|^2M^2 \ud v$, does not define a norm. Also, the operator $Q+\ui \eta v_1$ is not invertible for sufficiently small $\eta$, since the kernel of $Q$ contains $M$.
\end{remark}
\begin{lemma}\label{equivalence des normes}
The two norms $\|\cdot\|_{\mathcal{H}_0}$ and $\|\cdot\|_*$ defined by
$$ \| \psi \|_{{\mathcal{H}}_0}^2 := \int_{\RR^d} \bigg|\nabla_v\bigg(\frac{\psi}{M}\bigg) M \bigg|^2 + \frac{|\psi|^2}{\langle v \rangle^2} \ \ud v   \quad \mbox{ and } \quad   \|\psi\|_*^2 := \int_{\RR^d} \bigg|\nabla_v\bigg(\frac{\psi}{M}\bigg) M \bigg|^2 + V|\psi|^2 \ \ud v$$
are equivalent, i.e., there are two positive constants $C_1$ and $C_2$ such that
$$ C_1 \| \psi \|_{\mathcal{H}_0} \leqslant \| \psi \|_* \leqslant C_2 \| \psi \|_{\mathcal{H}_0} , \quad \forall \psi \in \mathcal{H}_0 . $$
\end{lemma}
\begin{remark}
Regarding the choice of the function $V$ introduced at the beginning of this section, it can be any function such that
 \begin{itemize}
 \item Lemma \ref{equivalence des normes} is true.
 \item The function $\langle v \rangle^2 V$ belongs to $C^1_0(\RR^d,\RR)$.
 \end{itemize}
\end{remark}
To prove Lemma \ref{equivalence des normes}, we need the Hardy-Poincar\'e inequality that we recall in the following 
\begin{lemma}[Hardy-Poincar\'e inequality \cite{BDGV}]\label{Hardy-Poincare}
Let $d\geqslant 1$ and $\alpha_*=\frac{2-d}{2}$. For any $\alpha<0$,  and $\alpha \in (-\infty,0)\setminus \{ \alpha_*\}$ for $d\geqslant 3$, there is a positive constant $\Lambda_{\alpha,d}$ such that
\begin{equation}
\Lambda_{\alpha,d} \int_{\RR^d} |f|^2(D+|x|^2)^{\alpha-1}\mathrm{d}x \leqslant \int_{\RR^d} |\nabla f|^2 (D+|x|^2)^{\alpha}\mathrm{d}x
\end{equation}
holds for any function $f \in H^1\big((D+|x|^2)^{\alpha}\mathrm{d}x\big)$ and any $D \geqslant 0$, under the additional condition $\int_{\RR^d} f (D+|x|^2)^{\alpha-1}\mathrm{d}x = 0$ and $D>0$ if $\alpha<\alpha_*$.
\end{lemma}
\begin{remark}
For $\ds f=\frac{g}{M}$, $D=1$ and $\alpha = -\gamma$ in the previous lemma, the inequality becomes
\begin{equation}\label{ineg.  Hardy-Poincare}
\Lambda_{\alpha,d} \int_{\RR^d} \frac{|g|^2 }{\langle v\rangle^2} \mathrm{d}v \leqslant \int_{\RR^d} \bigg|\nabla_v \bigg(\frac{g}{M}\bigg)\bigg|^2 M^2 \ \mathrm{d}v ,
\end{equation}
under the orthogonality condition
\begin{equation}\label{condition d'orthogonalite}
 \int_{\RR^d}  \frac{g M}{\langle v \rangle^2} \mathrm{d}v = 0 ,
\end{equation}
since $\  -\gamma < \frac{2-d}{2} =: \alpha_* \ $ for $\ \gamma \in (\frac{d}{2},\frac{d+4}{2})$.  Moreover, if we  denote by $\mathcal{P}(g)$ the following scalar:
$$ \mathcal{P}(g) := \left(\int_{\RR^d} \frac{M^2}{\langle v \rangle^2} \mathrm{d}v\right)^{-1} \int_{\RR^d}  \frac{g M}{\langle v \rangle^2} \mathrm{d}v ,$$
then,  inequality \eqref{ineg. Hardy-Poincare} can be written, for all $g \in \mathcal{H}_0$, as follow
\begin{equation}\label{Hardy-Poincare avec orthogonalite}
\Lambda_{\gamma,d} \int_{\RR^d}  \frac{\big| g -\mathcal{P}(g)M \big|^2}{\langle v\rangle^2} \mathrm{d}v \leqslant \int_{\RR^d} \bigg|\nabla_v \bigg(\frac{g}{M}\bigg)\bigg|^2 M^2 \ \mathrm{d}v .
\end{equation}
\end{remark}
\bpl \ref{equivalence des normes}.  The inequality $\| \psi \|_* \leqslant C_2 \| \psi \|_{\mathcal{H}_0}$ is trivial since $V(v) \leqslant \langle v \rangle^{-2}$ for $\sigma >2$.  It remains to show the second inequality, i.e. there exists $C_1>0$ such that  $\| \psi \|_{{\mathcal{H}}_0} \leqslant C_1 \| \psi \|_*$.  Let $\psi$ such that $\| \psi \|_* < \infty$.   First, since $M \in L^2(\RR^d)$ then, by the Cauchy-Schwarz inequality we get
$$\bigg| \int_{\RR^d} \frac{\psi M}{\langle v \rangle^2} \mathrm{d}v \bigg| \leqslant \left(\int_{\RR^d} \frac{|\psi|^2}{\langle v \rangle^\sigma} \mathrm{d}v\right)^{\frac{1}{2}}  \left(\int_{\RR^d} \frac{M^2}{\langle v\rangle^{4-\sigma}} \mathrm{d}v\right)^{\frac{1}{2}} \leqslant C \| \psi \|_*  ,$$
where 
$ C:=\left(\int_{\RR^d} \frac{M^2}{\langle v\rangle^{4-\sigma}} \mathrm{d}v\right)^{\frac{1}{2}} \lesssim\left(\int_{\RR^d} \langle v\rangle^{-2\gamma-4+\sigma} \mathrm{d}v\right)^{\frac{1}{2}} \lesssim 1 $ for $\sigma < 2\gamma -d + 4   $.
Now, since the function $\psi - \mathcal{P}(\psi)M$ satisfies condition \eqref{condition d'orthogonalite} then,   inequality \eqref{Hardy-Poincare avec orthogonalite} can be used and therefore
\begin{align*}
\int_{\RR^d} \frac{|\psi|^2}{\langle v \rangle^2} \mathrm{d}v &= \int_{\RR^d} \frac{|\psi- \mathcal{P}(\psi)M + \mathcal{P}(\psi)M|^2}{\langle v \rangle^2} \mathrm{d}v \\ 
&\leqslant 2 \left( \Lambda_{\gamma,d}^{-1} \int_{\RR^d} \bigg|\nabla_v \bigg(\frac{\psi}{M}\bigg)\bigg|^2 M^2 \ \mathrm{d}v+ |\mathcal{P}(\psi)|^2 \int_{\RR^d} \frac{M^2}{\langle v\rangle^2}\mathrm{d}v \right)  \\
&\leqslant 2 \left( \Lambda_{\gamma,d}^{-1}  \| \psi \|_*^2 + C^2 \left(\int_{\RR^d} \frac{M^2}{\langle v \rangle^2} \mathrm{d}v\right)^{-1} \| \psi \|_*^2 \right) .
\end{align*}
Hence, $  \|\psi\|_{\mathcal{H}_0}^2 \leqslant C_1 \|\psi\|_*^2 \ $  with $C_1$ a positive constant that depends only on $\gamma$ and $d$. 
\epl

To invert the operator $L_{\lambda,\eta}$, we will employ an elaborate version of the Lax-Milgram theorem. More precisely, we will demonstrate the continuity of the sesquilinear form $a$ on the Hilbert space $\mathcal{H}_\eta$. Subsequently,  by virtue of the Riesz representation theorem, this form will correspond to a linear operator $A$. Finally, we will proceed to invert the operator $A$. For this final step, we will follow the same procedure as outlined in \cite{DP-2}. However, note that the proof is slightly modified since we are not working with the same norm here. It is important to observe that the norm $\|\cdot\|_{\mathcal{H}_\eta}$ defined above does not depend on the expression of $W$.  \\

Let's start with the following lemma, which gives us the continuity of $a$.
\begin{lemma}\label{Poincare1} Let $\eta>0$ be fixed. Then, there exists a constant $C_0>0$, independent of $\eta$ such that the following inequality holds true
\begin{equation}\label{Poincare1 ineq}
\| \psi \|_{L^2(\RR^d)} \leqslant C_0 \eta^{-\frac{1}{3}} \| \psi \|_{\mathcal{H}_\eta}  , \quad \forall \psi \in \mathcal{H}_\eta .
\end{equation}
\end{lemma}
\bp First, we set $v :=(v_1,v')\in \RR\times\RR^{d-1}$ and we define the following sets:
$$ A_\eta := \{|v_1|\leqslant \eta^{-\frac{1}{3}}\} , \quad B_\eta := \{|v|\leqslant \eta^{-\frac{1}{3}}\} \quad \mbox{ and } \quad C_\eta := \{|v_1|\leqslant \eta^{-\frac{1}{3}} \leqslant |v'|\} . $$
Now,  using the fact that $A_\eta = B_\eta \cup C_\eta$, we decompose the integral on $\RR^d$ as follows:
$$ \int_{\RR^d} |\psi|^2 \ud v = \int_{A_\eta} |\psi|^2 \ud v + \int_{A_\eta^c} |\psi|^2 \ud v = \int_{B_\eta} |\psi|^2 \ud v + \int_{C_\eta} |\psi|^2 \ud v + \int_{A_\eta^c} |\psi|^2 \ud v .$$
(i) We have on $A_\eta^c$:
$$ \int_{A_\eta^c} |\psi|^2 \ud v \leqslant \int_{A_\eta^c} \eta^{\frac{1}{3}} |v_1| |\psi|^2 \ud v \leqslant \eta^{-\frac{2}{3}}\int_{A_\eta^c} \eta |v_1| |\psi|^2 \ud v .$$
(ii) For the part $B_\eta$, we write:
$$ \int_{B_\eta} |\psi|^2 \ud v = \int_{B_\eta} \langle v \rangle^2 \frac{|\psi|^2}{\langle v \rangle^2} \ud v \lesssim \eta^{-\frac{2}{3}} \int_{B_\eta} \frac{|\psi|^2}{\langle v \rangle^2} \ud v . $$
(iii) Finally, for the part $C_\eta$, it is summarized in the Poincar\'e-type inequality given in the following 
\begin{lemma}\label{Poincare2} Let $\eta>0$ be fixed. Then, there exists a constant $C>0$, independent of $\eta$ such that the following inequality holds for all $\psi \in \mathcal{H}_\eta$ 
\begin{equation}\label{Poincare2 ineq}
\int_{C_\eta} |\psi|^2 \ud v \leqslant C \eta^{-\frac{2}{3}} \left( \int_{\RR^d} \bigg| \na_v\bigg(\frac{\psi}{M}\bigg)\bigg|^2 M^2  \ \ud v + \eta \int_{A_\eta^c} |v_1| |\psi|^2 \ \ud v \right) .
\end{equation}
\end{lemma}
\bp  Let $\chi \in C^\infty(\mathbb{R})$ such that $0 \leqslant \chi \leqslant 1$, $\chi \equiv 1$ in $B(0,1)$, and $\chi \equiv 0$ outside of $B(0,2)$. We define $\chi_\eta$ by $\chi_\eta(v_1) := \chi(\eta^\frac{1}{3} v_1)$, and we define the set $\tilde C_\eta$ as the extension of $C_\eta$ in the direction of $v_1$, $\tilde C_\eta := \{|v_1| \leqslant 2\eta^{-\frac{1}{3}} \leqslant 2|v'|\}$. We have: $\ \int_{C_\eta} |\psi|^2 \ud v \leqslant \int_{\tilde C_\eta} |\chi_\eta\psi|^2 \ud v $.  Now, by writing $\frac{\chi_\eta \psi}{M} = \int_{-2\eta^{-\frac{1}{3}}}^{v_1} \partial_{w_1} \big(\frac{\chi_\eta \psi}{M}\big) \ud w_1$, and using the fact that $\frac{M(v_1,v')}{M(w_1,v')} \lesssim 1$ on $\tilde C_\eta$, since $\langle v \rangle^{\gamma} M(v)$ is bounded above and below by two positive constants, we write:
\begin{align*}
|\chi_\eta \psi |^2 = \bigg| M(v_1,v')\int_{-2\eta^{-\frac{1}{3}}}^{v_1} \partial_{w_1}\bigg(\frac{\chi_\eta \psi}{M}\bigg) \ud w_1 \bigg|^2 &=\bigg| \int_{-2\eta^{-\frac{1}{3}}}^{v_1}\frac{M(v_1,v')}{M(w_1,v')} \partial_{w_1} \bigg(\frac{\chi_\eta \psi}{M}\bigg) M(w_1,v') \ud w_1 \bigg|^2 \\
&\lesssim  \left( \int_{-2\eta^{-\frac{1}{3}}}^{v_1} \bigg|\partial_{w_1} \bigg(\frac{\chi_\eta \psi}{M}\bigg) M(w_1,v')\bigg| \ud w_1 \right)^2 \\
&\lesssim  \eta^{-\frac{1}{3}}  \int_{-2\eta^{-\frac{1}{3}}}^{2\eta^{-\frac{1}{3}}} \bigg|\partial_{v_1} \bigg(\frac{\chi_\eta \psi}{M}\bigg)\bigg|^2 M^2  \ud v_1 .
\end{align*}
Hence,
\begin{equation}\label{P1}
\int_{\tilde C_\eta} |\chi_\eta \psi|^2 \ud v \lesssim \eta^{-\frac{2}{3}} \int_{\tilde C_\eta} \bigg|\partial_{v_1} \bigg(\frac{\chi_\eta \psi}{M}\bigg)\bigg|^2 M^2  \ud v .
\end{equation}
Now from the equality $\partial_{v_1}\big(\frac{\chi_\eta \psi}{M}\big) M = \chi_\eta' \psi + \chi_\eta \partial_{v_1} \big(\frac{\psi}{M}\big)M$,  we obtain
\begin{equation}\label{P2}
 \int_{\tilde C_\eta} \bigg|\partial_{v_1} \left(\frac{\chi_\eta \psi}{M}\right)\bigg|^2 M^2  \ud v \lesssim \int_{\tilde C_\eta} |\chi_\eta'|^2 |\psi|^2  \ud v + \int_{\tilde C_\eta} \bigg|\partial_{v_1} \left(\frac{\psi}{M}\right)\bigg|^2 M^2  \ud v .
\end{equation}
Finally, since $\chi_\eta' \equiv 0$ except on $\tilde C_\eta \setminus C_\eta = \{\eta^{-\frac{1}{3}} \leqslant |v_1| \leqslant 2 \eta^{-\frac{1}{3}} \leqslant 2|v'|\} \subset A_\eta^c$ then,
$$
|\chi_\eta'(v_1)| = \left|\eta^{\frac{1}{3}} \chi'(\eta^{\frac{1}{3}}v_1)\right| \leqslant \eta^{\frac{1}{3}} \underset{t\in[1,2]}{\sup} |t\chi'(t)| \quad \mbox{ and } \quad |\chi_\eta'(v_1)|^2 \leqslant \eta |v_1| \bigg(\underset{t\in[1,2]}{\sup} |t\chi'(t)| \bigg)^2 1_{A^c_\eta}.
$$
Therefore,
\begin{equation}\label{P3}
\int_{\tilde C_\eta} |\chi_\eta'|^2 |\psi|^2  \ud v \lesssim \eta \int_{A_\eta^c} |v_1| |\psi|^2  \ud v . 
\end{equation}
Hence, inequality \eqref{Poincare2 ineq} holds by combining \eqref{P1}, \eqref{P2}, and \eqref{P3}.
\epl 
To conclude the proof of Lemma \ref{Poincare1}, by summing up (i), (ii), and (iii), we obtain 
\begin{equation}\label{*}
\int_{\RR^d} |\psi|^2 \ud v \leqslant C \eta^{-\frac{2}{3}} \left( \int_{\RR^d}\bigg|\na_v\bigg(\frac{\psi}{M}\bigg)\bigg|^2M^2 \ud v + \int_{\RR^d} \frac{|\psi|^2}{\langle v\rangle^2} \ud v + \eta \int_{A_\eta^c} |v_1| |\psi|^2 \ud v\right) .
\end{equation}
\ep

As a consequence of Lemma \ref{Poincare1} applied to the term $\lambda \eta^{\frac{2}{3}}\int_{\RR^d} |\psi|^2 \ud v$, with the Cauchy-Schwarz inequality, we have
\begin{lemma} The sesquilinear form $a$ is continuous on $\mathcal{H}_\eta\times\mathcal{H}_\eta$. Moreover, there exists a constant $C>0$, independent of $\lambda$ and $\eta$ such that, for all $\psi,\phi \in \mathcal{H}_\eta$ one has
$$
|a(\psi,\phi)| \leqslant C \| \psi \|_{\mathcal{H}_\eta}\| \phi \|_{\mathcal{H}_\eta} .
$$
\end{lemma}

\begin{remark}\label{A representant de a}
By application of Riesz's theorem to continuous sesquilinear forms, there exists a continuous linear map $ A_{\lambda,\eta}\in \mathcal{L}({\mathcal {H}_\eta})$ such that $\ a(\psi,\phi) = \langle A_{\lambda,\eta}\psi, \phi\rangle_{\mathcal{H}_\eta}$ for all $\psi, \phi \in \mathcal{H}_\eta$.  \\
Note that $A_{\lambda,\eta}$ depends on $\lambda$ and $\eta$ since the form $a$ depends on these last parameters.
\end{remark}
\begin{lemma}\label{coercivite de a} Let $\eta>0$ and $\lambda\in\CC$ fixed, such that $|\lambda|\leqslant \lambda_0$ with $\lambda_0$ small enough. Let $A_{\lambda,\eta}$ be the linear operator representing the sesquilinear form $a$. Then, there exists a constant $C>0$, independent of $\lambda$ and $\eta$ such that
\begin{equation}\label{psi < A psi}
\| \psi \|_{\mathcal{H}_\eta} \leqslant C \| A_{\lambda,\eta} \psi \|_{\mathcal{H}_\eta} , \quad \forall \psi \in \mathcal{H}_\eta .
\end{equation}
\end{lemma}
\begin{remark}
Note that for $\lambda = 0$,  $ \ds a(\psi,\psi)\neq \|\psi\|_{*}^2 + \eta \int_{\RR^d} |v_1| |\psi|^2 \ud v $.
\end{remark}
\bp We have for all $a, b \in \RR$ and $z\in \CC$: $ |a+ \mathrm{i} b + z| \geqslant |a|-|z|$. Now, applying this inequality to $|a(\psi,\psi)|$,  using Lemma \ref{Poincare1} for the term $\lambda \eta^{\frac{2}{3}} \int |\psi|^2 \ud v$, and using Lemma \ref{equivalence des normes},  we write
\begin{align*}
|a(\psi,\psi)| &= \left| \int_{\RR^d} \left(\bigg|\nabla_v \bigg(\frac{\psi}{M}\bigg)\bigg|^2 M^2 + V|\psi|^2 +  \mathrm{i} \eta v_1 |\psi|^2 -\lambda \eta^{\frac{2}{3}}|\psi|^2 \right) \mathrm{d}v \right|  \\
& \geqslant \left| \int_{\RR^d} \left(\bigg|\nabla_v \bigg(\frac{\psi}{M}\bigg)\bigg|^2 M^2 + V|\psi|^2 \right)\ud v \right| - |\lambda| \eta^{\frac{2}{3}}\| \psi \|_{L^2}^2 \\
&\geqslant C_1^2 \|\psi\|^2_{\mathcal{H}_0} - C_0|\lambda| \|\psi\|^2_{\mathcal{H}_\eta} .
\end{align*}
Then,  since $|a(\psi,\psi)|=|\langle A_{\lambda,\eta} \psi,\psi \rangle_{\mathcal{H}_\eta}|\leqslant \|A_{\lambda,\eta} \psi\|_{\mathcal{H}_\eta} \|\psi\|_{\mathcal{H}_\eta}$, we get
\begin{equation}\label{a}
C_1^2 \|\psi\|^2_{\mathcal{H}_0}  \leqslant \|A_{\lambda,\eta} \psi\|_{\mathcal{H}_\eta} \|\psi\|_{\mathcal{H}_\eta} + C_0 |\lambda| \|\psi\|^2_{\mathcal{H}_\eta} .
\end{equation}
Let denote 
$$ I^\eta_1 := \int_{\{ |v_1| \leqslant \eta^{-\frac{1}{3}}\}} \eta |v_1||\psi|^2\mathrm{d}v \quad \mbox{ and } \quad  I^\eta_2 := \int_{\{ |v_1| \geqslant \eta^{-\frac{1}{3}}\}} \eta |v_1||\psi|^2\mathrm{d}v . $$
Note that $\|\psi\|^2_{\mathcal{H}_\eta} = \|\psi\|^2_{\mathcal{H}_0} + I^\eta_1 + I^\eta_2$.  To establish  inequality \eqref{psi < A psi}, we will control $I_1^\eta + I^\eta_2$ by $\|A_{\lambda,\eta} \psi\|_{\mathcal{H}_\eta} \|\psi\|_{\mathcal{H}_\eta} + |\lambda| \|\psi\|^2_{\mathcal{H}_\eta}$ and combine it with \eqref{a}.
The estimation of $I^\eta_1$ and $I^\eta_2$ follows the same spirit as that of Lemma \ref{Poincare1}. With the same notations for the sets $A_\eta$, $B_\eta$, and $C_\eta$, we have:\\

\ni \textbf{Step 1: Estimation of $I^\eta_1$.} We have thanks to inequality \eqref{*}
$$ I^\eta_1 := \int_{A_\eta} \eta |v_1| |\psi|^2 \ud v  \leqslant \eta^{\frac{2}{3}} \int_{A_\eta} |\psi|^2 \ud v \leqslant C \left( \int_{\RR^d}\bigg|\na_v\bigg(\frac{\psi}{M}\bigg)\bigg|^2M^2 \ud v + \int_{\RR^d} \frac{|\psi|^2}{\langle v\rangle^2} \ud v + I^\eta_2 \right)  .$$
Hence,
\begin{equation}\label{I^eta_1 < I^eta_2 + |psi|_H_0^2}
I^\eta_1 \leqslant C \left(\|\psi\|_{\mathcal{H}_0}^2 + I^\eta_2 \right) .
\end{equation}
\textbf{Step 2: Estimation of $I^\eta_2$.} Let $\chi_\eta$ this time be the function defined by $\chi_\eta(v_1):=\chi(\eta^{\frac{1}{3}}v_1)$ with $\chi \in C^\infty(\RR)$ such that: $-1 \leqslant \chi \leqslant 1$, $\chi \equiv -1$ on $]-\infty,-1]$,  $\chi \equiv 1$ on $[1,+\infty[$ and $\chi \equiv 0$ in $B(0,\frac{1}{2})$.  Let $D_\eta$ denote the set $D_\eta := \{ |v_1| \geqslant \frac{1}{2}\eta^{-\frac{1}{3}}\}$. Then,
$$ I^\eta_2 := \int_{\{ |v_1| \geqslant \eta^{-\frac{1}{3}}\}} \eta |v_1||\psi|^2  \mathrm{d}v \leqslant  \int_{D_\eta} \eta v_1 \chi_\eta \psi \bar \psi \ \mathrm{d}v .$$
By integrating the equation of $\psi$ against $\chi_\eta \bar \psi$ over $D_\eta$ and taking the imaginary part, we obtain
$$ \int_{D_\eta} \eta v_1 \chi_\eta \psi \bar \psi \ \mathrm{d}v = \mathrm{Im}\left( a(\psi,\chi_\eta \psi) - \int_{D_\eta} \bigg[ \na_v \bigg(\frac{\psi}{M}\bigg)\cdot \na_v\bigg(\frac{\chi_\eta \bar \psi}{M}\bigg) M^2 - \lambda \eta^{\frac{2}{3}}  \chi_\eta \psi \bar \psi  \bigg] \mathrm{d}v \right) . $$  
$\bullet$ Let's start with the last term. By Lemma \ref{Poincare1} we get 
\begin{equation} \label{J1}
\bigg|\mathrm{Im} \lambda \eta^{\frac{2}{3}} \int_{D_\eta} \chi_\eta \psi \bar \psi \ \mathrm{d}v \bigg| \leqslant C_0 |\lambda| \ \|\psi\|^2_{\mathcal{H}_\eta} .
\end{equation}
$\bullet$ For the first term, by Cauchy-Schwarz: $|\mathrm{Im}\ a(\psi,\chi_\eta \psi)|\leqslant \|A_{\lambda,\eta} \psi\|_{ \mathcal{H}_\eta}\|\chi_\eta\psi\|_{\mathcal{H}_\eta}$. We still need to estimate $\|\chi_\eta\psi\|_{\mathcal{H}\eta}$. For this purpose, we have
$$ \|\chi_\eta\psi\|_{\mathcal{H}\eta}^2 = \int_{\RR^d} \bigg| \na_v \bigg(\frac{\chi_\eta \psi}{M}\bigg)\bigg|^2 M^2 \ud v +  \int_{\RR^d} \frac{|\chi_\eta \psi|^2}{\langle v \rangle^2} \ud v + \eta \int_{\RR^d} |v_1| |\chi_\eta \psi|^2\ud v .$$
For the two terms on the right,  since $|\chi_\eta| \leqslant 1$ then
$$ \int_{\RR^d} \frac{|\chi_\eta \psi|^2}{\langle v \rangle^2} \ud v + \eta \int_{\RR^d} |v_1| |\chi_\eta \psi|^2\ud v \leqslant \int_{\RR^d} \frac{|\psi|^2}{\langle v \rangle^2} \ud v + \eta \int_{\RR^d} |v_1| |\psi|^2\ud v .$$
For the left-hand term, since $\big| \na_v \big(\frac{\chi_\eta \psi}{M}\big)\big|^2 M^2 = |\chi_\eta'\psi|^2 + \big| \na_v \big(\frac{\psi}{M}\big)\big|^2 M^2 + 2\Re \big[\chi_\eta'\psi\pa_{v_1}\big(\frac{\psi}{M}\big)M\big]$,  and as in the proof of Lemma \ref{Poincare2}, since $\chi_\eta' \equiv 0$ except on the set $\{\frac{1}{2}\eta^{-\frac{1}{3}} \leqslant |v_1| \leqslant \eta^{\frac{1}{3}}\} \subset A_\eta$ where $|\chi_\eta' \psi|^2 \leqslant C \eta |v_1| |\psi|^2 1_{A_\eta}$, with $C := \big(\underset{t\in[\frac{1}{2},1]}{\sup} |t \chi'(t)|\big)^2$,  then we write
\begin{align*}
\int_{\RR^d} \bigg| \na_v \bigg(\frac{\chi_\eta \psi}{M}\bigg)\bigg|^2 M^2 \ud v &\leqslant 2 \left( \int_{A_\eta} |\chi_\eta' \psi|^2 \ud v + \int_{\RR^d} \bigg| \na_v \bigg(\frac{\psi}{M}\bigg)\bigg|^2 M^2 \ud v   \right) \\
&\leqslant 2C I^\eta_1 + 2 \|\psi\|_{\mathcal{H}_0}^2 \leqslant C' \|\psi\|_{\mathcal{H}_\eta}^2 .
\end{align*}
Hence $\|\chi_\eta\psi\|_{\mathcal{H}_\eta} \leqslant C \|\psi\|_{\mathcal{H}_\eta}$, and therefore
\begin{equation}\label{J2}
|\mathrm{Im}\ a(\psi,\chi_\eta \psi)|\leqslant \|A_{\lambda,\eta} \psi\|_{ \mathcal{H}_\eta}\|\chi_\eta\psi\|_{\mathcal{H}_\eta} \leqslant C  \|A_{\lambda,\eta} \psi\|_{ \mathcal{H}_\eta} \|\psi\|_{\mathcal{H}_\eta}.
\end{equation}
$\bullet$ Finally, for the second term, we write 
\begin{align*}
\bigg| \mathrm{Im}\int_{D_\eta} \bigg| \na_v \bigg(\frac{\psi}{M}\bigg)\cdot \na_v\bigg(\frac{\chi_\eta \bar \psi}{M}\bigg) M^2 \ud v \bigg| &= \bigg| \mathrm{Im} \int_{D_\eta}  \chi_\eta'\frac{\psi}{M} \pa_{v_1}\bigg(\frac{\bar\psi}{M}\bigg) M^2 \ud v \bigg| \\
&\leqslant \bigg(\int_{D_\eta} |\chi_\eta' \psi|^2 \ud v \bigg)^{\frac{1}{2}} \left\| \na_v\bigg(\frac{\psi}{M}\bigg)M\right\|_{L^2(\RR^d)} .
\end{align*}
As in the previous point,  since $|\chi_\eta' \psi|^2 \leqslant C \eta |v_1| |\psi|^2 1_{A_\eta}$,  then
\begin{align}\label{J3}
\bigg| \mathrm{Im}\int_{D_\eta} \bigg| \na_v \bigg(\frac{\psi}{M}\bigg)\cdot \na_v\bigg(\frac{\chi_\eta \bar \psi}{M}\bigg) M^2 \ud v \bigg| &\leqslant C_1 \left(I^\eta_1 \right)^{\frac{1}{2}}  \| \psi\|_{\mathcal{H}_0} \nonumber \\
&\leqslant C_2 \left(I^\eta_2 + \| \psi\|_{\mathcal{H}_0}^2  \right)^{\frac{1}{2}}  \| \psi\|_{\mathcal{H}_0} \nonumber \\
&\leqslant \frac{1}{4} I^\eta_2 + C \| \psi\|_{\mathcal{H}_0}^2  .
\end{align}
Where we used inequality \eqref{I^eta_1 < I^eta_2 + |psi|_H_0^2} in the second line and the Young's inequality in the third line. Now, by summing up inequalities \eqref{J1}, \eqref{J2}, and \eqref{J3}, we obtain:
$$  I^\eta_2 \leqslant \int_{D_\eta} \eta v_1 \chi_\eta \psi \bar \psi \ud v \leqslant \frac{1}{4} I^\eta_2 + C \left( \|A_{\lambda,\eta} \psi\|_{ \mathcal{H}_\eta} \|\psi\|_{\mathcal{H}_\eta} + \|\psi\|_{\mathcal{H}_0}^2 + |\lambda| \|\psi\|_{\mathcal{H}_\eta}^2  \right) . $$
Finally,  we have from inequality \eqref{a}, $\|\psi\|_{\mathcal{H}_0}^2 \lesssim \|A_{\lambda,\eta} \psi\|_{ \mathcal{H}_\eta} \|\psi\|_{\mathcal{H}_\eta} +  |\lambda| \|\psi\|_{\mathcal{H}_\eta}^2$. Therefore,
$$
I^\eta_2 \leqslant C \left( \|A_{\lambda,\eta} \psi\|_{ \mathcal{H}_\eta} \|\psi\|_{\mathcal{H}_\eta} + |\lambda| \|\psi\|^2_{\mathcal{H}_\eta} \right) .
$$
Incorporating this last inequality into \eqref{I^eta_1 < I^eta_2 + |psi|_H_0^2},  and using \eqref{a} we obtain that
$$ I^\eta_1 \leqslant C \left( \|A_{\lambda,\eta} \psi\|_{ \mathcal{H}_\eta} \|\psi\|_{\mathcal{H}_\eta} + |\lambda| \|\psi\|^2_{\mathcal{H}_\eta} \right) .$$
Thus, by summing the last two inequalities on $I^\eta_1$ and $I^\eta_2$ with $\|\psi\|_{\mathcal{H}_0}^2$,  we obtain
$$ \|\psi\|_{\mathcal{H}_\eta}^2 \leqslant C \bigg( \|A_{\lambda,\eta} \psi\|_{ \mathcal{H}_\eta} \|\psi\|_{\mathcal{H}_\eta} + |\lambda| \|\psi\|^2_{\mathcal{H}_\eta} \bigg) .$$
Hence inequality \eqref{psi < A psi} holds, thanks to the inequality $ab \leqslant C a^2+\frac{b^2}{4C}$ applied to the term $\|A_{\lambda,\eta} \psi\|_{ \mathcal{H}_\eta} \|\psi\|_{\mathcal{H}_\eta}$, and with $\lambda$ small enough, $|\lambda|\leqslant \frac{1}{4C}$. 
\ep

\begin{lemma}[Complementary lemma]\label{lemme compl.}
Let $\eta, \lambda_0>0$ fixed small enough. Let $\lambda \in \CC$ such that $|\lambda|\leqslant \lambda_0$. Then, for all $\psi, F \in \mathcal{H}_\eta$ such that $|a(\psi,\psi)| \leqslant C \| F \|_{\mathcal{H}_\eta} \| \psi \|_{\mathcal{H}_\eta}$, the following inequality holds
\begin{equation}\label{psi < C F dans H_eta}
\| \psi \|_{\mathcal{H}_\eta} \leqslant \tilde C \| F \|_{\mathcal{H}_\eta} ,
\end{equation}
where $C$ and $\tilde C$ are two positive constants that do not depend on $\lambda$ and $\eta$.
\end{lemma}
\bp 
The proof is identical to that of the previous lemma,  just replace the inequality $|a(\psi,\psi)| \leqslant \|A_{\lambda,\eta} \psi\|_{\mathcal{H}_\eta} \|\psi \|_{\mathcal{H}_\eta}$ by $|a(\psi,\psi)| \leqslant C \| F \|_{\mathcal{H}_\eta} \| \psi \|_{\mathcal{H}_\eta}$.
\ep

Let denote by $\mathcal{H}_\eta'$ the topological dual of $\mathcal{H}_\eta$. By the Riesz representation theorem, for all $F \in \mathcal{H}_\eta'$,  there exists a unique $f \in \mathcal{H}_\eta$ such that 
$$ (F,\phi) = \langle f,\phi\rangle_{\mathcal{H}_\eta} , \quad \forall \phi \in \mathcal{H}_\eta , $$ 
where $(F,\phi)$ denotes the value taken by $F \in \mathcal{H}_\eta'$ in $\phi \in \mathcal{H}_\eta$. Then, by Remark \ref{A representant de a}, the problem  
\begin{equation}\label{formulation fable}
a(\psi,\phi) = (F,\phi) , \quad \forall \phi \in \mathcal{H}_\eta
\end{equation}
is equivalent to the problem $A_{\lambda,\eta} \psi = f$ for $f \in \mathcal{H}_\eta$. Therefore, equivalent to the invertibility of the operator $A_{\lambda,\eta}$.
\begin{proposition}\label{inversibilite de A} Let $\eta_0>0$ and $\lambda_0>0$ small enough. Let $\eta \in [0,\eta_0]$ and $\lambda \in \CC$ fixed, with $|\lambda|\leqslant \lambda_0$. 
For all $F \in \mathcal{H}_\eta'$,  equation \eqref{formulation fable} admits a unique solution $\psi^{\lambda,\eta} \in \mathcal{H}_\eta\subset  \mathcal{H}_0$, satisfying the following estimate
\begin{equation}\label{psi_H-eta < C F_H-eta'}
\| \psi^{\lambda,\eta} \|_{\mathcal{H}_0} \leqslant \| \psi^{\lambda,\eta} \|_{\mathcal{H}_\eta} \leqslant C \| F \|_{\mathcal{H}_\eta'} ,
\end{equation}
where $C$ is a positive constant that does not depend on $\lambda$ and $\eta$. Moreover, for $F \in L^2_{\langle v\rangle^2} \subset \mathcal{H}_\eta'$ we have
\begin{equation}\label{psi_H-eta < C <v>F_L^2}
\| \psi^{\lambda,\eta} \|_{\mathcal{H}_0} \leqslant\| \psi^{\lambda,\eta} \|_{\mathcal{H}_\eta} \leqslant C \| F \|_{L^2_{\langle v\rangle^2}} ,
\end{equation}
where $L^2_{\langle v\rangle^2}$ denote the weighted $L^2$ space: $\ds L^2_{\langle v\rangle^2}:=\left\{f: \RR^d \longrightarrow \CC \ ; \ \int_{\RR^d} |f|^2 \langle v\rangle^2 \mathrm{d}v < \infty \right\}$.
\end{proposition}
\begin{remark}
The sesquilinear form $a$ depends continuously on $\eta$ and holomorphically on $\lambda$. \\
The solution in the previous proposition, is for $\lambda$ and $\eta$ fixed, and it depends on $\lambda$ and $\eta$ since $a$ depends on these last parameters.
\end{remark}

 We denote by $T_{\lambda,\eta}$ the inverse operator of $L_{\lambda,\eta}$ for $\lambda$ and $\eta$ fixed, i.e. the operator which associates to $F$ the solution $\psi^{\lambda,\eta} =: T_{\lambda,\eta}(F)$.
\subsection{Existence of solutions via the implicit function theorem}
The results of this section are the same as those of \cite[Subsection 2.2]{DP-2}. We recall them here for the self-containedness of the paper. The idea is as follows, we use the operator $T_{\lambda,\eta}$ to rewrite equation \eqref{eq penalisee2} as a fixed point problem for the identity plus a compact map. Then, the Fredholm alternative will allow us to apply the implicit function theorem in order to get the existence of solutions.  For this purpose,  we define $F: \{\lambda\in\CC;|\lambda|\leqslant\lambda_0\}\times[0,\eta_0]\times \mathcal{H}_0\longrightarrow \mathcal{H}_0$ by 
$$F(\lambda,\eta,h):= h - \mathcal{T}_{\lambda,\eta}(h) ,$$ 
with
$$\mathcal{T}_{\lambda,\eta}(h):= T_{\lambda,\eta}\big[Vh-\langle h-M,\Phi \rangle \Phi \big] .$$
Note that finding a solution $h(\lambda,\eta)$, solution to $F\big(\lambda,\eta,h(\lambda,\eta)\big)=0$,  gives a solution to the penalized equation by taking $M_{\lambda,\eta}=h(\lambda,\eta)$.\\

Recall that the function $\Phi$ is given by $\Phi := \left(\int_{\RR^d} \frac{M^2}{\langle v \rangle^2} \ud v \right)^{-1} \frac{M}{\langle v \rangle^2}$ and satisfies:
\begin{enumerate}
\item For all $v$ in $\RR^d$,  one has: $\quad \Phi(v) \lesssim \frac{M(v)}{\langle v \rangle^2}$.  
\item The function $ \Phi$ belongs to the weighted Lebesgue space $L^2_{\langle v\rangle^2} := L^2\big(\RR^d,\langle v\rangle^2\ud v\big)$  thanks to Assumption \ref{a1}.
\item Its dot product, in $L^2(\RR^d;\RR)$,  with $M$ is equal to $1$: $\quad \langle \Phi, M \rangle = 1$.
\item Assumption \ref{a2} means that $\langle v \rangle \left|\na_v \Phi\right|$ belongs to the space $ L^2(\RR^d)$. Thus, $\langle v \rangle \Phi$ belongs to the Sobolev space $H^1(\RR^d;\RR)$.
\end{enumerate}
\begin{remark}
Note that the operator $\mathcal{T}_{\lambda,0}$ does not depend on $\lambda$ since $T_{\lambda,0}$ does not. Let's denote it by $\mathcal{T}_0$. Also, $\mathcal{T}_{\lambda,\eta}$ is affine with respect to $h$, we denote by $\mathcal{T}^l_{\lambda,\eta}$ its linear part.
\end{remark}

With the results from the previous section, the operator $\mathcal{T}_{\lambda,\eta}$ is continuous with respect to the function $h \in \mathcal{H}_0$ and with respect to the parameters $\lambda$ and $\eta$. Furthermore, the proof provided in \cite[Lemma 2.15]{DP-2} remains unchanged and does not depend on the expression of the equilibrium $M$ or the potential $W$. We recall the statement in the following lemma. 
\begin{lemma}\label{continuite de T_eta}
Let $\eta_0, \lambda_0>0$ small enough. Let $\eta \in [0,\eta_0]$ and $\lambda \in \mathbb{C}$ such that $|\lambda|\leqslant\lambda_0$. Then,  \\
\ni 1.  The map $\mathcal{T}_{\lambda,\eta} : \mathcal{H}_0 \longrightarrow \mathcal{H}_\eta$ is continuous.  Moreover, there exists a constant $C>0$, independent of $\lambda$ and $\eta$ such that
\begin{equation}\label{continuite de H_0 dans H_eta}
\| \mathcal{T}_{\lambda,\eta}^l(h)\|_{\mathcal{H}_\eta} \leqslant C \| h \|_{\mathcal{H}_0} , \quad \forall h \in \mathcal{H}_0 ,
\end{equation}
and the embedding $ \mathcal{T}_{\lambda,\eta}^l(\mathcal{H}_0) \subset \mathcal{H}_\eta \subset \mathcal{H}_0 $ holds for all $ \eta \in [0,\eta_0]$ and for all $\lambda \in \{|\lambda|\leqslant \lambda_0\}$.  Consequently,  the map $\mathcal{T}_{\lambda,\eta} : \mathcal{H}_0 \longrightarrow \mathcal{H}_0$ is continuous.  \smallskip

\ni 2.  The map $\mathcal{T}_{\lambda,\eta}$ is continuous with respect to $\lambda$ and $\eta$.  Moreover, there exists a constant $C>0$, independent of $\lambda$ and $\eta$ such that, for all $\eta' \in [0,\eta_0]$ and for all $|\lambda'| \leqslant \lambda_0$
\begin{equation}\label{T_eta(h)-T_0(h) -->0}
\left\|\mathcal{T}_{\lambda,\eta}(h)-\mathcal{T}_{\lambda,\eta'}(h)\right\|_{\mathcal{H}_0} \leqslant C \bigg(\bigg|1-\frac{\eta'}{\eta} \bigg|+\bigg|1-\bigg|\frac{\eta'}{\eta}\bigg|^{\frac{2}{3}} \bigg| \bigg) \big( \|h \|_{\mathcal{H}_0} + \left\|\langle v \rangle\Phi \right\|_{L^2}\big)
\end{equation}
and
\begin{equation}
\left\|\mathcal{T}_{\lambda,\eta}(h)-\mathcal{T}_{\lambda',\eta}(h)\right\|_{\mathcal{H}_0} \leqslant C  |\lambda-\lambda'| \big( \|h \|_{\mathcal{H}_0} + \left\|\langle v \rangle\Phi\right\|_{L^2}\big) 
\end{equation}
for all $h \in \mathcal{H}_0$.
\end{lemma}
For the compactness of the operator $\mathcal{T}_{\lambda,\eta}$ at the point $(\lambda,\eta) = (0,0)$, we have 
\begin{lemma}\label{compacite de T_0^l}
Assume \ref{a1} and \ref{a2}. Then, the map $\mathcal{T}^l_0$ is compact.
\end{lemma}
 \bp  We use the fact that we can approximate both functions $ g_1 := \langle v \rangle^2 V $ and $ g_2 := \langle v \rangle \Phi $ by functions with compact support in $ W^{1,\infty} $ and $ H^1 $ respectively, since $ g_1 \in C^1_0(\mathbb{R}^d, \mathbb{R}) $; the space of functions that tend to $0$ at infinity along with their first derivatives; and $ g_2 \in H^1(\mathbb{R}^d, \mathbb{R}) $. Once we are on a compact set, we apply the Rellich theorem. The key point here is that the compact set we work on does not depend on the function $h$.  See \cite[Lemma 2.16]{DP-2} for more details.
\ep

As a consequence of the two previous lemmas, we have the following proposition.
\begin{proposition}[Assumptions of the implicit function theorem \cite{DP-2}]\label{hypotheses de TFI} 
\item 1. The map $F(\lambda,\eta,\cdot)=Id-\mathcal{T}_{\lambda,\eta}$ is continuous on $\mathcal{H}_0$ uniformly with respect to $\lambda$ and $\eta$. Moreover, there exists $c>0$, independent of $\lambda$ and $\eta$ such that
$$ \|F(\lambda,\eta,h_1)-F(\lambda,\eta,h_2)\|_{\mathcal{H}_0} \leqslant c \|h_1-h_2\|_{\mathcal{H}_0}, \qquad \forall h_1, h_2 \in \mathcal{H}_0,  \ \forall\eta\in[0,\eta_0],  \ \forall |\lambda|\leqslant \lambda_0 .$$
2.  The map $F$ is continuous with respect to $\lambda$ and $\eta$ and we have
$$ \underset{\eta \to \eta'}{\lim}\|F(\lambda,\eta,h)-F(\lambda,\eta',h)\|_{\mathcal{H}_0} = \underset{\lambda \to \lambda'}{\lim} \|F(\lambda,\eta,h)-F(\lambda',\eta,h)\|_{\mathcal{H}_0} = 0,  \qquad \forall h \in \mathcal{H}_0.$$
3.  The map $F(\lambda,\eta,\cdot)$ is differentiable on $\mathcal{H}_0$. Moreover,
 $$ \frac{\partial F}{\partial h}(\lambda,\eta,\cdot)=Id-\mathcal{T}_{\lambda,\eta}^l, \qquad \forall |\lambda|\leqslant\lambda_0,\ \forall \eta\in[0,\eta_0] .$$
4.  We have $\ F(0,0,M)=0 \ $ and $\ \frac{\partial F}{\partial h}(0,0,M)$ is invertible.
\end{proposition}

\begin{theorem}[Existence of solutions with constraint]\label{thm d'existence} 
Assume \ref{a1}-\ref{a2}.  Let $\eta_0>0$ and  $\lambda_0>0$ small enough. Let $\eta \in [0,\eta_0]$ and $\lambda \in \CC$ fixed, with $|\lambda|\leqslant \lambda_0$. Then, there is a unique function $M_{\lambda,\eta}$ in $\mathcal{H}_0$, solution to the penalized equation
\begin{equation}\label{eq de M_eta + cont.}
\big[-\Delta_v+ W(v)+  \mathrm{i}  \eta v -\lambda\eta^{\frac{2}{3}}\big]M_{\lambda,\eta}(v)= b(\lambda,\eta)\Phi(v), \qquad v \in \RR^d ,
\end{equation} 
where $\ds b(\lambda,\eta):=\langle N_{\lambda,\eta},\Phi \rangle$ with $N_{\lambda,\eta} := M_{\lambda,\eta}-M$. Moreover,
\begin{equation}\label{M_lambda,eta-M_0-->0 dans H_0}
 \| N_{\lambda,\eta} \|_{\mathcal{H}_0} = \| M_{\lambda,\eta} - M \|_{\mathcal{H}_0} \underset{\eta \rightarrow 0}{\longrightarrow}0.
\end{equation}
\end{theorem}
The proof of this theorem is a direct consequence of the implicit function theorem applied to the function $F$ of the previous proposition, at the point $(\lambda,\eta,h)=(0,0,M)$.
\begin{remark}\label{rmq sur thm d'existence}
The sequence $|b(\lambda,\eta)|$ is uniformly bounded with respect to $\lambda$ and $\eta$ since $ \ds |b(\lambda,\eta)| \underset{\eta \rightarrow 0}{\longrightarrow}0$, which we obtain by the Cauchy-Schwarz inequality and limit \eqref{M_lambda,eta-M_0-->0 dans H_0}:
\begin{equation}\label{b(lambda,eta) --> 0}
|b(\lambda,\eta)|= |\langle N_{\lambda,\eta},\Phi \rangle| \leqslant \left\|\frac{N_{\lambda,\eta}}{\langle v \rangle} \right\|_{L^2} \| \langle v \rangle \Phi\|_{L^2} \leqslant \| N_{\lambda,\eta} \|_{\mathcal{H}_0} \| \langle v \rangle \Phi\|_{L^2} \underset{\eta \rightarrow 0}{\longrightarrow}0 .
\end{equation}
\end{remark}
\section{Existence of the eigen-solution $(\mu(\eta),M_{\mu,\eta})$}\label{section sol-propre}
The aim of this section is to prove Theorem \ref{main}. It is composed of three subsections. In the first, we recall a result on $L^2$ estimates for the function $M_{\lambda,\eta}$. These estimates allow us to apply the implicit function theorem to the penalized term to get the existence of an eigenpair. The second subsection is devoted to the study of the rescaled eigenfunction, which is itself used to investigate the behavior of the eigenvalue in the final subsection.

\subsection{$L^2$ estimates and study of the penalized term}\label{subsection L2}
We start this subsection with the following proposition.

\begin{proposition}[$L^2$ estimates \cite{DP-2}]\label{propostion estimation L^2} Assume \ref{a1}-\ref{a2}.  Let $ \eta_0, \lambda_0>0$ small enough.  Then, for all $\eta \in [0,\eta_0]$ and for all $\lambda \in \CC$ such that $|\lambda|\leqslant \lambda_0$,  the solution $M_{\lambda,\eta}$ of the penalised equation \eqref{eq de M_eta + cont.} satisfies the following estimates \smallskip

\ni 1.  For all $\gamma > \frac{d}{2}$, one has
\begin{equation}\label{estimation de N_lambda,eta dans L^2}
\|N_{\lambda,\eta}\|_{L^2(\RR^d)}^2 := \|M_{\lambda,\eta}-M\|_{L^2(\RR^d)}^2 \lesssim |\lambda| + |b(\lambda,\eta)| + \eta^{\frac{\delta}{3}} ,
\end{equation}
and 
\begin{equation}\label{M_eta sur A^c fin}
\|M_{\lambda,\eta}\|^2_{L^2(\{ |v_1| \geqslant s_0\eta^{-\frac{1}{3}}\})} \leqslant \frac{1}{s_0^2} \left\| \eta^{\frac{1}{3}} v_1 M_{\lambda,\eta} \right\|_{L^2(\{ |v_1| \geqslant s_0\eta^{-\frac{1}{3}}\})} \lesssim \eta^{\frac{\delta}{3}} ,
\end{equation}
where $\delta := \frac{1}{2}(\gamma-\frac{d}{2})$ and $s_0$ is a positive constant.  \smallskip

\ni 2.  For all $\gamma > \frac{d+1}{2}$, one has 
\begin{equation}\label{estimation de v_1^1/2 N_eta dans L^2}
\left\|  |v_1|^{\frac{1}{2}} N_{\lambda,\eta} \right\|_{L^2(\RR^d)}^2 := \left\| |v_1|^{\frac{1}{2}}(M_{\lambda,\eta}-M) \right\|_{L^2(\RR^d)}^2 \lesssim 1 .
\end{equation}
In particular, the functions $M_{\lambda,\eta}$ and $|v_1|^{\frac{1}{2}} M_{\lambda,\eta}$ are uniformly bounded, with respect to $\lambda$ and $\eta$, in $L^2(\RR^d,\CC)$ for $\gamma > \frac{d}{2}$ and $\gamma > \frac{d+1}{2}$ respectively. 
\end{proposition}

\ni \bp The proof of this proposition is given in \cite[Proposition 3.1]{DP-2} and does not depend on the symmetry of the equilibrium $M$, we only need the inequality $\langle v \rangle^{-\gamma} \lesssim M(v) \lesssim \langle v \rangle^{-\gamma}$ in order to recover the Hardy-Poincar\'e inequality and an inequality of the Poincar\'e type proved in \cite[Lemma 3.2]{DP-2}.  Here we only recall the idea of the proof.  We are only discussing the first point, the second is done in the same way.  Let denote $v:=(v_1,v') \in \RR\times\RR^{d-1}$.  The proof is given in four steps and the idea is as follows: first, we decompose $\RR^d$ into two parts, $\RR^d = \{ |v_1|\leqslant s_0\eta^{-\frac{1}{3}}\} \cup \{|v_1| \geqslant s_0\eta^{-\frac{1}{3}}\}$, small/medium and large velocities.  In the first step, using the equation of $M_{\lambda,\eta}$, we estimate the norm of $M_{\lambda,\eta}$ for large velocities to get
$$ \|M_{\lambda,\eta} \|^2_{L^2(\{|v_1| \geqslant s_0\eta^{-\frac{1}{3}}\})}  \leqslant \nu_1 \|M_{\lambda,\eta} \|^2_{L^2(\{ |v_1|\leqslant s_0\eta^{-\frac{1}{3}}\})}  + c_1(\eta)  , $$
where $\nu_1$ and $c_1$ depend on $s_0$, $\lambda$ and $\eta$. To estimate $\|M_{\lambda,\eta} \|_{L^2(\{ |v_1|\leqslant s_0\eta^{-\frac{1}{3}}\})}$, it is enough to estimate $\| N_{\lambda,\eta} \|_{L^2(\{ |v_1|\leqslant s_0\eta^{-\frac{1}{3}}\})}$ since $M$ belongs to $L^2$, which is the purpose of steps two and three. In step 2, using a Poincar\'e type inequality, we show that
$$\|N_{\lambda,\eta} \|^2_{L^2(\{ |v_1| \leqslant s_0\eta^{-\frac{1}{3}}\leqslant |v'|\})}  \leqslant C_1 \|M_{\lambda,\eta} \|^2_{L^2(\{|v_1| \geqslant s_0\eta^{-\frac{1}{3}}\})}  + c_2(\eta)  , $$
where $C_1$ is a positive constant and $c_2$ depends on $s_0$, $\lambda$ and $\eta$. Then, in the third step, using the Hardy-Poincar\'e inequality, we prove that
 $$\|N_{\lambda,\eta} \|^2_{L^2(\{|v|\leqslant s_0\eta^{-\frac{1}{3}}\})}  \leqslant \nu_2 \|N_{\lambda,\eta} \|^2_{L^2(\{ |v_1|\leqslant s_0\eta^{-\frac{1}{3}}\})} + \nu_3 \|M_{\lambda,\eta} \|^2_{L^2(\{|v_1| \geqslant s_0\eta^{-\frac{1}{3}}\})}  + c_3(\eta)  ,$$  
with $\nu_2$, $\nu_3$ and $c_3$ depend on $s_0$, $\lambda$ and $\eta$. The last step is left for the conclusion: we first fix $s_0$ large enough, then $|\lambda|$ small enough, then $\eta$ small enough, we obtain $\nu_2\leqslant \frac{1}{4}$, $\nu_3\leqslant \frac{1}{4}$ and $\nu_1 \big(C_1+\frac{\nu_3}{1-\nu_2}\big) \leqslant \frac{1}{2}$, which allows us to conclude.
\ep

\ni As a direct consequence of Proposition \ref{propostion estimation L^2}, we have the following two results: \\

In the first, we have two limits which prepare the assumptions of the implicit function theorem applied to the \emph{penalized term} $b(\lambda,\eta) := \langle M_{\lambda,\eta} - M,\Phi\rangle$.
\begin{corollary}[\cite{DP-2}]\label{lim int M_eta ...}
Let $M_{\lambda,\eta}$ be the solution to equation \eqref{eq de M_eta + cont.}. Then, for all $\lambda \in \CC$ such that, $|\lambda| \leqslant \lambda_0$ with $\lambda_0$ small enough, the following limit holds:
\begin{equation}\label{lim int eta^(1/3) v_1 M_eta M}
\underset{\eta \to 0}{\lim } \int_{\RR^d} \eta^{\frac{1}{3}} v_1 M_{\lambda,\eta}(v)M(v) \mathrm{d}v = 0 .
\end{equation}
For $\lambda =0$, one has
\begin{equation}\label{lim int  M_eta M}
\underset{\eta \to 0}{\lim } \int_{\RR^d} M_{0,\eta}(v)M(v) \mathrm{d}v = \int_{\RR^d} M^2(v) \mathrm{d}v .
\end{equation}
\end{corollary}
\bp For the first limit, we write with $\delta := \frac{1}{2}(\gamma - \frac{d}{2}) \in (0,1)$ for $\gamma \in (\frac{d}{2},\frac{d+4}{2})$:
\begin{align*}
\bigg| \int_{\RR^d} \eta^{\frac{1}{3}} v_1 M_{\lambda,\eta}(v) M(v) \mathrm{d}v \bigg| &\leqslant \eta^{\frac{1}{3}} \int_{\{|v_1| \leqslant s_0 \eta^{-\frac{1}{3}}\}}  |v_1|^{1-\delta} |v_1|^\delta M(v) \left|M_{\lambda,\eta}(v)\right| \mathrm{d}v \\
& \hspace{1.5cm} + \int_{\{|v_1| \geqslant s_0 \eta^{-\frac{1}{3}}\}} |v_1|^{-\delta} |v_1|^\delta M(v) \left|\eta^{\frac{1}{3}} v_1 M_{\lambda,\eta}(v)\right| \mathrm{d}v \\
&\lesssim_{s_0} \eta^{\frac{\delta}{3}} \left\| |v_1|^\delta M \right\|_2 \bigg( \| M_{\lambda,\eta} \|_{L^2(\RR^d)} + \left\| \eta^{\frac{1}{3}} v_1 M_{\lambda,\eta} \right\|_{L^2(\{|v_1| \geqslant s_0 \eta^{-\frac{1}{3}}\})} \bigg) \\
&\lesssim \eta^{\frac{\delta}{3}} \underset{\eta \to 0}{\longrightarrow} 0 ,
\end{align*}
thanks to \eqref{estimation de N_lambda,eta dans L^2} and \eqref{M_eta sur A^c fin}. \\
The second limit is a direct consequence of inequality \eqref{estimation de N_lambda,eta dans L^2} and  limit \eqref{b(lambda,eta) --> 0}.
\ep

In the second one, we have the existence of a $\mu$, a function of $\eta$, such that the constraint $\langle M_{\mu(\eta),\eta} - M,\Phi\rangle =0$ is satisfied. This gives us a pair of $(\mu(\eta),M_{\mu(\eta),\eta})$ solutions to the spectral problem \eqref{M_mu,eta}.

\begin{proposition}[Cancellation of the penalized term \cite{DP-2}]\label{contrainte}
Define $\  B(\lambda,\eta):= \eta^{-\frac{2}{3}} b(\lambda,\eta) $.  

\ni 1. The expression of $B(\lambda,\eta)$ is given by
\begin{equation}\label{b(lambda,eta)}
B(\lambda,\eta)= \eta^{-\frac{2}{3}}\langle N_{\lambda,\eta},\Phi\rangle = \int_{\RR^d}(\lambda- \mathrm{i} \eta^{\frac{1}{3}} v)M_{\lambda,\eta}(v)M(v)\mathrm{d}v .
\end{equation}
2.  The $\eta$ order of $B(\lambda,\eta)$ in its expansion with respect to $\lambda$ is given by
\begin{equation}\label{B(lambda,0)}
\underset{\eta \rightarrow 0}{\lim} \ \frac{\pa B}{\pa \lambda}(0,\eta) = \int_{\RR^d} M^2(v)\mathrm{d}v .
\end{equation}
3.  There exists $\tilde\eta_0, \tilde\lambda_0>0$ small enough, a function $\tilde{\lambda}: \{|\eta|\leqslant\tilde\eta_0\} \longrightarrow \{|\lambda|\leqslant\tilde\lambda_0\}$ such that,
for all $(\lambda,\eta)\in [0,\tilde\eta_0[\times\{|\lambda|<\tilde\lambda_0\}$,  $\lambda = \tilde{\lambda}(\eta)$ and the constraint is satisfied: $$B(\lambda,\eta)= B(\tilde{\lambda}(\eta),\eta)=0.$$
Consequently, $\mu(\eta) := \eta^{\frac{2}{3}} \tilde{\lambda} (\eta)$ is the eigenvalue associated to the eigenfunction $M_\eta: = M_{\tilde{\lambda}(\eta),\eta}$ for the operator $\mathcal{L}_\eta$, and the couple $\big(\mu(\eta),M_\eta \big)$ is solution to the spectral problem \eqref{M_mu,eta}.
\end{proposition} 
\bp 1. We obtain equation \eqref{b(lambda,eta)} by integrating the penalized equation \eqref{eq de M_eta + cont.} against $M$ and using the facts that $Q(M) = 0$ and $\langle M, \Phi \rangle_{L^2} = 1$. \\
2. Since
$$ B(\lambda, \eta) = \lambda \langle M_{\lambda, \eta}, M \rangle - \ui \eta^\frac{1}{3} \langle v_1 M_{\lambda, \eta}, M \rangle ,  $$
we have
$$ \frac{\partial B}{\partial \lambda}(0, \eta) = \langle M_{0, \eta}, M \rangle,  $$
and the limit as $\eta$ tends to $0$ is given by limit \eqref{lim int  M_eta M} from the previous corollary. \\
3. The existence of the function $\tilde{\lambda}$ is a direct consequence of the implicit function theorem applied to $B$ at the point $(\lambda, \eta) = (0, 0)$.  Note that $ B(0,0) = 0$ thanks to  limit \eqref{lim int Im M_0,eta/<v> v_1 petit}.
\ep

\subsection{Preliminary estimates for the rescalated solution}\label{subsection H_eta}
The objective of this subsection is to prepare the necessary arguments to study the approximation of the eigenvalue and to compute the diffusion coefficient. To achieve this, by writing $\mu(\eta) = \eta^\alpha \kappa_\eta$, we will take the limit as $\eta \to 0$ and prove that $\kappa_\eta$ converges to a positive constant. For this purpose, we will use a weak-strong convergence argument in the Hilbert space $L^2$ to compute the limit of $\kappa_\eta$. Specifically, we interpret $\kappa_\eta$ as $\langle \mathsf{H}_\eta , \mathsf{I}_\eta \rangle$,  where $\mathsf{I}_\eta \rightarrow \mathsf{I}_0$ in $L^2$ and $\mathsf{H}_\eta \rightharpoonup \mathsf{H}_0$. To establish this convergence, we first prove that $\|\mathsf{H}_\eta\|_{L^2}$ is uniformly bounded with respect to $\eta$, ensuring convergence up to a subsequence.  Subsequently, to establish convergence of the entire sequence, we identify the limit $\mathsf{H}_0$ and demonstrate its uniqueness. This subsection is devoted to estimates of $\mathsf{H}_\eta$.  \\

We will first start by introducing the rescaled function of $M_{0,\eta}$ as well as the equation satisfied by this function, and set some notations to avoid long expressions. Recall that $M_{0,\eta}$ satisfies the equation:
$$ \big[ Q + \mathrm{i} \eta v_1 \big] M_{0,\eta}(v) = - b(0,\eta) \Phi(v) ,  \quad v \in \RR^d ,$$
with $Q = -\frac{1}{M}\nabla_v \cdot \big(M^2\nabla_v\big(\frac{\cdot}{M}\big)\big)$ and $b(0,\eta) = \langle M_{0,\eta}-M,\Phi \rangle$. 
Then, the rescaled function $H_\eta$ defined by $H_\eta(s):=\eta^{-\frac{\gamma}{3}}M_{0,\eta}(\eta^{-\frac{1}{3}}s)$ is solution to the rescaled equation
\begin{equation}\label{eq H_eta}
\big[ Q_\eta + \mathrm{i} s_1 \big] H_\eta(s) = - \eta^{-\frac{\gamma+2}{3}} b(0,\eta) \Phi_\eta(s) , \quad s \in \RR^d ,
\end{equation}
where 
$$ Q_\eta := -\frac{1}{m_\eta}\nabla_s \cdot \left(m_\eta^2\nabla_s\bigg(\frac{\cdot}{m_\eta}\bigg)\right)  ,  \qquad  m_\eta(s) :=\eta^{-\frac{\gamma}{3}}M(\eta^{-\frac{1}{3}}s) $$
and 
\begin{equation}\label{Phi_eta < ...}
\Phi_\eta(s) :=\Phi(\eta^{-\frac{1}{3}}s) \lesssim \eta^{\frac{\gamma+2}{3}} |s|_\eta^{-\gamma-2}  \leqslant  \eta^{\frac{\gamma+2}{3}} |s|^{-\gamma-2} ,
\end{equation}
Note that: $Q(M)=0$ implies that $Q_\eta(m_\eta)=0$. \\

Now, we will show a series of lemmas.  The first concerns the estimation of the penalized term $b(0,\eta) = \langle M_{0,\eta}-M,\Phi\rangle$.

\begin{lemma}\label{estimation c_eta -1} 
Let $\nu \in (0,\frac{1}{2})$ and let $\delta := 2\nu(\gamma- \frac{d}{2})$. Then, the following estimate holds:
$$ |b(0,\eta)| = |c_\eta - 1| \lesssim \left\{\begin{array}{l}  \eta^{\frac{2+\delta}{3}},  \quad \mbox{ if } \quad \gamma\in (\frac{d}{2},\frac{d+1}{2}),  \\
 \\
  \eta^{\frac{5-\nu}{6}} ,  \quad   \mbox{ if } \quad \gamma = \frac{d+1}{2},  \\
 \\
 \eta ,   \hspace{0.97cm}  \mbox{ if } \quad  \gamma\in (\frac{d+1}{2},\frac{d+4}{2}) ,
\end{array}\right.  
$$
where 
$$ c_\eta := \left(\int_{\RR^d} \frac{M^2}{\langle v\rangle^2}\ud v\right)^{-1} \int_{\RR^d} \frac{M_{0,\eta}M}{\langle v\rangle^2}\ud v  .  $$
\end{lemma}
\bp  This lemma is a corollary of Proposition \ref{propostion estimation L^2}.  Indeed,  we have
$$ c_\eta - 1 = \left(\int_{\RR^d} \frac{M^2}{\langle v\rangle^2}\ud v\right)^{-1} \int_{\RR^d} \frac{(M_{0,\eta}-M)M}{\langle v\rangle^2}\ud v = \langle N_{0,\eta},\Phi \rangle = b(0,\eta).$$
Then,  by \eqref{b(lambda,eta)} we get
$$ |c_\eta - 1| = \eta \left| \int_{\RR^d} v_1 M_{0,\eta}(v)M(v) \ud v\right| \leqslant \eta \int_{\{|v_1|\leqslant \eta^{-\frac{1}{3}}\}} \big| v_1 M_{0,\eta} M \big| \ud v +  \eta\int_{\{|v_1| \geqslant \eta^{-\frac{1}{3}}\}} \big| v_1 M_{0,\eta}M \big|\ud v.$$
\textbf{Case 1: $\gamma\in (\frac{d}{2},\frac{d+1}{2})$.} For this range of $\gamma$, we write
\begin{align*}
 |c_\eta - 1| &\leqslant \eta \int_{\{|v_1|\leqslant \eta^{-\frac{1}{3}}\}} |v_1|^{1-\delta} |M_{0,\eta}| |v_1|^{\delta}M  \ud v +  \eta^{\frac{2}{3}}  \int_{\{|v_1| \geqslant \eta^{-\frac{1}{3}}\}} \left|\eta^{\frac{1}{3}} v_1 M_{0,\eta}\right| \eta^{\frac{\delta}{3}}|v_1|^{\delta}M \ud v \\
 &\leqslant  \eta^{\frac{2+\delta}{3}} \big\| |v_1|^{\delta}M \big\|_{L^2(\RR^d)} \left( \|M_{0,\eta}\|_{L^2(\{|v_1|\leqslant \eta^{-\frac{1}{3}}\})} + \big\| \eta^{\frac{1}{3}} v_1 M_{0,\eta}\big\|_{L^2(\{|v_1| \geqslant \eta^{-\frac{1}{3}}\})} \right) \\
 &\lesssim \eta^{\frac{2+\delta}{3}},
\end{align*}
thanks to the first item of Proposition \ref{propostion estimation L^2}. \\

\ni \textbf{Case 2: $\gamma = \frac{d+1}{2}$.} Similarly, we have 
\begin{align*}
 |c_\eta - 1| &\leqslant \eta \int_{\{|v_1|\leqslant \eta^{-\frac{1}{3}}\}} |v_1|^{\frac{1+\nu}{2}} |M_{0,\eta}| |v_1|^{\frac{1-\nu}{2}}M  \ud v + \eta^{\frac{2}{3}} \int_{\{|v_1| \geqslant \eta^{-\frac{1}{3}}\}} |\eta^{\frac{1}{3}} v_1 M_{0,\eta}| \eta^{\frac{1-\nu}{6}}|v_1|^{\frac{1-\nu}{2}}M \ud v \\
 &\leqslant   \eta^{\frac{5-\nu}{6}}  \big\| |v_1|^{\frac{1-\nu}{2}}M \big\|_{L^2(\RR^d)} \left( \|M_{0,\eta}\|_{L^2(\RR^d)}  + \big\| \eta^{\frac{1}{3}} v_1 M_{0,\eta}\big\|_{L^2(\{|v_1| \geqslant \eta^{-\frac{1}{3}}\})} \right) \\
 &\lesssim  \eta^{\frac{5-\nu}{6}} ,
\end{align*}
thanks to the first item of Proposition \ref{propostion estimation L^2} and since $|v_1|^{\frac{1-\nu}{2}}M \in L^2(\RR^d)$ for $\gamma = \frac{d+1}{2}$ and $\nu > 0$.  \\
\ni \textbf{Case 3: $\gamma\in (\frac{d+1}{2},\frac{d+4}{2})$.} By the second item of Proposition \ref{propostion estimation L^2},  $|v_1|^{\frac{1}{2}} |M_{0,\eta}|$ is uniformly bounded in $L^2(\RR^d)$ for all $\gamma>\frac{d+1}{2}$. Therefore, 
$$
 |c_\eta - 1| \leqslant \eta  \left\| |v_1|^{\frac{1}{2}}M_{0,\eta} \right\|_{L^2(\RR^d)} \left\| |v_1|^{\frac{1}{2}}M \right\|_{L^2(\RR^d)}  \lesssim  \eta .
$$
\ep

 The second one show that the (very) small velocities do not participate in the limit of the diffusion coefficient.
\begin{lemma}[Small velocities]\label{petites vitesses}
\item 1.  For all $\gamma \in(\frac{d+1}{2},\frac{d+4}{2})$,  the following estimate holds
\begin{equation}\label{Im M_0,eta/<v> v_1 petit} \int_{\{|v_1|\leqslant R\}}  \bigg|\frac{  M_{0,\eta} - c_\eta M}{\langle v \rangle}\bigg|^2 \mathrm{d}v \lesssim \eta . 
\end{equation}
2.  For all $\gamma \in (\frac{d}{2},\frac{d+1}{2})$, one has
\begin{equation}\label{lim int Im M_0,eta/<v> v_1 petit}
\underset{\eta \to 0}{\lim} \ \eta^{1-\alpha}  \int_{\{|v_1|\leqslant R\}} v_1 M_{0,\eta}(v)M(v) \ \mathrm{d}v = 0 .
\end{equation}
3.  For all $\gamma \in [\frac{d+1}{2},\frac{d+4}{2})$, one has
\begin{equation}\label{lim int Im M_0,eta/<v> v_1 petit}
\underset{\eta \to 0}{\lim} \ \eta^{1-\alpha}  \int_{\{|v_1|\leqslant R\}} v_1 [M_{0,\eta}(v)- c_\eta M(v)]M(v) \ \mathrm{d}v = 0 .
\end{equation}
\end{lemma}

\ni \bp \textbf{1.} 
The function $M_{0,\eta} - c_\eta M$ satisfies condition \eqref{condition d'orthogonalite}. Then, by the Hardy-Poincar\'e inequality \eqref{Hardy-Poincare avec orthogonalite},  there exists a positive constant $C_{\gamma,d}$ such that: $$ \int_{\{|v_1|\leqslant R\}}  \bigg|\frac{M_{0,\eta}- c_\eta M}{\langle v \rangle}\bigg|^2 \mathrm{d}v \leqslant \bigg\|\frac{M_{0,\eta} - c_\eta M}{\langle v \rangle}\bigg\|^2_{L^2(\RR^d)} \leqslant C_{\gamma,d} \bigg\|\nabla_v\bigg(\frac{M_{0,\eta}}{M}\bigg)M\bigg\|^2_{L^2(\RR^d)} .$$
Now,  we have on the one hand
$$ 
\bigg\|\nabla_v\bigg(\frac{M_{0,\eta}}{M}\bigg)M\bigg\|^2_{L^2(\RR^d)} = \Re\bigg( \int_{\RR^d} Q(M_{0,\eta})\overline M_{0,\eta} \mathrm{d}v \bigg) .
$$
On the other hand,
\begin{align*}
\Re\bigg( \int_{\RR^d} Q(M_{0,\eta}) \bar M_{0,\eta} \ud v \bigg) &= \Re \bigg[ -\ui \eta  \int_{\RR^d} v_1  |M_{0,\eta}|^2 \ud v - b(0,\eta)  \int_{\RR^d} \Phi \bar M_{0,\eta} \ud v \bigg] \\
& = \Re \bigg[ - b(0,\eta)  \int_{\RR^d} \Phi \bar M_{0,\eta} \ud v \bigg] \\
&=  \Re \bigg[ \ui \eta \int_{\RR^d} v_1  M_{0,\eta} M \ud v  \int_{\RR^d} \Phi \bar M_{0,\eta} \ud v \bigg]  . 
\end{align*}
Which implies that, 
$$
 \bigg|\Re\bigg( \int_{\RR^d} Q(M_{0,\eta}) \bar M_{0,\eta} \ud v \bigg)  \bigg|  \leqslant \eta \big\| |v_1|^\frac{1}{2} M_{0,\eta} \big\|_{L^2}  \big\| |v_1|^\frac{1}{2} M \big\|_{L^2} \| M_{0,\eta} \|_{L^2} \| \Phi \|_{L^2} \lesssim \eta .
$$
Hence \eqref{Im M_0,eta/<v> v_1 petit} holds true, since  $\big\|(1+ |v_1|^\frac{1}{2}) M_{0,\eta}\big\|_{L^2} \lesssim 1$ for $\gamma > \frac{d+1}{2}$ by Proposition \ref{propostion estimation L^2}.\\

\ni \textbf{2.  Case 1: $\gamma \in (\frac{d}{2},\frac{d+1}{2})$.}  We have
$$
 \eta^{1-\alpha}\bigg| \int_{\{|v_1|\leqslant R\}}  v_1   M_{0,\eta}(v)M(v) \mathrm{d}v \bigg| \leqslant  R \eta^{1-\alpha} \| M_{0,\eta} \|_{L^2} \| M \|_{L^2}  \lesssim \eta^{1-\alpha} \underset{\eta \to 0}{\longrightarrow} 0,
$$
since $1-\alpha = \frac{1+d-2\gamma}{3} > 0$ for all $\gamma < \frac{d+1}{2}$,  and $\| M_{0,\eta} \|_{L^2} \lesssim 1$ by Proposition \ref{propostion estimation L^2}.  \\

 \textbf{Case 2: $\gamma = \frac{d+1}{2}$.} We have $1-\alpha = 0$.  By writing $M_{0,\eta} - c_\eta M = N_{0,\eta} + (c_\eta - 1)M$, we get  
$$ \bigg| \int_{\{|v_1|\leqslant R\}} v_1 [M_{0,\eta}- c_\eta M]M\  \mathrm{d}v \bigg|   \leqslant R\| M \|_{L^2} \big( \| N_{0,\eta} \|_{L^2} + |c_\eta-1| \| M \|_{L^2} \big) \underset{\eta \to 0}{\longrightarrow} 0 ,$$ 
since $\| N_{0,\eta} \|_{L^2}  \underset{\eta \to 0}{\longrightarrow} 0$ and $ |c_\eta-1|  \underset{\eta \to 0}{\longrightarrow} 0$. \\

\ni \textbf{3.} Let $\gamma \in (\frac{d+1}{2},\frac{d+4}{2})$.  Similary, we have by Cauchy-Schwarz,
\begin{align*}
\eta^{1-\alpha} \bigg| \int_{\{|v_1|\leqslant R\}} v_1 [M_{0,\eta}- c_\eta M]M\  \mathrm{d}v \bigg|  & \leqslant \eta^{1-\alpha}  \big\| v_1 \langle v \rangle M \big\|_{L^2(\{|v_1|\leqslant R\})} \bigg\| \frac{M_{0,\eta}-c_\eta M}{\langle v \rangle } \bigg\|_{L^2(\{|v_1|\leqslant R\})} \\
&\lesssim \eta^{2-\alpha} \underset{\eta \to 0}{\longrightarrow} 0,
\end{align*}
thanks to inequality \eqref{Im M_0,eta/<v> v_1 petit} and since $\alpha < 2$ for $\gamma < \frac{d+4}{2}$,  and $v_1 \langle v \rangle M \in L^2(\{|v_1|\leqslant R\})$ for $\gamma \in (\frac{d+1}{2},\frac{d+4}{2})$.
\ep

The third lemma contains some important estimates on the rescaled solution $H_\eta$. 

\begin{lemma}\label{grandes vitesses} Let $s_0>0$ fixed, large enough. Then, the following uniform estimates, with respect to $\eta$,  hold: \smallskip

\ni 1.  For all $\gamma \in (\frac{d}{2},\frac{d+1}{2}]$ and for $\nu \in (\gamma-\frac{d+1}{2},\gamma-\frac{d-1}{2})$,  one has
\begin{equation}\label{estimation gamma<(d+1)/2}
\left\| |s_1|^{\frac{1}{2}}|s|_\eta^{-\nu} (H_\eta - c_\eta m_\eta) \right\|_{L^2(\{|s_1| \leqslant s_0\})} + \big\| s_1 H_\eta \big\|_{L^2(\{|s_1| \geqslant s_0\})} \lesssim 1 ,
\end{equation}
2.  For all $\gamma \in (\frac{d+1}{2},\frac{d+4}{2})$,  one has 
\begin{equation}\label{estimation pour gamma>(d+1)/2}
\left\|\frac{H_\eta - c_\eta m_\eta}{|s|_\eta} \right\|_{L^2(\RR^d)} + \big\| s_1 H_\eta \big\|_{L^2(\{|s_1| \geqslant s_0\})} \lesssim 1 .
\end{equation}
\end{lemma}
Before giving the proof of this lemma, we have the following remark:

\begin{remark}\label{remarque c_eta}
Note that, by making the change of variables $v=\eta^{-\frac{1}{3}}s$ with $\eta>0$,
$$ c_\eta =  \bigg(\int_{\RR^d} \frac{m_\eta^2}{|s|_\eta^2} \mathrm{d}s \bigg)^{-1}\int_{\RR^d}   \frac{m_\eta(s) H_\eta(s)}{|s|_\eta^2} \mathrm{d}s .$$
Hence,  $H_\eta - c_\eta m_\eta$ satisfies the orthogonality condition \eqref{condition d'orthogonalite} of the Hardy-Poincar\'e inequality 
$$ \int_{\RR^d} \frac{m_\eta(s) \big[H_\eta(s) - c_\eta m_\eta(s)\big]}{|s|_\eta^2} \ \mathrm{d}s =0 .$$  
\end{remark} 

\noindent \bpl \ref{grandes vitesses}. We will establish estimates on different ranges of (rescalated) velocities, and in order to avoid long expressions in the proof, we will fix some notations of sets.  Let denote $s:=(s_1,s') \in \RR\times\RR^{d-1}$.  For $s_0>0$, we set: $A :=\{|s_1|\leqslant s_0\}$, $B:=\{|s_1|\leqslant s_0 , |s'|\leqslant s_0\}$ and $C:=\{|s_1|\leqslant s_0 \leqslant |s'|\}$.  Also, for $\eta>0$, we denote by $K_\eta$ and $\tilde K_\eta$ the functions defined by $K_\eta := H_\eta - m_\eta$ and  $\tilde K_\eta := H_\eta - c_\eta m_\eta$ respectively.  \\

\noindent \textbf{1.} Let $\gamma \in (\frac{d}{2},\frac{d+1}{2}]$ and let $\nu \in (\gamma-\frac{d+1}{2},\gamma-\frac{d-1}{2})$.  To prove the first point of this lemma, we will estimate $\big\| |s_1|^{\frac{1}{2}} |s|_\eta^{-\nu} \tilde K_\eta \big\|_{L^2(B)}$ using the Hardy-Poincar\'e inequality, $\big\| |s_1|^{\frac{1}{2}} |s|_\eta^{-\nu} \tilde K_\eta \big\|_{L^2(C)}$ using the Poincar\'e-type inequality \eqref{Poincare2 ineq} after the change of variable $v=\eta^{-\frac{1}{3}}s$, and estimate $\|s_1 H_\eta\|_{L^2(A^c)}$ using the equation for $H_\eta$. Thus, we obtain inequality \eqref{estimation gamma<(d+1)/2} by combining these estimates, and since $|s_1|^{\frac{1}{2}} |s|_\eta^{\nu} m_\eta(s) \lesssim |s_1|^{\frac{1}{2}} |s|^{\nu-\gamma} \in L^2(A)$ for $\gamma \in (\frac{d}{2},\frac{d+1}{2}]$ and $\nu \in (\gamma-\frac{d+1}{2},\gamma-\frac{d-1}{2})$.   \\

\noindent \textbf{Estimation of $\big\| |s_1|^{\frac{1}{2}} |s|_\eta^{-\nu} \tilde K_\eta \big\|_{L^2(B)}$.} Recall that $B :=\{|s| \leqslant s_0\}$. On the one hand, we have thanks to the Hardy-Poincar\'e inequality \eqref{Hardy-Poincare avec orthogonalite} after the change of variable $v=\eta^{-\frac{1}{3}}s$, 
$$
 \left\| \frac{\tilde K_\eta}{|s|_\eta} \right\|^2_{L^2(\RR^d)}  \lesssim_{\gamma,d}  \bigg\| \nabla_s \bigg(\frac{\tilde K_\eta}{m_\eta}\bigg) m_\eta \bigg\|^2_{L^2(\RR^d)} .
$$
Therefore,
\begin{equation}\label{s_1^(1/2)K_eta sur B 1}
 \big\| |s_1|^{\frac{1}{2}} |s|_\eta^{-\nu} \tilde K_\eta \big\|^2_{L^2(B)}  \lesssim s_0^{3-2\nu} \left\| \frac{\tilde K_\eta}{|s|_\eta} \right\|^2_{L^2(B)}  \lesssim s_0^{3-2\nu}  \bigg\| \nabla_s \bigg(\frac{\tilde K_\eta}{m_\eta}\bigg)m_\eta \bigg\|^2_{L^2(\RR^d)}   .
\end{equation}
On the other hand,  since $\nabla_s \big(\frac{\tilde K_\eta}{m_\eta}\big) = \nabla_s \big(\frac{ K_\eta}{m_\eta}\big)$ then, by integrating the equation of $K_\eta$ against $\bar K_\eta$ and taking the real part, and using the fact that $-b(0,\eta)\int_{\RR^d} \Phi_\eta\bar K_\eta \ud s = -\eta^{\frac{d-\gamma}{3}} |b(0,\eta)|^2 \leqslant 0$, we write:
$$  \bigg\| \nabla_s \bigg(\frac{\tilde K_\eta}{m_\eta}\bigg)m_\eta \bigg\|^2_{L^2(\RR^d)}  =  \Re \int_{\RR^d} Q_\eta(K_\eta)\bar K_\eta \ \mathrm{d}s \leqslant \bigg| \Re \int_{\RR^d} \mathrm{i} s_1 m_\eta \bar K_\eta \ \mathrm{d}s  \bigg| =  \bigg|  \int_{\RR^d} s_1 m_\eta \mathrm{Im} K_\eta \ \mathrm{d}s\bigg| .$$
Thus, by splitting the integral above into two parts, on $A :=\{|s_1|\leqslant s_0\}$ and on $A^c$, and by writing $K_\eta = \tilde K_\eta + (c_\eta-1)m_\eta$ on $A$, we obtain:
\begin{align} \label{nabla K/m}
 \bigg\| \nabla_s \bigg(\frac{\tilde K_\eta}{m_\eta}\bigg) m_\eta \bigg\|^2_{L^2(\RR^d)} \leqslant \left\| |s_1|^\frac{1}{2} |s|_\eta^{\nu} m_\eta \right\|_{L^2(A)} & \left\| |s_1|^\frac{1}{2} |s|_\eta^{-\nu} \tilde K_\eta \right\|_{L^2(A)} + \left\| m_\eta \right\|_{L^2(A^c)} \left\| s_1 H_\eta \right\|_{L^2(A^c)} \nonumber \\
  &+ |c_\eta -1| \int_A |s_1|  \ m_\eta^2(s) \ \ud s .
\end{align}
Hence, returning to \eqref{s_1^(1/2)K_eta sur B 1}, we get:
\begin{align}\label{s_1^(1/2)K_eta sur B}
 \big\| |s_1|^{\frac{1}{2}} |s|_\eta^{-\nu} \tilde K_\eta \big\|^2_{L^2(B)} &\leqslant C s_0^{3-2\nu} \bigg\| \nabla_s \bigg(\frac{\tilde K_\eta}{m_\eta}\bigg) m_\eta \bigg\|^2_{L^2(\RR^d)} \nonumber \\
 &\leqslant  \frac{1}{4} \big\| |s_1|^\frac{1}{2} |s|_\eta^{-\nu} \tilde K_\eta \big\|^2_{L^2(A)} + C s_0^{3-2\nu} \big\| m_\eta \big\|_{L^2(A^c)} \big\| s_1 H_\eta \big\|_{L^2(A^c)} + C_{s_0}  ,
\end{align}
where $C_{s_0} >0$  depends only on $s_0$,  since for $\gamma \in (\frac{d}{2},\frac{d+1}{2}]$ and $\nu \in (\gamma-\frac{d+1}{2},\gamma-\frac{d-1}{2})$,
$$ \int_A  |s_1| \ |s|_\eta^{2\nu} \ m_\eta^2(s) \ \ud s  \lesssim  \int_A |s_1| |s|^{2\nu-2\gamma} \ \ud s  \lesssim 1 ,  $$
and by Lemma \ref{estimation c_eta -1},  for $\nu \in (\gamma-\frac{d+1}{2},\gamma-\frac{d-1}{2}) \subset (0,\frac{1}{2})$, we write
\begin{equation} \label{c_eta-1 int_A}
|c_\eta -1| \int_A  |s_1| \ m_\eta^2(s) \ud s \lesssim \eta^{\min(\frac{2+\delta}{3},\frac{5-\nu}{6})} \int_A \eta^{-\frac{2\nu}{3}} |s_1| \ |s|_\eta^{2\nu} \ m_\eta^2(s) \ \ud s   \lesssim \eta^{\omega} ,
\end{equation}
with $\omega := \min\left(\frac{2+\delta}{3},\frac{5-\nu}{6}\right) - \frac{2\nu}{3} = \min\left(\frac{2(1-\nu)+\delta)}{3},\frac{5(1-\nu)}{6}\right) > 0$ for $\gamma \in (\frac{d}{2},\frac{d+1}{2}]$, $\delta := 2\nu(\gamma-\frac{d}{2})$. \\

\noindent \textbf{Estimation of $\big\| |s_1|^{\frac{1}{2}} |s|_\eta^{-\nu} \tilde K_\eta \big\|_{L^2(C)}$.} Recall that $C:=\{|s_1|\leqslant s_0\leqslant |s'|\}$.  Let's begin by estimating $\big\|\tilde K_\eta \big\|^2_{L^2(C)}$.  Following the proof of Lemma \ref{Poincare2} up to inequality \eqref{P2}, with the change of variable $v=\eta^{-\frac{1}{3}}s$ and $\chi\left(\frac{s_1}{s_0}\right)$ instead of $\chi(\eta^\frac{1}{3}v_1)$, we obtain
$$ \big\| \tilde K_\eta \big\|^2_{L^2(C)} \lesssim s_0^2 \bigg\| \nabla_s \bigg(\frac{\tilde K_\eta}{m_\eta}\bigg) m_\eta \bigg\|^2_{L^2(\RR^d)} + \frac{1}{s_0^2} \int_{A^c} \big|s_1 H_\eta \big|^2 \ud s + \int_{A^c} m_\eta^2(s) \ud s ,  $$
where we used the inequality $|\tilde K_\eta| \leqslant |H_\eta| + m_\eta$ in the integral over $A^c$.  Now, since 
$$\big\| |s_1|^{\frac{1}{2}} |s|_\eta^{-\nu} \tilde K_\eta \big\|_{L^2(C)}^2 \leqslant s_0^{1-2\nu} \big\| \tilde K_\eta \big\|^2_{L^2(C)} $$ 
then, using inequality \eqref{s_1^(1/2)K_eta sur B} with Young's inequality for  $\left\| m_\eta \right\|_{L^2(A^c)} \left\| s_1 H_\eta \right\|_{L^2(A^c)}$, we obtain:
\begin{align}\label{s_1^(1/2)K_eta sur C}
\left\| |s_1|^{\frac{1}{2}} |s|_\eta^{-\nu} \tilde K_\eta \right\|_{L^2(C)}^2 &\lesssim s_0^{3-2\nu} \bigg\| \nabla_s \bigg(\frac{\tilde K_\eta}{m_\eta}\bigg) m_\eta \bigg\|^2_{L^2(\RR^d)} + \frac{1}{s_0^{1+2\nu}} \int_{A^c} \big|s_1 H_\eta \big|^2 \ud s + s_0^{1-2\nu} \int_{A^c} m_\eta^2(s) \ud s \nonumber \\
 &\leqslant  \frac{1}{4} \left\| |s_1|^{\frac{1}{2}} |s|_\eta^{-\nu} \tilde K_\eta \right\|^2_{L^2(A)} + \frac{1}{ s_0} \left\| s_1 H_\eta \right\|_{L^2(A^c)}^2 + C_{s_0} ,
\end{align}
where $C_{s_0}$ is a positive constant that depends only on $s_0$,  since $\int_{A^c} m_\eta^2(s) \ud s$ is uniformly bounded with respect to $\eta$.  \\

\noindent \textbf{Estimation of $\big\| s_1 H_\eta \big\|_{L^2(A^c)}$.}  Let $\chi \in C^\infty(\RR)$ such that $0 \leqslant \chi \leqslant 1$, $\chi \equiv 0$ on $B(0,\frac{1}{2})$ and $\chi \equiv 1$ outside $B(0,1)$. We set $\chi_{s_0}(s_1):= \chi\left( \frac{s_1}{s_0}\right)$. Then, we have $\big\| s_1 H_\eta \big\|_{L^2(A^c)} \leqslant \big\| s_1 \chi_{s_0} H_\eta \big\|_{L^2(D)}$, where $D:=\{|s_1|\geqslant \frac{s_0}{2}\}$. Now, by integrating the equation of $H_\eta$ against $ s_1 \bar H_\eta \chi_{s_0}^2$ and taking the imaginary part, we obtain:
\begin{equation}\label{int_D s_1 zeta H_eta}
\int_D \big| s_1 \chi_{s_0} H_\eta \big|^2 \mathrm{d}s = -\mathrm{Im} \bigg(\int_D Q_\eta(H_\eta) s_1 \bar H_\eta \chi_{s_0}^2 \mathrm{d}s\bigg) - \eta^{-\frac{\gamma+2}{3}}\mathrm{Im} \bigg( b(0,\eta) \int_D \Phi_\eta s_1 \bar H_\eta \chi_{s_0}^2 \mathrm{d}s\bigg) 
\end{equation}
Let's start with the second term which is simpler. By \eqref{Phi_eta < ...} we have,
\begin{align*}
\eta^{-\frac{\gamma+2}{3}}\bigg|\mathrm{Im} \bigg( b(0,\eta) \int_D \Phi_\eta s_1 \bar H_\eta \chi_{s_0}^2 \mathrm{d}s\bigg) \bigg| &\lesssim \frac{1}{s_0^2} |b(0,\eta)| \left\| |s|^{-\gamma} \right\|_{L^2(D)} \left\| s_1 \chi_{s_0} H_\eta \right\|_{L^2(D)}  \\
&\lesssim \frac{1}{s_0^2} |b(0,\eta)| \left(\left\| |s|^{-\gamma} \right\|_{L^2(D)}^2 + \left\| s_1 \chi_{s_0} H_\eta \right\|^2_{L^2(D)} \right) .
\end{align*}
For the first term,  by integration by parts,  we write:
\begin{align*}
\bigg| \mathrm{Im} \int_D Q_\eta(H_\eta) s_1 \bar H_\eta \chi_{s_0}^2 \mathrm{d}s \bigg| &= \bigg| \mathrm{Im} \int_D \pa_{s_1}\bigg(\frac{H_\eta}{m_\eta} \bigg)m_\eta\chi_{s_0}\big[\chi_{s_0}\bar H_\eta + 2 s_1 \bar H_\eta \chi_{s_0}'\big] \mathrm{d}s \bigg|  \\
&\leqslant \bigg\| \pa_{s_1}\bigg(\frac{H_\eta}{m_\eta} \bigg)m_\eta \bigg\|_{L^2(D)} \bigg( \| \chi_{s_0} H_\eta\|_{L^2(D)} + 2 \big\| s_1 \chi_{s_0}' \chi_{s_0} H_\eta \big\|_{L^2(D)} \bigg) .
\end{align*}
Then, since $\chi_{s_0}'\equiv 0$ except on: $D \setminus A^c:=\{\frac{s_0}{2} \leqslant |s_1|\leqslant s_0\} \subset A$ where $\underset{D}{\sup} \left| s_1\chi_{s_0}'(s_1)\right| = \underset{[\frac{1}{2},1]}{\sup} | t\chi '(t)|\lesssim 1$, we get
\begin{align*}
\bigg| \mathrm{Im} \int_D Q_\eta(H_\eta) s_1 \bar H_\eta \chi_{s_0}^2 \mathrm{d}s \bigg| &\leqslant \frac{2}{s_0} \bigg\| \na_s \bigg(\frac{H_\eta}{m_\eta} \bigg)m_\eta \bigg\|_{L^2(\RR^d)} \left\| s_1 \chi_{s_0} H_\eta \right\|_{L^2(D)} \bigg( 1 + \left\| s_1 \chi_{s_0}' \right\|_{L^\infty(D)} \bigg)   \\
&\leqslant    \frac{C}{s_0} \bigg(\bigg\| \na_s \bigg(\frac{H_\eta}{m_\eta} \bigg)m_\eta \bigg\|_{L^2(\RR^d)}^2 + \left\| s_1 \chi_{s_0}' H_\eta \right\|^2_{L^2(D)} \bigg)  .
\end{align*}
Thus, using the previous estimates with $|b(0,\eta)| \lesssim 1$, returning to \eqref{int_D s_1 zeta H_eta} we obtain:
$$
\left\| s_1 \chi_{s_0} H_\eta \right\|^2_{L^2(D)} \leqslant   \frac{C}{s_0} \left\| s_1 \chi_{s_0} H_\eta \right\|^2_{L^2(D)} + \frac{C}{s_0} \bigg\| \na_s \bigg(\frac{H_\eta}{m_\eta} \bigg)m_\eta \bigg\|_{L^2(\RR^d)}^2 + C. 
$$
Hence,  for $s_0$ large enough, 
\begin{equation}\label{s_1 H_eta sur A^c 1}  
\left\| s_1  H_\eta \right\|^2_{L^2(A^c)} \leqslant \left\| s_1 \chi_{s_0} H_\eta \right\|^2_{L^2(D)} \leqslant    \frac{C}{s_0} \bigg\| \na_s \bigg(\frac{H_\eta}{m_\eta} \bigg)m_\eta \bigg\|_{L^2(\RR^d)}^2 + C.  
\end{equation}
Therefore, as for the two previous steps, by inequality \eqref{nabla K/m} combined with Young's inequality and  \eqref{c_eta-1 int_A}, we write:
\begin{align*}
\frac{1}{s_0} \bigg\| \na_s \bigg(\frac{H_\eta}{m_\eta} \bigg)m_\eta \bigg\|_{L^2(\RR^d)}^2  \lesssim   \frac{1}{s_0^2}  &\left\| |s_1|^{\frac{1}{2}} |s|_\eta^{-\nu} \tilde K_\eta \right\|^2_{L^2(A)} +  \frac{1}{s_0^2} \left\| s_1 H_\eta \right\|^2_{L^2(A^c)}  \\
&+  \left\| |s_1|^{\frac{1}{2}} |s|^{\nu-\gamma} \right\|^2_{L^2(A)} + \left\| |s|^{-\gamma} \right\|^2_{L^2(A^c)} + \eta^\omega .
\end{align*}
This gives, returning to \eqref{s_1 H_eta sur A^c 1},  and for $s_0$ sufficiently large 
\begin{equation}\label{s_1 H_eta sur A^c fin}
\left\| s_1  H_\eta \right\|^2_{L^2(A^c)}  \lesssim  \frac{1}{s_0^2}   \left\| |s_1|^{\frac{1}{2}} |s|_\eta^{-\nu} \tilde K_\eta \right\|^2_{L^2(A)} + 1 .
\end{equation}
Finally, summing \eqref{s_1^(1/2)K_eta sur B} and \eqref{s_1^(1/2)K_eta sur C}, and using the previous inequality, we obtain:
$$ \left\| |s_1|^{\frac{1}{2}} |s|_\eta^{-\nu} \tilde K_\eta \right\|^2_{L^2(A)}   \lesssim 1 .$$
Hence,  inequality \eqref{estimation gamma<(d+1)/2} follows from the previous inequality with \eqref{s_1 H_eta sur A^c fin}.

\begin{remark} Note that we can take $s_0$ sufficiently large such that all the previous inequalities hold; more precisely, all the norms present next to $\frac{1}{s_0}$ can be absorbed.
\end{remark}

\noindent \textbf{2.} Let $\gamma \in (\frac{d+1}{2},\frac{d+4}{2})$.  Recall that $\tilde K_\eta:=H_\eta-c_\eta m_\eta$ satisfies the orthogonality condition \eqref{condition d'orthogonalite} of Hardy-Poincar\'e inequality \eqref{Hardy-Poincare avec orthogonalite}. Then, we get on the one hand
\begin{equation}\label{tilde K_eta / |s|}
\bigg\| \frac{\tilde K_\eta}{|s|_\eta} \bigg\|^2_{L^2(A)}   \leqslant \bigg\| \frac{\tilde K_\eta}{|s|_\eta} \bigg\|^2_{L^2(\RR^d)}   \lesssim_{\gamma,d} \bigg\| \nabla_s\bigg(\frac{\tilde K_\eta}{m_\eta}\bigg) m_\eta \bigg\|^2_{L^2(\RR^d)} = \bigg\| \nabla_s\bigg(\frac{K_\eta}{m_\eta}\bigg) m_\eta \bigg\|^2_{L^2(\RR^d)} .
\end{equation}
On the other hand, proceeding as for obtaining inequality \eqref{s_1^(1/2)K_eta sur B}, we write:
\begin{align*}
 \bigg\| \nabla_s \bigg(\frac{K_\eta}{m_\eta}\bigg)m_\eta \bigg\|^2_{L^2(\RR^d)}  &=  \Re \int_{\RR^d} Q_\eta(K_\eta)\bar K_\eta \ \mathrm{d}s \leqslant \bigg| \Re \int_{\RR^d} \mathrm{i} s_1 m_\eta \bar K_\eta \ \mathrm{d}s  \bigg| =  \bigg|  \int_{\RR^d} s_1 m_\eta \mathrm{Im} K_\eta \ \mathrm{d}s\bigg|  \\ 
&\leqslant  \int_{A} \left| s_1 m_\eta \mathrm{Im} \tilde K_\eta\right| \ud s + \int_{A^c} \left| s_1 m_\eta \mathrm{Im} H_\eta\right| \ud s  + |c_\eta -1| \int_A |s_1|  \ m_\eta^2(s) \ \ud s .
\end{align*}
Hence,
\begin{align}\label{grad K_eta pour 2gamma>(d+1)}
 \bigg\| \nabla_s \bigg(\frac{K_\eta}{m_\eta}\bigg)m_\eta \bigg\|^2_{L^2(\RR^d)}  &\leqslant \left\| s_1 |s|_\eta m_\eta \right\|_{L^2(A)} \bigg\| \frac{\tilde K_\eta}{|s|_\eta} \bigg\|_{L^2(A)} + s_0^{-\frac{1}{2}} \left\| |s_1|^{\frac{1}{2}} m_\eta \right\|_{L^2(A^c)} \left\|s_1 H_\eta \right\|_{L^2(A^c)} \nonumber \\
& \hspace*{3cm} + |c_\eta -1| \int_A |s_1|  \ m_\eta^2(s) \ \ud s.
\end{align}
It remains to estimate the norm $\|s_1 H_\eta\|_{L^2(A^c)}$,  and then we simply return to \eqref{tilde K_eta / |s|} and use Young's inequality to complete the estimates as in the first point.  For this purpose, by revisiting the proof from point 1 concerning the estimation of $\|s_1 H_\eta\|_{L^2(A^c)}$ up to inequality \eqref{s_1 H_eta sur A^c 1}, since the proof does not depend on the phase of $\gamma$, we obtain:
\begin{equation}\label{s_1 H gamma>}
\left\| s_1  H_\eta \right\|^2_{L^2(A^c)}  \leqslant    \frac{C}{s_0} \bigg\| \na_s \bigg(\frac{H_\eta}{m_\eta} \bigg)m_\eta \bigg\|_{L^2(\RR^d)}^2 + C.  
\end{equation}
Also, by Lemma \ref{estimation c_eta -1},  and since $m_\eta(s) \lesssim |s|_\eta^{-\gamma}$ and $\gamma \in (\frac{d+1}{2},\frac{d+4}{2})$, we have
$$|c_\eta -1| \int_A |s_1|  \ m_\eta^2(s) \ \ud s \lesssim \eta  \int_A |s_1|  |s|_\eta^{-2\gamma} \ \ud s \lesssim \int_A |s_1|  |s|_\eta^{3-2\gamma} \ \ud s \lesssim 1 .$$
Therefore, by inserting inequality \eqref{s_1 H gamma>} into \eqref{grad K_eta pour 2gamma>(d+1)}, we obtain
\begin{align*}
 \bigg\| \nabla_s \bigg(\frac{K_\eta}{m_\eta}\bigg)m_\eta \bigg\|^2_{L^2(\RR^d)}  
&\leqslant \left\| s_1 |s|_\eta m_\eta \right\|_{L^2(A)} \bigg\| \frac{\tilde K_\eta}{|s|_\eta} \bigg\|_{L^2} + \frac{C}{s_0} \left\| |s_1|^{\frac{1}{2}} m_\eta \right\|_{L^2(A^c)}  \bigg\| \nabla_s \bigg(\frac{\tilde K_\eta}{m_\eta}\bigg)m_\eta \bigg\|_{L^2} + C \\
&\leqslant \left\| s_1 |s|_\eta m_\eta \right\|_{L^2(A)} \bigg\| \frac{\tilde K_\eta}{|s|_\eta} \bigg\|_{L^2(A)} + \frac{1}{2} \bigg\| \nabla_s \bigg(\frac{\tilde K_\eta}{m_\eta}\bigg)m_\eta \bigg\|_{L^2(\RR^d)}^2 + C ,
\end{align*}
since $\big\| |s_1|^{\frac{1}{2}} m_\eta \big\|_{L^2(A^c)} \lesssim \big\| |s_1|^{\frac{1}{2}} |s|^{-\gamma} \big\|_{L^2(A^c)} \lesssim 1$ for $\gamma > \frac{d+1}{2}$.  Hence, 
\begin{equation}\label{grad K_eta fin}
\bigg\| \nabla_s \bigg(\frac{K_\eta}{m_\eta}\bigg)m_\eta \bigg\|^2_{L^2(\RR^d)}  
\leqslant C\left\| s_1 |s|_\eta m_\eta \right\|_{L^2(A)} \bigg\| \frac{\tilde K_\eta}{|s|_\eta} \bigg\|_{L^2(\RR^d)} + C . 
\end{equation}
Finally, by using Young's inequality and returning to \eqref{tilde K_eta / |s|}, we obtain:
\begin{equation}\label{tilde K_eta / |s| fin}
 \bigg\| \frac{\tilde K_\eta}{|s|_\eta} \bigg\|_{L^2(\RR^d)}^2 \lesssim \ \left\| s_1 |s|_\eta m_\eta \right\|_{L^2(A)}^2 + 1 \lesssim 1, 
\end{equation}
since $\big\| s_1 |s|_\eta m_\eta \big\|_{L^2(A)} \lesssim \big\| s_1 |s|^{1-\gamma} \big\|_{L^2(A)} \lesssim 1$ for $\gamma \in (\frac{d+1}{2},\frac{d+4}{2})$.  Thus, inequality \eqref{estimation pour gamma>(d+1)/2} follows from \eqref{tilde K_eta / |s| fin}, \eqref{grad K_eta fin}, and \eqref{s_1 H gamma>}.
\epl

As a direct consequence of Lemma \ref{grandes vitesses}, we have the following additional uniform estimates for all $\gamma \in \left(\frac{d}{2}, \frac{d+4}{2}\right)$.

\begin{lemma}[Complementary estimates] There exists a constant $C > 0$ such that, for all $\eta >0$ small enough,  and for all $\gamma \in \left(\frac{d}{2}, \frac{d+4}{2}\right)$, the following uniform estimates hold:
 \begin{equation}\label{estimations complementaires}
 \bigg\| \frac{H_\eta - c_\eta m_\eta}{|s|_\eta} \bigg\|_{L^2(\RR^d)} +  \bigg\| \nabla_s \bigg(\frac{H_\eta}{m_\eta}\bigg)m_\eta \bigg\|_{L^2(\RR^d)} +  \left\|s_1 H_\eta \right\|_{L^2(\{|s_1|\geqslant 1\})}  \leqslant C .
 \end{equation}
\end{lemma} 
\bp We have from the Hardy-Poincar\'e inequality:
$$ \bigg\| \frac{H_\eta - c_\eta m_\eta}{|s|_\eta} \bigg\|_{L^2(\RR^d)} \lesssim_{\gamma,d}  \bigg\| \nabla_s \bigg(\frac{H_\eta}{m_\eta}\bigg)m_\eta \bigg\|_{L^2(\RR^d)} .$$
For $\gamma \in \left( \frac{d}{2}, \frac{d+1}{2} \right]$, inequality \eqref{estimations complementaires} follows from inequality \eqref{nabla K/m} together with inequality \eqref{estimation gamma<(d+1)/2} from Lemma \ref{grandes vitesses}.  While for $\gamma \in \left(\frac{d+1}{2}, \frac{d+4}{2} \right)$,  it follows from inequality \eqref{grad K_eta fin} together with the second point of Lemma \ref{grandes vitesses}, specifically inequality \eqref{estimation pour gamma>(d+1)/2}.  For the norm $\|s_1 H_\eta\|_{L^2(\{|s_1|\geqslant s_0\})}$, we can actually revisit the estimates on $\|s_1 \chi_{s_0} H_\eta\|_{L^2(D)}$ and show that, by using Young's inequality instead of taking $s_0$ sufficiently large, for any $s_0 > 0$ we have
$$
\left\| s_1 \chi_{s_0} H_\eta \right\|^2_{L^2(D)} \leqslant \frac{1}{2} \left\| s_1 \chi_{s_0} H_\eta \right\|^2_{L^2(D)} + \frac{C}{s_0^2} \left\| \nabla_s \left(\frac{H_\eta}{m_\eta} \right)m_\eta \right\|_{L^2(\mathbb{R}^d)}^2 + \frac{C}{s_0^4}|b(0,\eta)| \left\| |s|^{-\gamma} \right\|_{L^2(D)}^2.
$$
This implies that $\|s_1 H_\eta\|_{L^2(\{|s_1|\geqslant s_0\})} \lesssim 1$ for any $s_0 > 0$.
\ep

We will complete the study of the rescaled solution $H_\eta$ and establish the link between the eigenvalue and the diffusion coefficient in the following subsection.
 
\subsection{Approximation of the eigenvalue and proof of Theorem \ref{main}}\label{subsection vp}
The aim of this subsection is to give an approximation for the eigenvalue $\mu(\eta)$ given in Proposition \ref{contrainte}. The study of this limit is based on the estimates obtained in the previous subsection, on $M_{0,\eta}$ the solution of the penalized equation \eqref{eq penalisee2} for $\lambda=0$, as well as the solution $H_\eta$ of the rescaled equation.  We have the following
\begin{proposition}[Approximation of the eigenvalue]\label{vp} 
\item Assume \ref{a1}--\ref{a5}. Let $\alpha:=\frac{2\gamma-d+2}{3}$ for all $\gamma \in (\frac{d}{2},\frac{d+4}{2})$.  Then, the eigenvalue $\mu(\eta)$ satisfies
\begin{equation}\label{mu(eta)}
\mu(\eta) - \ui \eta j^\eps_1 = \kappa \eta^{\alpha}\big(1+O(\eta^\alpha)\big) ,
\end{equation}
where we recall that 
\begin{equation}\label{j_eta}
j^\eps_1 := \left\{\begin{array}{l}  0,   \hspace{3.9cm} \mbox{ if } \quad \ds 2\gamma\in \left(d,d+1\right),  \\
 \\
 \ds \int_{\{|v|  \leqslant \eps^{-\frac{1}{3}}\}} v_1 M^2(v)  \ud v , \quad  \ \ \mbox{if } \quad 2\gamma = d+1 ,  \\
 \\
\ds \int_{\RR^d} v_1 M^2(v)  \ud v ,   \hspace{1.6cm}  \mbox{ if } \quad  2\gamma\in \left(d+1,d+4\right)  ,
\end{array}\right.  
\end{equation}
and $\kappa$ is a positive constant given by
\begin{equation}\label{kappa}
\kappa :=  - \int_{\{|s_1| >0\}} s_1 m(s)  \mathrm{Im} H_0(s) \ \mathrm{d}s ,
\end{equation}
and where $H_0$ is the unique solution to  
\begin{equation}\label{eq de H_0}
- \frac{1}{m} \na_s \cdot \bigg(m^2 \ \na_s\bigg(\frac{H_0}{m}\bigg)\bigg) + \mathrm{i} s_1 H_0 = 0 \quad \mbox{ in } \  \mathcal{D}'(\RR^d\setminus\{0\})  ,
\end{equation}
satisfying \begin{equation}\label{condition H_0}
\int_{\{|s_1|\geqslant 1\}}|H_0(s)|^2\mathrm{d}s <+\infty \qquad  \mbox{ and } \qquad   H_0(s)\underset{0}{\sim} m(s).
\end{equation}
\end{proposition}
\begin{remark}
Note that the existence of solutions for equation \eqref{eq de H_0} is obtained by passing to the limit in the rescaled equation \eqref{eq H_eta}, while the uniqueness is obtained by an integration by part on $\RR^d\setminus \{0\}$,  using the two conditions of \eqref{condition H_0}.
\end{remark}

To make the proof of this proposition shorter, we will first present a few lemmas. The first one concerns the convergence of the rescaled equilibrium $m_\eta$ in $H^1$ away from the origin.

\begin{lemma}\label{lemma cv m_eta}
Assume \ref{a1}, \ref{a3} and \ref{a5}. Let $\Omega$ be a compact set of $\mathbb{R}^d \setminus \{0\}$. Then, the rescaled equilibrium $m_\eta$ converges, when $\eta$ goes to $0$, to $m$ in $H^1(\Omega)$. The function $m$ is the equilibrium of the limiting operator.
\end{lemma}

\ni \bp The function $m_\eta$ satisfies the equation
$$
-\Delta_s m_\eta + W_\eta(s) \ m_\eta = 0 ,
$$
where $W_\eta$ is the rescaled potential defined by $W_\eta(s) := \eta^{-\frac{2}{3}} W(\eta^{-\frac{1}{3}} s)$.  Note that  \ref{a5} implies that $\Delta_s m_\eta$ is uniformly bounded in $L^2(\{|s|\geqslant r\})$ for all $r > 0$,  since $|W_\eta  m_\eta(s)| \lesssim |s|^{-\sigma-\gamma}$ with $\sigma \geqslant 2$ and $\gamma>\frac{d}{2}$. Which ensures weak convergence up to a subsequence. With Assumptions \ref{a1} and \ref{a3}, which in particular imply the convergence of $m_\eta$ to $m$ in the latter space, we obtain the weak convergence of $\Delta_s m_\eta$ to $\Delta_s m$ in $L^2(\{|s|\geqslant r\})$ for all $r > 0$. Then, the function $n_\eta := m_\eta - m$ satisfies the equation
\begin{equation}\label{eq n_eta}
-\Delta_s n_\eta = \Delta_s m - W_\eta m_\eta .
\end{equation}
Let $r, R > 0$ be such that $\Omega \subset\subset C := B(0,R) \setminus B(0,r)$, and let $\zeta \in C^\infty(\RR^d)$ be such that $0 \leqslant \zeta \leqslant 1$, $\zeta \equiv 1$ in $\Omega$, and $\zeta \equiv 0$ in $C^c$. By integrating equation \eqref{eq n_eta} against $n_\eta \zeta$, we get:
$$
\int_{C} \left| \nabla_s n_\eta \right|^2 \zeta \ \ud s + \frac{1}{2} \int_{C} \nabla_s |n_\eta|^2 \cdot \nabla_s \zeta \, \ud s = \int_{C} \Delta_s m n_\eta \zeta \ \ud s - \int_{C} W_\eta m_\eta n_\eta \zeta \, \ud s .
$$
Therefore,
$$
\int_{C} \left| \nabla_s n_\eta \right|^2 \zeta \, \ud s \leqslant \bigg(\| \Delta_s \zeta \|_{L^\infty(C \setminus \Omega)} \| n_\eta \|_{L^2(C)} + \| \Delta_s m \|_{L^2(C)} + \| W_\eta m_\eta \|_{L^2(C)} \bigg) \| n_\eta \|_{L^2(C)} .
$$
Hence,   $\ds \int_{\Omega} \left| \nabla_s n_\eta \right|^2  \ud s \underset{\eta \to 0}{\longrightarrow} 0$,  
since $\ds \int_C |n_\eta|^2 \ud s \underset{\eta \to 0}{\longrightarrow} 0$ thanks to   \ref{a3} which yields \ref{m_eta -->m L^2}.
\ep

The second lemma concerns the convergence of the rescaled solution $H_\eta$ to $H_0$, the unique solution of the limiting PDE that satisfies conditions \eqref{condition H_0}.

\begin{lemma}\label{lemma cv H_eta}
Under the assumptions of Lemma \ref{lemma cv m_eta}, the function $H_\eta$,  solution of the equation
$$ -\frac{1}{m_\eta} \na_s \cdot \bigg(m_\eta^2\na_s\bigg(\frac{H_\eta}{m_\eta}\bigg)\bigg) + \ui s_1 H_\eta = - \eta^{-\frac{\gamma+2}{3}} b(0,\eta) \Phi_\eta $$
converges, in $\mathcal{D}'(\RR^d\setminus\{0\})$, to $H_0$, the unique solution of the equation
$$ -\frac{1}{m} \na_s \cdot \bigg(m^2\na_s\bigg(\frac{H_0}{m}\bigg)\bigg) + \ui s_1 H_0 = 0 $$
that satisfies both conditions: $H_0 \in L^2(\{|s| \geqslant 1\})\ $ and $\ H_0(s) \underset{0}{\sim} m(s)$.
\end{lemma}

\ni \bp Let $r>0$.  From inequality \eqref{estimations complementaires}, we obtain that $\na_s\left(\frac{H_\eta}{m_\eta}\right) m_\eta \rightharpoonup \mathbb{H}$ in $L^2(\RR^d;\CC^d)$ and $s_1 H_\eta \rightharpoonup h$ in $L^2(\{|s|\geqslant r\})$. Therefore, by the uniqueness of the limit in $\mathcal{D}'(\RR^d\setminus\{0\})$ and since $m_\eta$ converges to $m$ in $L^2(\{|s|\geqslant r\})$, we deduce that $h = s_1 H_0$ and $\mathbb{H} = \na_s\left(\frac{H_0}{m}\right) m$, where $H_0 := \underset{\eta \to 0}{\lim} H_\eta$. Now, let $\varphi \in \mathcal{D}(\RR^d \setminus \{0\})$. Then, by integrating the equation of $H_\eta$ against $\varphi$, we obtain:
\begin{align*}
\int_{\RR^d \setminus \{0\}} \bigg[ \na_s\bigg(\frac{H_\eta}{m_\eta}\bigg) \cdot \na_s\bigg(\frac{\varphi}{m_\eta}\bigg) m_\eta^2 +\mathrm{i} s_1 \varphi(s)  H_\eta(s) \bigg]\mathrm{d}s &= - \eta^{-\frac{2+\gamma}{3}} b(0,\eta)  \int_{\RR^d \setminus \{0\}} \Phi_\eta(s) \varphi (s)\mathrm{d}s .
\end{align*}
Passing to the limit when $\eta$ goes to $0$ in the last equality, we obtain that $H_\eta$ converges to $H_0$ solution to the equation
\begin{equation}\label{equation H_0}
-\frac{1}{m}\na_s \cdot \bigg(m^2\ \na_s\bigg(\frac{H_0}{m}\bigg)\bigg)+ \ui s_1 H_0 = 0 .
\end{equation}
Indeed,  for the right-hand side of the equality, it converges to $0$ since $\Phi_\eta(s) \lesssim \eta^{\frac {2+\gamma}{3}}|s|^{-2-\gamma}$ by inequality \eqref{Phi_eta < ...},  $b(0,\eta)\to 0$ by inequality \eqref{b(lambda,eta) --> 0},  and we have:
$$\eta^{-\frac{2+\gamma}{3}} |b(0,\eta)| \int_{\RR^d \setminus \{0\}} \left|\Phi_\eta(s) \varphi (s)\right| \,\mathrm{d}s \lesssim |b(0,\eta)| \int_{\RR^d \setminus \{0\}} |s|^{-\gamma-2} \left| \varphi(s) \right| \,\mathrm{d}s \lesssim |b(0,\eta)|.$$
For the left-hand side of the equality, we write:
$$ \int_{\RR^d \setminus \{0\}} \na_s\bigg(\frac{H_\eta}{m_\eta}\bigg) \cdot \na_s \bigg(\frac{\varphi}{m_\eta}\bigg) m_\eta^2 \ \ud s = \int_{\RR^d \setminus \{0\}} \na_s\bigg(\frac{H_\eta}{m_\eta}\bigg) m_\eta \cdot \bigg[\na_s \varphi - \na_s m_\eta \frac{\varphi}{m_\eta} \bigg] \ud s . $$
Thus, we pass to the limit in the weak formulation above using Lemma \ref{lemma cv m_eta} and the weak-strong convergence argument in the Hilbert space $L^2$. This yields
$$ \int_{\RR^d \setminus \{0\}} \bigg[ \na_s\bigg(\frac{H_0}{m}\bigg) \cdot \na_s\bigg(\frac{\varphi}{m}\bigg) m^2 +\mathrm{i} s_1 \varphi(s)  H_0(s) \bigg]\mathrm{d}s = 0 .$$
Moreover, for all $\gamma \in (\frac{d} {2},\frac{d+4}{2})$,  the function $H_\eta$ satisfies the uniform estimate \eqref{estimations complementaires}:
$$ \bigg\| \frac{H_\eta - c_\eta m_\eta}{|s|_\eta}\bigg\|_{L^2(\RR^d)} + \left\|s_1 H_\eta \right\|_{L^2(\{|s_1|\geqslant 1\})} \lesssim 1 .$$
Therefore,  $H_0$ satisfies the estimate
$$ \bigg\| \frac{H_0 - m}{|s|}\bigg\|_{L^2(\RR^d)}  + \left\|s_1 H_0\right\|_{L^2(\{|s_1|\geqslant 1\})} \lesssim 1 . $$
Now,  $\|s_1 H_0\|_{L^2(\{|s_1|\geqslant 1\})} \lesssim 1$ implies that $H_0 \in L^2(\{|s_1|\geqslant 1\})$,  and $\big\| \frac{H_0 - m}{|s|}\big\|_{L^2} \lesssim 1 $ implies that $H_0(s) \underset{0}{ \sim} m(s)$,  since 
$$ \int_{\RR^d} \bigg|\frac{H_0(s) - m(s)}{|s|}\bigg|^2 \ud s =  \int_{\RR^d} \bigg|\frac{H_0(s)}{m(s)}-1\bigg|^2 \frac{m^2(s)}{|s|^2} \ \ud s \lesssim 1 , $$
and $C_1 |s|^{-\gamma} \leqslant m(s) \leqslant C_2 |s|^{-\gamma}$.  So,  a different behaviour near zero would make the latter integral infinite.  These two conditions imply that $H_0$ is the unique solution of equation \eqref{equation H_0}.  This can be seen by integrating the equation satisfied by the difference of two solutions against their complex conjugate over $\RR^d \setminus \{0\}$ and performing integration by parts, using the fact that both solutions are in $L^2$ for large $|s|$ and are equivalent to $m$ as $ |s| \to 0$.  Finally, thanks to the uniqueness of this limit,  the whole sequence $H_\eta$ converges weakly.
\ep

The last lemma in this subsection gives the formula of the diffusion coefficient.

\begin{lemma}\label{lim kappa_eta}
Assume \ref{a1}--\ref{a5}.  Let $\alpha := \frac{2\gamma-d+2}{3}$ for all $\gamma \in (\frac{d}{2},\frac{d+4}{2})$, and let $j^\eps_1$ be given by \eqref{j_eta}. We define $\kappa_\eta$ by
$$ \kappa_\eta := - \eta^{-\alpha} \left[ b(0,\eta) + \ui \eta j^\eps_1 \right] .$$
Then, we have the following limit
$$
\underset{\eta \rightarrow 0}{\lim}\  \kappa_\eta = \ui \left\{\begin{array}{l} \ds \int_{|s_1|>0} s_1 m(s)  H_0(s) \ud s  ,  \hspace{5.9cm} \mbox{if } \  2\gamma < d+1 ,\\
 \\
 \ds \int_{0<|s| \leqslant 1} s_1 \big[H_0(s)-m(s) \big]m(s)  \ud s + \int_{|s| \geqslant 1} s_1 H_0(s)m(s)  \ud s  ,  \ \mbox{if } \ 2\gamma = d+1 , \\
 \\
\ds  \int_{|s_1| > 0} s_1 \big[H_0(s)-m(s) \big]m(s)  \ud s  ,  \hspace{4.35cm} \mbox{if } \ 2\gamma > d+1 ,
\end{array}\right. 
$$
where $H_0$ is the unique solution to \eqref{eq de H_0} satisfying  \eqref{condition H_0}.
\end{lemma}
\begin{remark}
The different formulas for $\kappa$ come from different expressions that define the drift $j^\eps_1$, but we will show in point 3 of Subsection \ref{subsection comments} that the previous limit is real, which gives us a unified formula for all $\gamma \in (\frac{d}{2}, \frac{d+4}{2})$:
$$ \kappa :=  \underset{\eta \rightarrow 0}{\lim}\  \kappa_\eta = - \int_{\{|s_1|>0\}} s_1 m(s) \Imm H_0(s) \ \ud s .$$
\end{remark}

\noindent \bpl \ref{lim kappa_eta}.  First, recall that $b(0,\eta)$ is given by
$$ b(0,\eta) = - \ui \eta \int_{\RR^d} v_1 M_{0,\eta}(v) M(v) \, \ud v .$$
We will decompose the integral (or the difference of integrals, depending on $ j^\eps_1$) defining $ \kappa_\eta$ into three parts as follow: \\
$\bullet$ For $2\gamma \neq d+1$, we write:
$$ \kappa_\eta = \ui \eta^{1-\alpha} \int_{\RR^d} v_1 \left( M_\eta(v) - c_\eta 1_{\{2\gamma > d+1\}} M(v) \right) M(v) \ud v +  \ui \eta^{1-\alpha} (c_\eta -1) 1_{\{2\gamma > d+1\}} \int_{\RR^d} v_1 M^2 \ \ud v .$$
The term on the right tends to $0$ thanks to Lemma \ref{estimation c_eta -1}.  For the integral on the left, we decompose it into  $|v_1|\leqslant R$,  $R \leqslant |v_1| \leqslant \eta^{-\frac{1}{3}} $,  and $|v_1| \geqslant \eta^{-\frac{1}{3}} $.  For the first part, which corresponds to very small velocities, it tends to $0$ as $\eta \to 0$ thanks to Lemma \ref{petites vitesses}.  For the other velocities, medium and large, we proceed with the change of variable $v = \eta^{-\frac{1}{3}} s$, which leads to the calculation of a limit of an integral that can be seen as the inner product between two sequences, and whose limit is computed using the weak-strong convergence argument in $ L^2$ as mentioned earlier.   \\
$\bullet$ For $2\gamma = d+1$, we have $\alpha =1$. Then, we write:
$$ \kappa_\eta = \ui \int_{\{|v| \leqslant \eps^{-\frac{1}{3}}\}} v_1 \left( M_\eta - c_\eta M \right) M \ud v +  \ui  (c_\eta -1) \int_{\{|v| \leqslant \eps^{-\frac{1}{3}}\}} v_1 M^2 \ud v + \ui \int_{\{|v| \geqslant \eps^{-\frac{1}{3}}\}} v_1 M_\eta M  \ud v .$$
Then, we proceed exactly as in the previous case. The middle term tends to $0$ thanks to Lemma \ref{estimation c_eta -1}, and for the other two terms, we perform the same change of variable and decompose the integrals over the same sets as before. More precisely, in both cases, this amounts to calculating the limit of $\langle \mathsf{H}_\eta, \mathsf{I}_\eta \rangle$, where the sequence $\mathsf{H}_\eta$ is defined by
$$
\mathsf{H}_\eta(s):=  \left\{\begin{array}{l} |s_1|^{\frac{1}{2}}  |s|_\eta^{-\nu} [ H_\eta (s) - c_\eta m_\eta(s) ],  \hspace{0.6cm} \gamma\in (\frac{d}{2},\frac{d+1}{2}],  \  0< |s_1| \leqslant  1 , \\
 \\
 |s|_\eta^{-1} [ H_\eta (s) - c_\eta m_\eta(s) ] , \hspace{1.25cm} \ \gamma\in (\frac{d+1}{2},\frac{d+4}{2}),  \ 0 < |s_1| \leqslant  1 , \\
 \\
 s_1  H_\eta (s) ,  \hspace{3.88cm}  \gamma\in (\frac{d}{2},\frac{d+4}{2}) , \  |s_1| \geqslant  1 ,
\end{array}\right.
$$
and it is bounded in $L^2(\RR^d)$, uniformly with respect to $\eta$,  thanks to the estimates of Lemma \ref{grandes vitesses}, which implies that $\mathsf{H}_\eta$ converges weakly in $L^2(\RR^d)$,  up to a subsequence.  Let's identify this limit that we denote by $\mathsf{H}_0 \in L^2(\RR^d)$.  We have from Lemma \ref{lemma cv H_eta}, $H_\eta$ converges to $H_0$ in $\mathcal{D}'(\RR^d\setminus\{0\})$. Therefore, thanks to the uniqueness of this limit in $\mathcal{D}'$ and the uniqueness of $H_0$,  the whole sequence $ \mathsf{H}_\eta$ converges  weakly to 
$$
\mathsf{H}_0(s):=  \left\{\begin{array}{l} |s_1|^{\frac{1}{2}} |s|^{-\nu} [H_0 (s)-m(s)] ,  \hspace{0.5cm} \gamma\in (\frac{d}{2},\frac{d+1}{2}],  \ 0< |s_1| \leqslant  1 , \\
 \\
 |s|^{-1} [H_0 (s)-m(s)] ,  \hspace{1.3cm} \gamma\in (\frac{d+1}{2},\frac{d+4}{2}),  \ 0 < |s_1| \leqslant  1 , \\
 \\
 s_1  H_0 (s) ,  \hspace{3.15cm}  \ \gamma\in (\frac{d}{2},\frac{d+4}{2}) , \  |s_1| \geqslant  1 .
\end{array}\right.
$$
For the sequence $\mathsf{I}_\eta$, it is defined by
$$
\mathsf{I}_\eta:=  \left\{\begin{array}{l} |s_1|^{\frac{1}{2}} |s|_\eta^{\nu} m_\eta(s),  \hspace*{0.555cm} \gamma\in (\frac{d}{2},\frac{d+1}{2}], \ 0< |s_1| \leqslant  1 , \\
 \\
 s_1 |s|_\eta m_\eta(s) , \quad \quad \ \ \gamma\in (\frac{d+1}{2},\frac{d+4}{2}), \ 0< |s_1| \leqslant  1 , \\
 \\
m_\eta(s) ,   \hspace{1.98cm}  \gamma\in (\frac{d}{2},\frac{d+4}{2}) , \ |s_1| \geqslant  1 ,
\end{array}\right.  
$$
and converges strongly in $L^2(\RR^{d}\setminus\{0\})$ to 
$$
\mathsf{I}_0:=  \left\{\begin{array}{l} |s_1|^{\frac{1}{2}} |s|^{\nu} m(s),  \hspace{0,53cm} \gamma\in (\frac{d}{2},\frac{d+1}{2}], \ 0< |s_1| \leqslant  1 , \\
 \\
 s_1 |s| \ m(s) , \quad \quad \ \ \gamma\in (\frac{d+1}{2},\frac{d+4}{2}), \ 0< |s_1| \leqslant 1 , \\
 \\
m(s) ,   \hspace{1.93cm}  \gamma\in (\frac{d}{2},\frac{d+4}{2}) , \ |s_1| \geqslant 1 .
\end{array}\right.  
$$
Note that the term $|s|_\eta^\nu$ was introduced to cover the critical case $2\gamma = d+1$, since $s_1 m^2$ is not in $L^1(\{|s_1|\leqslant 1\})$ in this case.  Hence the limit of Lemma \ref{lim kappa_eta} holds true.  
\epl

\begin{remark} In the critical case $2\gamma = d+1$, we can replace the $1$ in the bounds of the integrals of the diffusion coefficient $\kappa$ by any positive value $s_0$ thanks to \eqref{mesure invariante}. In other words, if we change the value of the drift term on a ring, i.e., we define $j^\eps$ by $\int_{\{|v| \leqslant s_0 \eps^{-\frac{1}{3}}\}} v M^2 \, \ud v$, nothing changes in the limit as $\eps \to 0$.  Otherwise, in any case, we will show that for all $\gamma \in (\frac{d}{2},\frac{d+4}{2})$
\begin{equation}\label{kappa reel}
\kappa = - \int_{\{|s_1| > 0\}} s_1 m(s) \Imm H_0(s) \ \ud s .
\end{equation}
\end{remark}
\begin{remark} If $m$ is symmetric then, the function $H_0$ is symmetric with respect to $s_1$ in the following sense:
$$ H_0(s_1, s') = \overline{H_0}(-s_1, s'), \quad \forall s = (s_1, s') \in \RR\times\RR^{d-1} \setminus \{0\}.  $$
This implies that the coefficient $ \kappa$,  the limit of $ \kappa_\eta$ given by Lemma \ref{lim kappa_eta}, is real. Moreover, for all $ \gamma \in \left( \frac{d}{2}, \frac{d+4}{2} \right)$, we have \eqref{kappa reel}.
\end{remark}

\begin{remark} The author cannot directly show from the formula of Lemma \ref{lim kappa_eta} that $\kappa$ is real, i.e., that the imaginary part of the integrals is zero. However, it is only from the limiting equation of fractional diffusion that $\kappa$ must be real.
\end{remark}

\noindent \bpp \ref{vp}. By doing an expansion in $\lambda$ for $B$ and by Proposition \ref{contrainte}, we get
$$
B(\lambda,\eta)=\eta^{-\frac{2}{3}}b(\lambda, \eta)=\eta^{-\frac{2}{3}}b(0,\eta)+ \lambda \int_{\RR^d} M_{0,\eta} M \mathrm{d}v+ O(\lambda^2).
$$ 
Then, for $\lambda=\tilde{\lambda}(\eta)$ and since $B(\tilde{\lambda}(\eta),\eta)=0$, we obtain:
$$ 
\tilde{\lambda}(\eta)= - \eta^{-\frac{2}{3}}b(0,\eta)\bigg(\int_{\RR^d} M_{0,\eta} M \mathrm{d}v\bigg)^{-1}+o\big(\eta^{-\alpha}b(0,\eta)\big) ,
$$ 
which implies that 
$$\eta^{-\alpha}\mu(\eta)=\eta^{\frac{2}{3}-\alpha}\tilde{\lambda}(\eta)= -\eta^{-\alpha}b(0,\eta)\bigg(\int_{\RR^d} M_{0,\eta} M \mathrm{d}v\bigg)^{-1} .
$$
By limit \eqref{lim int  M_eta M}, 
$$\underset{\eta \rightarrow 0}{\lim} \int_{\RR^d} M_{0,\eta}(v)M(v) \ \mathrm{d}v = \int_{\RR^d} M^2(v) \ \ud v , $$ 
and by Lemma \ref{lim kappa_eta},
$$ \underset{\eta \rightarrow 0}{\lim}\ \eta^{-\alpha} \big[ b(0,\eta) +  \ui \eta j^\eps_1 \big] =   \int_{ \{|s_1| > 0\}}  s_1 m(s) \ \mathrm{Im }H_0(s) \ \mathrm{d}s = -\kappa .   $$ 
Hence,  $\underset{\eta \rightarrow 0}{\lim} \eta^{-\alpha}[\mu(\eta) - \ui \eta j^\eps_1] = \kappa$.  So it remains to prove the positivity of $\kappa$.  By integrating the equation of $M_{\eta}:=M_{\tilde{\lambda} (\eta),\eta}$ against $\bar M_\eta$ we obtain:
$$ \int_{\RR^d} \bigg| \nabla_v\bigg(\frac{M_\eta}{M}\bigg) \bigg|^2 M^2 \ \mathrm{d}v + \mathrm{i} \eta \int_{\RR^d} v_1 |M_\eta|^2 \ \mathrm{d}v = \mu(\eta) \int_{\RR^d} |M_\eta|^2 \ \mathrm{d}v . $$
Now, taking the real part and using the equality $\ds \Re \mu(\eta)\|M_\eta\|_{L^2}^2=\kappa \eta^{\alpha}\big(1+O(\eta^\alpha)\big)$ we get:
\begin{equation}\label{kappa >0}
\int_{\RR^d} \bigg| \nabla_v\bigg(\frac{M_\eta}{M}\bigg) \bigg|^2 M^2 \ \mathrm{d}v  = \kappa \eta^{\alpha}\big(1+O(\eta^\alpha)\big) .
\end{equation}
Therefore, multiplying this last equality by $\eta^{-\alpha}$ and performing the change of variable $v=\eta^{-\frac{1}{3}}s$ we obtain:
$$ \int_{\RR^d} \bigg| \nabla_s\bigg(\frac{H_\eta}{m_\eta}\bigg) \bigg|^2 m_\eta^2 \ \mathrm{d}s  = \kappa \big(1+o_\eta(1)\big) . $$
Thus, $ \kappa\geqslant 0$. If $\kappa=0$ then, 
$$ \int_{\RR^d} \bigg| \nabla_s\bigg(\frac{H_0}{m}\bigg) \bigg|^2 m^2 \ \mathrm{d}s \leqslant \liminf  \int_{\RR^d} \bigg| \nabla_s\bigg(\frac{H_\eta}{m_\eta}\bigg) \bigg|^2 m_\eta^2 \ \mathrm{d}s = 0 .$$
Therefore, $H_0= m$. Which leads to a contradiction since $H_0$ is solution to equation \eqref{eq de H_0}.  Hence, the proof of Proposition \ref{vp} is complete.
\epp

\noindent \bpt \ref{main}. The existence of the eigensolution $(\mu(\eta),M_\eta)$ is given by Proposition \ref{contrainte}. The limit of the first point of Theorem \ref{main} follows from inequality \eqref{estimation de N_lambda,eta dans L^2} for $$ |\lambda| = |\tilde{\lambda}(\eta)| = \eta^{-\frac{2}{3}} |\mu(\eta)| \lesssim \eta^{\alpha-\frac{2}{3}} + \eta^{\frac{1}{3}} \underset{\eta \to 0}{\longrightarrow} 0, $$
thanks to \eqref{mu(eta)} and since $\alpha-\frac{2}{3} = \frac{2\gamma-d}{3} > 0$ for $\gamma > \frac{d}{2}$.  Finally, the second point of the theorem is given by Proposition \ref{vp}.

\section{Derivation of the fractional diffusion equation}\label{section diff frac}
The goal of this section is to prove Theorem \ref{main2}. The proof was taken from \cite[Section 4]{DP-2} and adapted here for the case of non-symmetric equilibrium. \\

Let's start by defining the two weighted $L^p$ spaces, $L^p_{F^{1-p}}(\RR^d)$ and $Y^p_F(\RR^{2d})$: 
$$ L^p_{F^{1-p}}(\RR^d) := \left\{f: \mathbb{R}^d\rightarrow \mathbb{R}; \int_{\mathbb{R}^{d}} |f|^p \ F^{1-p}\ \mathrm{d}v< \infty\right\} \  \mbox{ and } \  Y^p_F(\RR^{2d}):=L^p\big(\RR^d; L^p_{F^{1-p}}(\RR^d)\big). $$
Recall that our goal is to show that the solution $f^\eps$ of the Fokker-Planck equation \eqref{fp-theta} converges, weakly star in $L^\infty\left([0, T]; L^2_{F^{-1}}(\mathbb{R}^{2})\right)$, towards $\rho(t,x) F(v)$ when $\eps$ goes to $0$, where $\rho$ is the solution of the following \emph{fractional diffusion} equation:
\begin{equation}
\partial_t\rho +\kappa (-\Delta_x)^{\frac{\beta-d+2}{6}}\rho =0,\quad \rho(0,x)=\int_{\RR^{d}} f_0 \mathrm{d}v .
 \end{equation}
\begin{remark}
Note that we will work with the Fourier transform of $\rho$ and we will prove that 
$ \hat\rho(t,\xi)=\int_{\RR^d} \ue^{-\mathrm{i}x \cdot \xi} \rho(t,x)\mathrm{d}x$ satisfies
\begin{equation}\label{diff}
\partial_t\hat \rho+\kappa|\xi|^{\frac{\beta-d+2}{3}}\hat\rho=0.
\end{equation}
\end{remark}
\subsection{A priori estimates}
We start by recalling the following compactness lemma.
\begin{lemma}$\cite{LebPu,NaPu}$\label{norm}
For initial datum $f_0 \in Y^p_F(\RR^{2d})$ where $p \geqslant 2$ and a positive time $T$.  \smallskip

\ni 1.  The solution $f^\varepsilon$ of \eqref{fp-theta} is bounded in $L^\infty\big([0,T];Y^p_F(\RR^{2d})\big)$ uniformly with respect to $\eps$ since it satisfies 
\begin{equation}\label{estimate1.1}\|f^\varepsilon(T)\|^p_{Y^p_F(\RR^{2d})}+\frac{p(p-1)}{\theta(\varepsilon)}\int_0^T\int_{\mathbb{R}^{2d}}\bigg|\nabla_v\bigg(\frac{f^\varepsilon}{F}\bigg)\bigg|^2 \bigg|\frac{f^\varepsilon}{F}\bigg|^{p-2}\ F\ \mathrm{d}v \mathrm{d}x \mathrm{d}t \leqslant \|f_0\|^p_{Y^p_F(\RR^{2d})}.
\end{equation}
2.  The density $\rho^\varepsilon(t,x)=\int_{\mathbb{R}^d}f^\varepsilon \ \mathrm{d}v$ is such that
\begin{equation} 
\|\rho^\varepsilon(t)\|^p_{L^p(\RR^d)} \leqslant \|F\|_{L^1(\RR^d)}^{p-1} \|f_0\|^p_{Y^p_F(\RR^{2d})}\quad \mbox{ for all } \quad t\in [0,T].
\end{equation}
3.  Up to a subsequence, the density $\rho^\varepsilon$ converges weakly star in $L^\infty \big([0,T]; L^p(\mathbb{R}^d)\big)$ to $\rho$.  \smallskip

\ni 4.  Up to a subsequence, the sequence $f^\varepsilon$ converges weakly star in $L^\infty \big([0,T]; Y^p_F(\RR^{2d})\big)$ to the function $f=\rho(t,x)F(v)$.
\end{lemma}
\begin{remark}
The proof of this lemma remains unchanged in the presence of the drift term, as the drift term only appears in the position variable, $f(\eps^{-\alpha}t,\eps^{-1}x+\eps^{-\alpha}j^\eps t,v)$, which is integrated over $\RR^d$.  See \cite[Lemma 6.1]{NaPu} for the proof without the drift term.
\end{remark}

As a consequence, we have the following estimate.
\begin{corollary}[\cite{LebPu,DP-2}] Let $F = M^2$ with $M$ satisfying \ref{a1} and $\beta=2\gamma\in (d,d+4)$.
Let $f^\eps$ solution to \eqref{fp-theta} with $\theta(\varepsilon)=\varepsilon^{\frac{\beta-d+2}{3}}$.  Assume that $\|f_0/F\|_{\infty} \leqslant C$. Then,  $g^\eps=f^\eps F^{-\frac{1}{2}}$ satisfies the following estimate
\begin{equation}\label{estblan}
\int_0^T \int_{\RR^d}\left(\int_{\RR^d} \left| g^\eps-\rho^\eps F^{\frac{1}{2}}\right|^2 \mathrm{d}v \right)^{\frac{\beta-d+2}{\beta-d}} \ud x \ud t \leqslant C \eps^\frac{\beta-d+2}{3}.
\end{equation}
\end{corollary}

The proof of the previous estimate is based on inequality \eqref{estimate1.1} and a Nash type inequality. See \cite[Corollary 4.3]{DP-2} for details. Note that Assumption \ref{a1} is crucial, see \cite{CGGR} and the references therein.

\subsection{Weak limit and proof of Theorem \ref{main2}}
This subsection is an adaptation of \cite[Subsection 4.2]{DP-2} to the non-symmetric case. By solving equation (\ref{rescaled}), we write
$$\hat g^\eps(t,\xi,v)= \ue^{-t\theta(\eps)\mathcal{L}_\eta}\hat g(0,\xi,v) ,
$$
which gives going back to the rescaled space variable $x$
$$
 g^\eps(t,x,v)=\frac{1}{(2\pi)^d} \int_{\RR^d} \ue^{\mathrm{i} x \cdot \xi} \ \hat g^\eps(t,\xi,v) \ \mathrm{d}\xi  .
$$
Our purpose is to pass to the limit when $\eps\rightarrow 0$.  Recall that $f^\eps(t,x,v) \geqslant 0$ and we have for all $t\geqslant 0$,  $\int f^\eps(t,x,v)  \ud x \ud v=\int f_0(x,v)\ud x \ud v$.
Let $\hat \rho^\eps (t,\xi)= \int \ue^{-\mathrm{i}x \cdot \xi}\rho^\eps (t,x) \ud x $ be the Fourier transform in $x$ of  $\rho^\eps = \int f^\eps \mathrm{d}v= \int g^\eps F^{\frac{1}{2}} \mathrm{d}v$.

\begin{proposition}\label{lem:equilibrium}
For  all $\xi\in \RR^d$,  $\hat\rho^\varepsilon(\cdot,\xi)$ converges  to $\hat \rho(\cdot,\xi)$, 
unique solution  to the ODE
\begin{equation}\label{fracdifeq}
\partial_t\hat\rho+\kappa|\xi|^\alpha \hat \rho=0, \quad \hat \rho_0 = \int_{\RR^d}\hat  f_0 \ \mathrm{d}v   .
\end{equation}
\end{proposition}

\noindent \bp Let $\xi \in \RR^d $ and let
$M_\eta$ be the unique solution in $L^2(\RR^d,\CC)$ of $\ \mathcal{L}_\eta(M_\eta)=\mu(\eta)M_\eta \ $ given in Theorem \ref{main}. One has  
$$
\begin{array}{rcl}
\dps \frac{\mathrm{d}}{\mathrm{d}t} \int_{\RR^d}  \hat g^\eps(t,\xi,v) M_\eta(v) \ \mathrm{d}v&=&\dps  \int_{\RR^d} \partial_t\hat g^\eps  M_\eta \ \mathrm{d}v  \\
&=& \ds -\eps^{-\alpha}\int_{\RR^d} \left[\mathcal{L}_\eps(\hat g^\eps) - \ui \eta j^\eps_1 \hat g^\eps \right] M_\eta \ \mathrm{d}v\\
&=&
\dps -\eps^{-\alpha}\int_{\RR^d} \hat g^\eps   \left[ \mathcal{L}_\eps(M_\eta)- \ui \eta j^\eps_1 M_\eta \right] \mathrm{d}v \\
&=& \ds -\eps^{-\alpha}\left[\mu(\eta)- \ui \eta j^\eps_1 \right] \int_{\RR^d} \hat g^\eps (t,\xi,v) M_\eta(v) \ \mathrm{d}v  .  
\end{array}
$$
Therefore one has, with $\ds F^\eps (t,x) :=  \int_{\RR^d}   g^\eps(t,x,v) M_\eta \ \mathrm{d}v$, 
\begin{equation}\label{gl9}
 \hat F^\eps (t,\xi) = \ue^{-t\eps^{-\alpha} \left[\mu(\eta)- \ui \eta j^\eps_1 \right]}  \hat F^\eps (0,\xi) , \qquad \forall t \geqslant 0.
\end{equation} 
By Theorem \ref{main}, we have
$ \eps^{-\alpha} \left[\mu(\eps |\xi|)- \ui \eps |\xi| j^\eps_1 \right] \rightarrow \kappa  |\xi|^\alpha$. 
Moreover, the following limit holds true:
\begin{equation}\label{gl11}
\forall \xi \in \RR^d, \qquad \hat F^\eps (0,\xi) =  \int_{\RR^d} \hat g^\eps(0,\xi,v)  M_\eta \ \ud v  \longrightarrow  \hat \rho_0\left(\xi\right) . 
\end{equation}
The verification of \eqref{gl11} is easy. One has $ \hat g^\eps(0,v,\xi)=  \hat f_0(v,\xi) F^{-\frac{1}{2}}(v) = \hat f_0(v,\xi)M^{-1}(v)$ and $M_\eta \to M$ in $L^2(\RR^d)$ thanks to the first point of Theorem \ref{main}. Thus, \eqref{gl11} holds true by Cauchy-Schwarz inequality by writing: 
$$ \bigg|  \int_{\RR^d} \hat g^\eps(0,\xi,v)  M_\eta \ \mathrm{d}v - \hat \rho_0(\xi) \bigg| \leqslant \bigg(\int_{\RR^d} \frac{f_0^2}{F} \ \mathrm{d}v \bigg)^{\frac{1}{2}}\bigg(\int_{\RR^d} |M_\eta-M|^2\ \mathrm{d}v \bigg)^{\frac{1}{2}} .$$
It remains to verify 
\begin{equation}\label{gl10}
\forall \xi \in {\RR^d}, \qquad  \int_{\RR^d} \hat g^\eps(t,\xi,v)  M_\eta \ \mathrm{d}v  \longrightarrow  \hat \rho(t,\xi)  
\quad \text{ in } \ \mathcal D'\big(]0,\infty[\times\RR^{d}\big).
\end{equation}
By \eqref{gl9} and \eqref{gl11}, for all $\xi\in \RR^d $ and $t\geqslant 0$, one has
$\underset{\eps\rightarrow 0}{\lim}\ \hat F^\eps (t,\xi)= \ue^{-t\kappa |\xi|^\alpha}\hat\rho_0(\xi)$, thus \eqref{gl10}
will be consequence of the weaker
\begin{equation}\label{gl10bis}
 \int_{\RR^d}  g^\eps(t,x,v)  M_\eta \ \mathrm{d}v  \longrightarrow   \rho(t,x)  
\quad \text{ in } \ \mathcal D'\big(]0,\infty[\times\RR^d\big) .
\end{equation}

\noindent Let us now verify \eqref{gl10bis}. For that purpose, we write
$$
\int_{\RR^d}  g^\eps  M_\eta \ \mathrm{d}v  -\rho= 
\int_{\RR^d} \left( g^\eps- \rho^\eps F^{\frac{1}{2}} \right) M_\eta \ \mathrm{d}v +  \rho^\eps \int_{\RR^d} (M_\eta-F^{\frac{1}{2}})F^{\frac{1}{2}} \ \mathrm{d}v 
+\rho^\eps -\rho  .
$$
By using (\ref{estblan}) and the limit of the first point of Theorem \ref{main},  we pass to the limit.
The proof of Proposition \ref{lem:equilibrium} is complete.
\ep

\noindent \bpt \ref{main2}.  From the two last items in Lemma \ref{norm}, we have just to prove that 
for any given $\xi$, the Fourier transform $\hat\rho(t,\xi)$ of the weak limit
$\rho(t,y)$, is solution of equation \eqref{diff}, which is precisely Proposition \ref{lem:equilibrium}.  \ept

\begin{remark}
The term $j^\eps$ that appears in the development of $\mu(\eta)$ corresponds to a drift in the kinetic equation \eqref{fp-theta} for $\beta \geqslant d+1$. Indeed, in the case $\beta > d+1$,  we have 
$$ \ui \eta j^\eps_1 = \ui  \eps |\xi| \int_{\RR^d} v \cdot \frac{\xi}{|\xi|} M^2(v) \ \ud v = \ui \eps \left(\int_{\RR^d} v F(v) \ \ud v \right) \cdot \xi  .  $$
In the case $\beta = d+1$, this corresponds to a truncation of the previous integral since it is not finite at $+\infty$. Moreover, in this case, we have $\eta \int_{\{|v| \leqslant \eps^{-\frac{1}{3}}\}} v_1 M^2(v) \ud v \underset{\eps \to 0}{\sim} \eps \ln\big(\eps^{-\frac{1}{3}}\big) \mathfrak{j}_m \cdot \xi$,  which can be seen as a kind of correction of order $|\ln(\eps)|$.  More precisely, we have the following
\end{remark}
\begin{lemma}\label{j_eps equiv ln eps}  One has the following limit:
$$
 \underset{\eps \to 0}{\lim} \  \frac{3}{\left| \ln (\eps) \right|} \int_{\{|v| \leqslant \eps^{-\frac{1}{3}}\}} v M^2(v) \ \ud v  =  \underset{R \to \infty}{\lim} \ \frac{1}{\ln R} \int_{ \{|v| \leqslant R\}} v  M^2(v)\  \ud v = \int_{\mathbb{S}^{d-1}} s \ m^2(s) \ \ud \sigma (s) =: \mathfrak{j}_m ,  $$
where $\ud \sigma (s)$ denotes the measure on $\mathbb{S}^{d-1}$.
\end{lemma}
\bp The limit of this lemma is obtained by L'H\^opital's rule. Indeed, let $\psi$ be the function defined by
$$ \psi(\eps) := \int_{\{|v| \leqslant \eps^{-\frac{1}{3}}\}} v M^2(v) \, \ud v  .$$
On one hand, by L'H\^opital's rule, we have
$$ \underset{\eps \to 0}{\lim} \  \frac{\psi (\eps)}{\ln(\eps^{-\frac{1}{3}})}  = \underset{\eps \to 0}{\lim} \  \frac{\psi' (\eps)}{-\frac{1}{3}\eps^{-1}} = \underset{\eps \to 0}{\lim} \  -3 \eps \psi' (\eps) .$$
On the other hand, by performing the change of variables $v = \eps^{-\frac{1}{3}} s$, we write
$$
 \psi' (\eps) = -\frac{1}{3}  \eps^{-\frac{4}{3}} \int_{|v| = \eps^{-\frac{1}{3}}} v  M^2(v) \ \ud \sigma (v)  = -\frac{1}{3}  \eps^{-1} \int_{|s| = 1} s \ m_{\eps}^2(s) \ \ud \sigma (s) .
$$
Hence, 
$$ \underset{\eps \to 0}{\lim} \  \frac{\psi (\eps)}{\ln(\eps^{-\frac{1}{3}})}  = \underset{\eps \to 0}{\lim} \ \int_{|s| = 1} s \ m_{\eps}^2(s) \ \ud \sigma (s) = \int_{|s| = 1} s \ m^2(s) \ \ud \sigma (s) .$$
\ep

\section{Comments and remarks}\label{section comments}
We will conclude this paper with some comments on the assumptions, along with examples of equilibria and potentials, and a discussion on the case of a potential with a scaling greater than that of the Laplacian.

\subsection{Comments on assumptions and examples}\label{subsection comments}
In summary, the most important assumptions are \ref{a1} and \ref{a3}. Assumption \ref{a2} is technical. Regarding Assumption \ref{a5}, it is elaborated upon in the following subsection. \smallskip

\ni \textbf{1.} Assumption \ref{a1} is weaker and more general than the existence of a limit for $|v|^\gamma M(v)$ as $|v|$ tends to infinity, because with  \ref{a1}, the function $|v|^\gamma M(v)$ may not have a limit at infinity, and its limit may depend on the directions in the case where $M$ is not radial.  For $M$ given by
$$ M(v) := C_\gamma (1+2v_1^2+ |v'|^2)\langle v\rangle^{-\gamma-2} ,$$
where $C_\gamma$ is a normalization constant such that $\|M\|_{L^2(\RR^d)} = 1$, all the assumptions are satisfied and the limit of $|v|^\gamma M(v)$, when $|v| \to \infty$, does not exists. \\
\ni \textbf{2.} Within a ball of radius $R$, as large as desired, the equilibrium $M$ can be any positive bounded function. In particular, $M$ may exhibit oscillations, such as
$$ M(v) := C_\gamma \left(2+ \frac{\cos\langle v\rangle}{\langle v\rangle^\sigma}\right)\langle v\rangle^{-\gamma} .$$
With this example, all assumptions are satisfied for $\sigma \geqslant 2$, and the potential $W$ satisfies
$$ W(v) \underset{|v| \to \infty}{\sim} \frac{2\gamma(\gamma-d+2)\langle v\rangle^{-2}}{2+\langle v\rangle^{-\sigma} \cos\langle v\rangle} - \frac{\langle v\rangle^{-\sigma} \cos\langle v\rangle}{2+\langle v\rangle^{-\sigma} \cos\langle v\rangle} \cdot $$
For $\sigma \geqslant 2$, we get $W(v) \sim \gamma(\gamma-d+2)\langle v\rangle^{-2}$ while for $\sigma < 2$, we obtain $ W(v) \sim - \frac{\langle v\rangle^{-\sigma} \cos\langle v\rangle}{2+\langle v\rangle^{-\sigma} \cos\langle v\rangle} \cdot$ In the latter case, we also have $\langle v\rangle^{\sigma} W \in L^\infty(\mathbb{R}^d)$. \\
\ni \label{comment sym m}\textbf{3.} The diffusion coefficient $\kappa$ is real. Indeed,  for a real initial condition $\rho_0$, the fractional diffusion equation
$$ \pa_{t} \rho + \kappa (-\Delta_x)^{\frac{\alpha}{2}} \rho = 0 $$
admits a real solution $\rho$ if and only if the coefficient $\kappa$ is real. This is the case since $\rho$ is the limit of the real sequence $\rho^\eps$. This result is proved in Appendix \ref{sec equa diff frac}.  \\
\ni \textbf{4.} Assuming only \ref{a1} and \ref{a2} for the existence of the eigen-pair $(\mu(\eta),M_\eta)$ allows us to cover a wide class of potentials $W$.  For examples: 
    \begin{itemize}
    \item For any $M$ satisfying \ref{a1} such that $\langle v \rangle^{-2} W \in L^\infty(\RR^d)$, Assumption \ref{a2} is satisfied. Indeed, since $\Delta_v M = W M$ then,  by integrations by parts we write
    \begin{align*}
    \int_{\RR^d} \langle v \rangle^2 \bigg|\na_v\bigg(\frac{M}{\langle v \rangle^2}\bigg)\bigg|^2 \ud v &= - \int_{\RR^d} \frac{M}{\langle v \rangle^2} \na_v \cdot \bigg( \langle v\rangle^2  \na_v\bigg(\frac{M}{\langle v \rangle^2}\bigg)\bigg) \ud v \\
    &= - \int_{\RR^d}  \frac{M}{\langle v \rangle^2} \Delta_v M\ \ud v + \int_{\RR^d} \frac{M}{\langle v \rangle^2}  \na_v \cdot \bigg( \frac{2v}{\langle v\rangle^2} M\bigg) \ud v  \\
      &= - \int_{\RR^d} \frac{M}{\langle v \rangle^2} W M\ \ud v - \int_{\RR^d} 2v \frac{M}{\langle v \rangle^2} \cdot  \na_v \bigg( \frac{M}{\langle v\rangle^2} \bigg) \ud v .
    \end{align*}
    Hence, 
    $$ \int_{\RR^d} \langle v \rangle^2 \bigg|\na_v\bigg(\frac{M}{\langle v \rangle^2}\bigg)\bigg|^2 \ud v = - \int_{\RR^d}  \frac{W}{\langle v \rangle^2} M^2 \ud v + d \int_{\RR^d} \frac{M^2}{\langle v \rangle^4} \ \ud v \lesssim 1 .$$
    \item For $M(v) = C_\gamma \left(2+\frac{\sin(\langle v \rangle^3)}{\langle v \rangle}\right) \langle v \rangle^{-\gamma} $, Assumptions \ref{a1} and \ref{a2} are satisfied and one has
    $$ W(v) \underset{|v| \to \infty}{\sim} -\frac{9}{2} |v|^3 \sin(\langle v \rangle^3) . $$ 
\end{itemize}     
\ni \textbf{5.} Theorem \ref{main} is valid for $\beta \geqslant d+4$, and there is no restriction on $\beta$ in the proof of the existence of the pair $(\mu(\eta),M_\eta)$; it is valid for any $\beta > d$. Similarly, for Theorem \ref{main2} on the diffusion limit, some additional estimates on $M_{0,\eta}$ of the type given in Lemmas \ref{petites vitesses} and \ref{grandes vitesses} are required. Indeed, in the case $\beta \geqslant d+4$ we have
$$ -\eta^{-2} \left[ b(0,\eta) +  \ui \eta \int_{\RR^d} v_1 M^2 \ud v \right] = \kappa_\eta + r_\eta ,$$ 
where 
$$ \kappa_\eta := \ui \eta^{-1} \int_{|v_1| \leqslant \eta^{-\frac{1}{3}}} v_1 \big[M_{0,\eta}-M\big]M \ud v \ \mbox{ and } \ r_\eta := \ui  \eta^{\frac{\beta-d-4}{3}} \int_{|s_1| \geqslant 1} s_1 \big[H_\eta(s)-m_\eta(s)\big] m_\eta(s) \ud s .$$
$\bullet$ For $\beta > d+4$, the remainder $r_\eta$ tends to $0$, since it can be shown that
$$ |r_\eta| \leqslant \eta^{\frac{\beta-d-4}{3}} \left\|m_\eta\right\|_{L^2(\{|s_1| \geqslant 1\})} \left(\left\|s_1 H_\eta\right\|_{L^2(\{|s_1| \geqslant 1\})} + \left\|s_1 m_\eta\right\|_{L^2(\{|s_1| \geqslant 1\})} \right) \lesssim \eta^{\frac{\beta-d-4}{3}} .$$
It is in the passage to the limit in $\kappa_\eta$ that the estimates on $M_{0,\eta}$ are needed, and the limit $\kappa$ in this case will be given by an integral equivalent to the fourth moment of $M^2$. \\
$\bullet$ The critical case $\beta = d+4$ requires the introduction of the term $|\ln(\eps)|$ in $\theta(\eps)$ and $\mu(\eta)$. The diffusion coefficient $\kappa$ will be given by a principal value.  See for instence \cite{CNP,FT,BM}. \medskip

\ni \textbf{Connection with the work of probabilists.  }We conclude this subsection with a comparison to the results of the probabilists in \cite{FT-d1} and \cite{FT}. It is important to recall that these results concern only the diffusion limit. Our focus is on the case $d \geqslant 2$, as studied in \cite{FT}, since for the case $d=1$, the authors assume symmetry of the equilibrium along with additional hypotheses on the potential, and the result of our paper remains more general in this case. Let us begin by clarifying the notations used. The Fokker-Planck operator studied in \cite{FT} is given by
$$
\frac{1}{2} \left(\Delta_v f + \beta \mathrm{div}_v\big[ \mathcal{F}(v) f \big] \right),
$$
which corresponds to $\frac{1}{2} \mathsf{Q}$ defined in \eqref{defQ}, with $\mathcal{F}(v) = -\frac{\nabla_v F}{\beta F}$, where $F$ is the equilibrium in that paper. In \cite{FT}, the authors assume a spherical decomposition of the equilibrium.  More precisely,  there exists a potential $U : \mathbb{R}^d \setminus \{0\} \longrightarrow (0,\infty)$ of the form
$$
U(v) = \Gamma(|v|)\gamma\left(\frac{v}{|v|}\right),
$$
where $\gamma : \mathbb{S}^{d-1} \longrightarrow (0,\infty)$ and $\Gamma : \mathbb{R}^+ \longrightarrow (0,\infty)$ are of class $C^\infty$, satisfying $\Gamma(r) \sim r$ as $r \to \infty$ and such that,  for any $v \in \mathbb{R}^d \setminus \{0\}$, one has
$$
\mathcal{F}(v) = \nabla[\log U(v)] = [U(v)]^{-1} \nabla U(v).
$$
With these assumptions, the field $\mathcal{F}$ is of class $C^\infty$ on $\mathbb{R}^d \setminus \{0\}$, and the equilibrium $F$, in these notations, is given by $F(v) = [U(v)]^{-\beta}$. Thus, all the assumptions \ref{a1}--\ref{a5} are satisfied. Finally,  for $s \in \mathbb{R}^d \setminus \{0\}$, by taking $\lambda = |s|^{-1}$ in \eqref{m sur la sphere}, the limit of the rescaled equilibrium $m$ satisfies $m(s) = |s|^{-\gamma} m(s/|s|)$.
 
\subsection{The case of a potential stronger than the Laplacian}\label{subsection W}
In this subsection, we discuss the case where the Laplacian and the potential have different regimes, specifically the case where the potential dominates, since the case where the potential has the same scaling or is weaker is discussed throughout the paper.  This case is precisely  Assumption \ref{a5}.

An example where all the terms of the Fokker-Planck operator, taking into account the advection term, are of the same order is given by the classical distribution considered in \cite{LebPu, BM, DP-2} and \cite{NaPu}, and which corresponds to $\Gamma(r)=\sqrt{1+r^2}$ and $\gamma \equiv 1$ in \cite{FT-d1,FT}, is given by
$$ M(v) := C_\gamma \langle v \rangle^{-\gamma}, $$
where $ C_\gamma $ is a normalization constant such that $ \|M\|_{L^2} = 1 $. With this equilibrium, the potential $ W$ is given by
 $$ W(v) = \frac{\gamma(\gamma-d+2)|v|^2-\gamma d}{\langle v \rangle^4} \underset{|v| \to \infty}{\sim} \frac{\gamma(\gamma-d+2)}{|v|^2} .$$
For $\gamma = d-2$ with $d \in \{5,6,7\}$, we get $ W(v) = -\frac{d(d-2)}{\langle v \rangle^4}$, which disappears in the limit after rescaling. The limit of the rescaled equilibrium is  $m(s) = C_\gamma |s|^{-\gamma}$, which is a solution of the equation $\ \Delta_s m = 0$ on $\RR^d\setminus \{0\}$.\\ 

In the case where $ |W(v)| \leqslant C \langle v \rangle^{-\sigma} $ with $ \sigma \in (-1,2) $, the Laplacian is negligible compared to the potential, and the appropriate scaling is given by $ \eta^{-\frac{1}{\sigma+1}} $ instead of $\eta^{-\frac{1}{3}} $. Thus, the rescaled solution $ H_\eta(s) := \eta^{-\frac{\gamma}{\sigma+1}}M_{0,\eta}\big(\eta^{-\frac{1}{\sigma+1}}s\big) $ satisfies the equation
$$ \left[-\eta^{\frac{2-\sigma}{\sigma+1}} \Delta_s + W_\eta(s) + \ui s_1 \right] H_\eta = -\eta^{-\frac{\gamma+\sigma}{\sigma+1}} b(0,\eta) \Phi_\eta, $$
where $ \Phi_\eta(s) := \Phi\big(\eta^{-\frac{1}{\sigma+1}}s\big) $ and $W_\eta$ is the rescaled potential given by $$ W_\eta(s) := \eta^{-\frac{\sigma}{\sigma+1}} W\big(\eta^{-\frac{1}{\sigma+1}}s\big) ,  $$
satisfying the uniform estimate $|W_\eta(s)| \lesssim |s|^{-\sigma}$. \\

One can prove that the estimates of Lemma \ref{grandes vitesses} are replaced here by
\begin{itemize}
\item For all $\gamma \in (\frac{d}{2},\frac{d+1}{2}] $ and $ \nu \in (\gamma-\frac{d+1}{2},\gamma-\frac{d-1}{2}) $, we have
$$ \eta^{\frac{2-\sigma}{\sigma+1}} \left\| |s_1|^{\frac{1}{2}} |s|_\eta^{-\nu} (H_\eta-c_\eta m_\eta) \right\|_{L^2(\{|s_1| \leqslant 1\})}  \lesssim 1 .$$
\item For all $ \gamma \in (\frac{d}{2},\frac{d+4}{2}) $, we have
$$ \eta^{\frac{2-\sigma}{\sigma+1}} \left\|\frac{H_\eta-c_\eta m_\eta}{|s|_\eta} \right\|_{L^2(\mathbb{\RR}^d)} + \eta^{\frac{2-\sigma}{\sigma+1}} \bigg\|\na_s\bigg(\frac{H_\eta}{m_\eta}\bigg)m_\eta\bigg\|_{L^2(\mathbb{\RR}^d)} +  \left\|s_1 H_\eta \right\|_{L^2(\{|s| \geqslant 1\})} \lesssim 1.$$
\end{itemize}
The previous estimates with Assumption \ref{a3}, which implies that $ m_\eta \to m $ in $L^2(\{|s| \geqslant r\}) $ for any $ r>0 $, and with the uniform bound for $ W_\eta m_\eta $ in $ L^2(\{|s| \geqslant r\})$, allow us to show that $ W_\eta \to 0 $ as $\eta \to 0 $, at least in the sense of distributions.\medskip

\ni\textbf{Conclusion:} In the case where the Laplacian is negligible compared to the potential, if we impose the condition of the convergence of the rescaled equilibrium $m_\eta$, then we do not capture anything at the limit, i.e., we obtain a trivial limit for $ H_\eta $.

\appendix

\section{On the solution of the fractional diffusion equation}\label{sec equa diff frac}
The solution of the fractional diffusion equation
$$ \pa_t \rho + \kappa (-\Delta_x)^\alpha \rho = 0 $$ 
is given by
$$ \rho(t,x) = e^{-\kappa t (-\Delta_x)^\alpha} \rho(0,x),$$
where $\rho(0,\cdot)$ is a real function.  For $\kappa = \kappa_R + i \kappa_I$, we get
$$
e^{-\kappa t (-\Delta)^\alpha} = e^{-\kappa_R t (-\Delta)^\alpha} \cdot e^{-i \kappa_I t (-\Delta)^\alpha}.
$$
The part $e^{-i \kappa_I t (-\Delta)^\alpha}$ can be decomposed as
$$
e^{-i \kappa_I t (-\Delta)^\alpha} = \cos(\kappa_I t (-\Delta)^\alpha) - i \sin(\kappa_I t (-\Delta)^\alpha).
$$
For $\rho(t,x)$ to be real, the imaginary part introduced by $-i \sin(\kappa_I t (-\Delta)^\alpha)$ must vanish. This is true if and only if $\kappa_I = 0$, because otherwise, $\sin(\kappa_I t (-\Delta)^\alpha)$ can be non-zero at certain times, making the solution complex. 

\subsection*{Definition of $\sin(\kappa_I t (-\Delta)^\alpha)$}
The function $\sin(\kappa_I t (-\Delta)^\alpha)$ is defined using operator theory and analytic functions of operators. Here's how we can understand and define this expression:  \smallskip

\ni $\bullet$ 
For an operator $A$, we define
$$
\sin(A) = \frac{e^{iA} - e^{-iA}}{2i}.
$$
Applied to our case
$$
\sin(\kappa_I t (-\Delta)^\alpha) = \frac{e^{i \kappa_I t (-\Delta)^\alpha} - e^{-i \kappa_I t (-\Delta)^\alpha}}{2i}.
$$
\ni $\bullet$  The sine function can also be defined by its Taylor series
$$
\sin(A) = \sum_{n=0}^\infty \frac{(-1)^n}{(2n+1)!} A^{2n+1}.
$$
Therefore,
$$
\sin(\kappa_I t (-\Delta)^\alpha) = \sum_{n=0}^\infty \frac{(-1)^n}{(2n+1)!} (\kappa_I t (-\Delta)^\alpha)^{2n+1}.
$$ 
~~~The linear operator $\sin(\kappa_I t (-\Delta)^\alpha)$ is well-defined on test functions and, generally, on an appropriate functional space.  

For real values of $\kappa_I$, the term $\kappa_I t$ causes periodic oscillations, which can introduce imaginary components into the solution if $\kappa_I \neq 0$. 

\ni \textbf{Conclusion:} In summary, the definition of $\sin(\kappa_I t (-\Delta)^\alpha)$ is in terms of operator exponentials or Taylor series, and the condition for the solution to remain real is that $\kappa_I = 0$.

{\footnotesize
\bibliography{biblio-FP}
}
\bibliographystyle{abbrv}

\end{document}